\documentclass[12pt]{article}
\usepackage[utf8]{inputenc}
\usepackage[T1]{fontenc}
\usepackage[includeheadfoot,top=1 cm, bottom=1.5 cm, left=2 cm, right=2 cm]{geometry}
\usepackage[english]{babel}
\usepackage{lmodern}
\usepackage{graphicx}
\usepackage{subcaption}
\usepackage{wrapfig}
\usepackage{colortbl} 
\usepackage[svgnames,table]{xcolor}
\usepackage{array}
\usepackage{amsfonts}
\usepackage{amssymb}
\usepackage{amsmath}
\usepackage{setspace}
\usepackage{hyperref}
\usepackage{multirow}
\usepackage{fancybox}
\usepackage{upgreek}
\usepackage{comment}
\usepackage{xcolor}
\usepackage{pifont}
\usepackage{makecell}
\usepackage[overload]{empheq}
\usepackage{cite}

\usepackage{algorithm}
\usepackage{algorithmicx}
\usepackage{algcompatible}
\usepackage{amsthm}
\usepackage{bm}
\usepackage{color}
\usepackage{soul}
\soulregister\cite7
\soulregister\eqref7
\soulregister\ref7

\theoremstyle{remark}
\newtheorem{remark}{Remark}

\newcommand{\bdelta}{\boldsymbol{\delta}}
\newcommand{\blambda}{\boldsymbol{\lambda}}

\newcommand{\bg}{\boldsymbol{g}}

\newcommand{\bB}{\boldsymbol{B}}
\newcommand{\bK}{\boldsymbol{K}}
\newcommand{\bU}{\boldsymbol{U}}
\newcommand{\bF}{\boldsymbol{F}}
\newcommand{\bD}{\boldsymbol{D}}

\newcommand{\bzero}{\boldsymbol{0}}

\newcommand{\cmark}{\textcolor{green!60!black}{\ding{51}}} 
\newcommand{\xmark}{\textcolor{red}{\ding{55}}}           
\definecolor{mygray}{gray}{0.85}
\arrayrulecolor{black}

\title{\bf Parameter-unbounded convergence of Crossed-Secant accelerated Uzawa and penalty-splitted algorithms for frictionless contact}
\author{{KOLIESNIKOVA Daria}\textsuperscript{1}, {RAMIÈRE Isabelle}\textsuperscript{1} \newline
	\\\small{\textsuperscript{1}CEA, DES, IRESNE, DEC, SESC, F-13108 Saint-Paul-lez-Durance, France}}

\begin{document}
\maketitle

\begin{abstract}
  
We propose a unified iterative framework for the solution of frictionless mechanical contact problems, which relies exclusively on the solution of standard stiffness systems. The framework is built upon a two-step fixed-point algorithm: first, the displacement (primal variable) is computed for given contact forces; second, the contact forces (dual variable) are updated based on the displacement solution.
		
The choice of the dual update scheme depends on the numerical contact formulation under consideration. Specifically, the Uzawa iterative scheme is obtained for the Lagrange multiplier formulation, while a penalty-based operator-splitting strategy is proposed for the penalty contact formulation. The main interest of such displacement-force splitting strategy is to involve only standard rigidity matrices in the solving step: no saddle-point or penalized ill-conditionned coefficient matrices have to be handled, so no specialized preconditioning is required. Moreover only the right-hand side of the system is updated throughout the iterations, which enables matrix factorization reuse or efficient iterative solvers initialization.

The main limitation of such splitting iterative strategies lies in the inherently slow convergence of the underlying fixed-point iterations. Moreover, convergence is guaranteed only within a narrow range of numerical parameter values (i.e., the augmentation parameter for the Uzawa method and the penalty parameter for the penalty-splitted formulation). This work addresses both issues by applying the Crossed-Secant fixed-point acceleration strategy, which substantially improves the convergence rate and renders the iterative schemes effectively parameter-unconstrained.
                
To the best of our knowledge, this contribution provides the first computational demonstration of efficient, parameter-unbounded convergence for such contact formulations. The substantial practical benefits of the proposed approach are illustrated through representative three-dimensional academic and industrial frictionless contact problems. Moreover, our findings pave the way for large-scale, parallel simulations of multi-contact systems.

\end{abstract}

              


\section{Introduction} \label{s:intro}
	
Many physical and engineering problems can be formulated as variational or hemivariational inequalities~\cite{Lions1967,Panagiotopoulos1993}, where inequality constraints arise either from admissible sets or from non-smooth constitutive laws.
Such formulations arise naturally in diverse fields including contact mechanics, phase transformations, plasticity, fluid-structure interaction, optimization, finance and economics,  etc.~\cite{Kinderlehrer2000,Sofonea2018,Migorski2013,Dafermos1992,Karatzas2000,ODEN1980,Han2019,Nocedal2006}.

\begin{figure}[!ht]
	\centering
	
	\begin{subfigure}[b]{0.48\textwidth}
		\centering
		\includegraphics[
		page=3,
		width=5.8cm,
		clip
		]{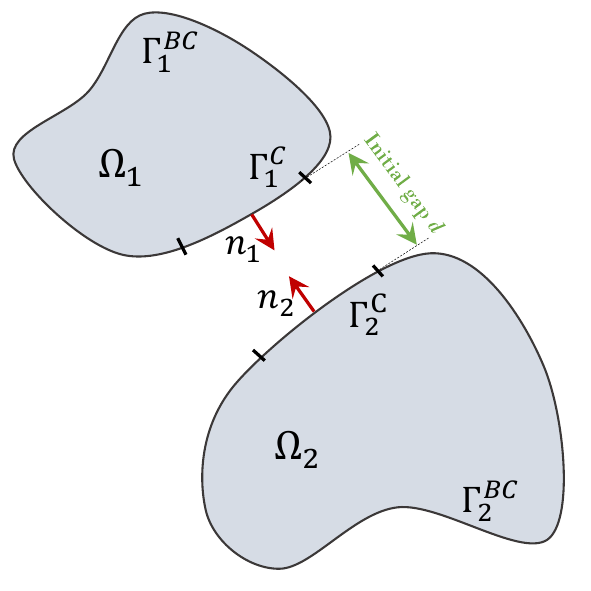}
		\caption{Schematic representation of a contact problem between two solids (with the notation $\Gamma^{BC}= \Gamma^D \cup \Gamma^N$)}
		\label{fig:contact_pb}
	\end{subfigure}
	\hfill
	\begin{subfigure}[b]{0.48\textwidth}
		\centering
		\includegraphics[
		page=4,
		width=7cm,
		trim=0.5cm 1cm 0.1cm 1.2cm,
		clip
		]{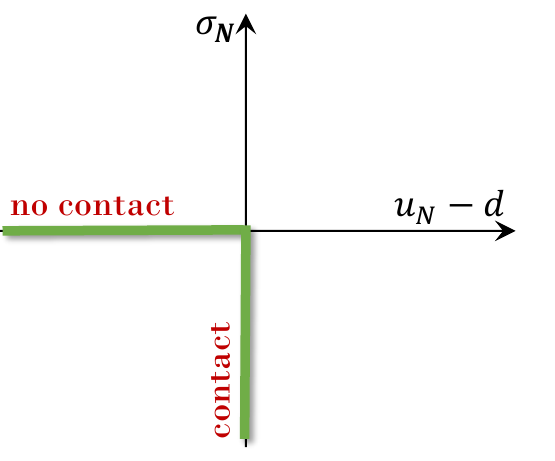}
		\caption{Signorini contact law: dual variable versus normal gap}
		\label{fig:Signorini}
	\end{subfigure}
	
	\caption{Contact problem and Signorini contact law}
	\label{fig:Signorini_contact}
      \end{figure}

In this work, we focus on a particularly important class of such problems, namely contact problems in solid mechanics, and more specifically on problems governed by frictionless Hertz-Signorini-Moreau contact conditions~\cite{Signorini1959,Fichera1972}, illustrated in Fig.~\ref{fig:Signorini}.

  We briefly recall the variational formulation of the contact problem for the unknown displacement field $\boldsymbol{u}=(\boldsymbol{u}_1,\boldsymbol{u}_2)$ between two bodies $\Omega_1$ and $\Omega_2$ potentially in contact along $\Gamma^C$, see Fig.~\ref{fig:contact_pb}, in the framework of linear elasticity under the small strain assumption; we refer the reader to standard textbooks for the formal details~\cite{Kikuchi1988,Wriggers2006,Yastrebov2013}.
      
Let $\boldsymbol{\sigma}(\boldsymbol{u})$ denote the associated Cauchy stress tensor.
 The boundary of $\Omega=\Omega_1\cup\Omega_2$ is partitioned into $\Gamma^D$, $\Gamma^N$ and $\Gamma^C$, on which a Dirichlet condition $\boldsymbol{u}=\boldsymbol{u}_0$ and a Neumann condition $\boldsymbol{\sigma}(\boldsymbol{u})\cdot\boldsymbol{n}=\boldsymbol{f}_S$ are respectively prescribed, together with the (frictionless) unilateral contact conditions on $\Gamma^C$
\begin{equation} \label{eq:cond_contact}
u_N\leq d, \qquad \sigma_N\leq 0, \qquad \sigma_N(u_N-d)=0, \qquad \sigma_T=0,
\end{equation}
where $u_N=\boldsymbol{u}_1\cdot\boldsymbol{n}_1+\boldsymbol{u}_2\cdot\boldsymbol{n}_2$ is the normal gap displacement, $\boldsymbol{n}_1,\boldsymbol{n}_2$ the outward unit normals to $\Omega_1,\Omega_2$ on $\Gamma^C$, $d$ the initial gap between $\Omega_1$ and $\Omega_2$ along $\Gamma^C$, and $\sigma_N$, $\sigma_T$ the normal and tangential components of the contact stress. Similarly, for any $\boldsymbol{v}\in(H^1(\Omega))^m$, we denote $\boldsymbol{v}=(\boldsymbol{v}_1,\boldsymbol{v}_2)$ its restrictions to $\Omega_1,\Omega_2$ and $v_N=\boldsymbol{v}_1\cdot\boldsymbol{n}_1+\boldsymbol{v}_2\cdot\boldsymbol{n}_2$ its normal gap. Let us note that, under the small strain hypothesis, $\Gamma^C$ is fixed and $\boldsymbol{n}_2=-\boldsymbol{n}_1$, so that $u_N=(\boldsymbol{u}_1-\boldsymbol{u}_2)\cdot\boldsymbol{n}_1$ and $v_N=(\boldsymbol{v}_1-\boldsymbol{v}_2)\cdot\boldsymbol{n}_1$.

Introducing the functional space
\begin{equation}
\mathcal{V}=\{\boldsymbol{v}\in (H^1(\Omega))^m \text{ such that } \boldsymbol{v}=\boldsymbol{u}_0 \text{ on } \Gamma^D\},
\end{equation}
and the convex set of admissible displacements
\begin{equation}
\mathcal{V}_C=\{\boldsymbol{v}\in\mathcal{V} \text{ such that } v_N\leq d \text{ on } \Gamma^C\},
\end{equation}
the variational formulation of this problem then reads~\cite{Duvaut1976}:
\begin{equation} \label{eq:ineg_variationnelle}
\begin{cases}
\text{Find } \boldsymbol{u}\in \mathcal{V}_C \text{ such that} \\
a(\boldsymbol{u},\boldsymbol{v}-\boldsymbol{u}) \geq l(\boldsymbol{v}-\boldsymbol{u}) \quad \forall \boldsymbol{v}\in\mathcal{V}_C.
\end{cases}
\end{equation}
The bilinear form $a(\cdot,\cdot)$ and the linear form $l(\cdot)$ are given by
\begin{equation}
a(\boldsymbol{u},\boldsymbol{v})=\int_\Omega \boldsymbol{\varepsilon}(\boldsymbol{v}):\mathbb{C}:\boldsymbol{\varepsilon}(\boldsymbol{u})\,\mathrm{d}\Omega, \qquad
l(\boldsymbol{v})=\int_\Omega \boldsymbol{f}_V\cdot\boldsymbol{v}\,\mathrm{d}\Omega + \int_{\Gamma^N} \boldsymbol{f}_S\cdot\boldsymbol{v}\,\mathrm{d}\Gamma,
\end{equation}
where $\mathbb{C}$ is the fourth-order elasticity tensor, $\boldsymbol{\varepsilon}(\boldsymbol{u})=\frac{1}{2}(\nabla\boldsymbol{u}+\nabla\boldsymbol{u}^T)$ the linearized strain tensor, and $\boldsymbol{f}_V$ the body force.

Problem~\eqref{eq:ineg_variationnelle} can equivalently be recast as a Lagrange multiplier (saddle-point) formulation~\cite{Kikuchi1988}. Let
\begin{equation}
  \mathcal{V}^0=\{\boldsymbol{v}\in (H^1(\Omega))^m \text{ such that } \boldsymbol{v}=\boldsymbol{0} \text{ on } \Gamma^D\},
\end{equation}
and
\begin{equation}
\mathcal{W}^+=\{\mu\in (H^{\frac{1}{2}}_{00}(\Gamma^C))' \text{ such that } \mu \geq 0\},
\end{equation}
where $(H^{\frac{1}{2}}_{00}(\Gamma^C))'$ denotes the dual of $H^{\frac{1}{2}}_{00}(\Gamma^C)$ (see e.g.~\cite{Dautray1988} for the definition of this space).
The saddle-point formulation then reads: find $(\boldsymbol{u},\lambda)\in \mathcal{V}\times\mathcal{W}^+$ such that
  \begin{equation} \label{eq:mult_point_selle}
    \begin{cases}  \begin{aligned}
      a(\boldsymbol{u},\boldsymbol{v}) +b(\boldsymbol{v},\lambda) &= l(\boldsymbol{v}),&\forall\boldsymbol{v}\in\mathcal{V}^0,\\
     b(\boldsymbol{u},\mu - \lambda) &\leq q(\mu - \lambda),&\forall \mu\in\mathcal{W}^+.
   \end{aligned}
 \end{cases}
\end{equation}
with the bilinear form $b(\cdot,\cdot)$ and the linear form $q(\cdot)$ defined through the duality pairing on $H^{\frac{1}{2}}_{00}(\Gamma^C)$
\begin{equation}
 b(\boldsymbol{v},\mu) = \langle v_N,\mu\rangle_{H^{\frac{1}{2}}_{00}(\Gamma^C),(H^{\frac{1}{2}}_{00}(\Gamma^C))'}, \qquad q(\mu)=  \langle d,\mu\rangle_{H^{\frac{1}{2}}_{00}(\Gamma^C),(H^{\frac{1}{2}}_{00}(\Gamma^C))'}.
 \label{eq:def_b_q}
\end{equation}
Problem~\eqref{eq:mult_point_selle} is a saddle-point problem, for which existence and uniqueness have been proven in the literature~\cite{Kikuchi1988}. It is worth noting that the Lagrange multiplier formulation respects all the contact conditions~\eqref{eq:cond_contact}. In particular, the complementarity condition is obtained by choosing $\mu=0$ and $\mu=2\lambda$, giving $b(\boldsymbol{u},\lambda)=q(\lambda)$.

In the discrete setting, the saddle-point frictionless contact problem leads to an inequality system (cf. the non-penetration condition), with a sign-constrained dual variable (cf. the non-adhesion condition), see Eq.~\ref{eq:saddle_point_system}.

\begin{equation} \label{eq:saddle_point_system}
     \begin{cases} 
       \text{Find $(\bU,\blambda) \in \mathbb{R}^N \times (\mathbb{R}^+)^{N_\lambda}$ such that}\\
       \begin{aligned}
	    \bK \bU &+ \bB^\top \blambda &=&\, \bF_{\text{ext}}, \\
        \bB \bU & &\leq& \, \bD,     \\
        \blambda \odot &  (\bB \bU - \bD) &=& \, \bzero
        \end{aligned}
     \end{cases}
   \end{equation}
where the solution vector consists of both the primal ($\bU$) and dual variables ($\blambda$), $N$ and ${N_\lambda}$ correspond, respectively, to the number of associated degrees of freedom. 
Moreover, $\bK$ is the global stiffness matrix, $\bB$ represents the pairing (or contact mapping) matrix between potential contact nodes, 
$\bD$ is the discrete initial gap vector between paired nodes and
$\bF_{\text{ext}}$ denotes the external load vector.
The $\leq$ symbol is used to represent the component-wise
inequality while $\odot$ denotes the Hadamard product (component-wise multiplication). {As a reminder, the Signorini conditions, which are a particular case of complementarity conditions in mathematical optimization, can be expressed as $\bzero \leq (\bD - \bB \bU) \perp \blambda \geq \bzero$, where $\perp$ denotes component-wise complementarity.}

 Set of contact conditions render contact problems inherently non-smooth and nonlinear, presenting significant challenges for their numerical solution.
Various formulations have been developed to numerically model contact problems, each offering distinct trade-offs in terms of accuracy, implementation complexity, and numerical performance \cite{Laursen2003,Wriggers2006,Yastrebov2013}. Among them, we find the penalty method~\cite{Kikuchi1981,WriggersZavarise2004}, Lagrange multipliers method~\cite{Kikuchi1988,Wriggers2006}, augmented Lagrangian technique~\cite{Glowinski1989,SIMO1992,Wriggers1985}, Nitsche method~\cite{Wriggers2008,Chouly2013Nitsche,Chouly2019}, Mortar-based approach~\cite{Puso2004_mortar3D,Popp2012}, interior-point method~\cite{TEMIZER2014, Acary2024}, or domain decomposition based techniques~\cite{Krause2002,MOTA2025}. \\
	
In this study, we particularly focus on Lagrange multipliers and penalty contact formulations.
The Lagrange multiplier (LM) formulation offers high accuracy by satisfying the contact conditions exactly through the introduction of additional unknowns (Lagrange multipliers). It theoretically corresponds to Formulation~\eqref{eq:saddle_point_system}. The  standard solution process consists in solving the system associated to the saddle-point coefficient matrix generally through an iterative active set strategy~\cite{Luenberger2008,Wriggers2006,Abide16} in order to rely on saddle-point equality systems and to verify all the Signorini contact conditions. Addressing large-scale contact problems remains the main limitation of the saddle-point matrix solution process, see this specific state-of-the-art point proposed in~\cite{Epalle2025}.\\
The penalty formulation transforms a constrained problem into an unconstrained one by adding a penalized term to the objective function associated to the variational inequality~\eqref{eq:ineg_variationnelle}. In the discrete framework~\eqref{eq:saddle_point_system}, this amounts to defining the equivalent nodal multiplier~\eqref{eq:lambda_penalty}, where the \texttt{max} function is applied component-wise, and solving only the first equation:
\begin{equation}
  \label{eq:lambda_penalty}
  \begin{aligned}
  \blambda_P = \max\left(0, k_N (\bB \bU - \bD)\right),\\
  \text{where $k_N$ is the penalty parameter ($\gg 1~\text{N/m}$).}
  \end{aligned}
\end{equation}
%
The penalty formulation has the advantage of relying practically only on the primal variable (without introducing additional contact unknowns), but it provides only an approximate enforcement of the non-penetration constraint. Moreover, the accuracy of the solution depends strongly on the choice of the penalty coefficient: larger values yield more accurate solutions.
The standard solution procedure consists of solving the penalized system, associated with the coefficient matrix $\bK + k_N \bB^\top \bB$, within an active set loop to handle the nonlinear projection in Eq.~\eqref{eq:lambda_penalty} (see~\cite{Wriggers2006,Epalle2025} and references therein). However, relatively large penalty coefficients may adversely affect the conditioning of the system~\cite{NOUROMID1986,Wriggers2006}, which represents the main limitation of the standard penalty-based solution strategy.\\


To overcome the main limitations of standard solution procedures based on Lagrange multipliers or penalty formulations, mentionned above, two main families of approaches are considered in the literature. The first one consists of using or developing tailored preconditioners for the standard solution strategies. For applications to the saddle-point matrices arising from the Lagrange multiplier method, see~\cite{Benzi2005,Siefert2006,Franceschini2022,KOTHARI2022}; for the penalized matrix formulation, refer to~\cite{Wriggers2006,DIONE2019,Epalle2025}. The second family relies on iterative splitting schemes. These methods use solve-update fixed-point iterations. Each iteration consists of computing the primal variable from the current dual iterate state, then updating the dual variable using the new primal iterate value. 
Formally speaking, the first step of this process corresponds to a block splitting of the first line of saddle-point system~\eqref{eq:saddle_point_system}, transferring the contact forces term ($\bB^\top \blambda$) to the right-hand side, while the second step depends on the considered numerical contact formulation.
For the Lagrange multiplier formulation, the best-known approach is the so-called Uzawa iterative scheme~\cite{ALART1991}. The update phase then corresponds to a projected dual gradient step involving an augmentation (or relaxation) parameter. For the penalty formulation, a similar two-level iterative scheme can be built, in which the dual update is directly given by Eq.~\eqref{eq:lambda_penalty}. 
To the best of our knowledge, the literature mainly focuses on iterative schemes of this type within the augmented Lagrangian framework combining penalty terms with Lagrange multipliers and leading to Uzawa-type methods~\cite{Glowinski1989,ALART1991,SIMO1992}, rather than on purely penalty-based formulations.
The main limitation of such iterative splitting schemes is the slow convergence of the fixed-point iterations, which is guaranteed only within a narrow parameter range~\cite{ALART1991,Benzi2005}. 

Although various enhancements, such as inexact solves, preconditioning strategies~\cite{Zulehner_2001,Bacuta2006}, and more recently Uzawa acceleration techniques~\cite{HO2017,KANNO2020,HUANG2022}, have been proposed, these methods remain rarely adopted in practical applications due to the aforementioned limitations. 

Nevertheless, displacement-force splitting strategies are particularly appealing, as they involve solving only standard rigidity systems at each iteration: neither saddle-point systems nor ill-conditioned penalized matrices need to be handled. This feature opens the way to the efficient solution of large-scale contact problems using existing scalable solvers. Moreover, only the right-hand side of the system is updated during iterations, which enables the reuse of matrix factorizations or the efficient initialization of iterative solvers, thereby significantly reducing the computational cost of iterations beyond the first.
\\

This work proposes a robust enhancement of iterative splitting schemes through the integration of the so-called Crossed-Secant acceleration/stabilization strategy~\cite{RAMIERE2015}. This approach has already demonstrated its effectiveness in overcoming divergent fixed-point iterations in various computational fields~\cite{Ruiz_2022,INTROINI2024}. The present article provides a theoretical perspective supporting these observations.
In contrast to previous Uzawa acceleration strategies applied to frictionless contact problems in~\cite{KANNO2020,HUANG2022}, the proposed scheme renders the iterative method efficient independently of the augmentation parameter value, even beyond the theoretical bounds of the classical Uzawa algorithm. Moreover, it ensures the numerical convergence of the penalty-based splitting procedure over a wide range of large penalty parameters, making the approach both accurate and practically attractive. 

As a result, two-step iterative schemes that are traditionally slow and highly sensitive to parameter choices are transformed into fast, robust, and essentially parameter-unconstrained solvers.
Finally, although this study focuses on contact mechanics, the proposed enhanced accelerated iterative two-step framework is generic and is expected to extend to a broad class of constrained problems across various application fields.
\\

%

The rest of the paper is organized as follows. Section~\ref{s:algo_unified} introduces a unified framework for fixed-point accelerated iterative splitting algorithms applied to frictionless contact problems.  Section~\ref{sec:std_methods} defines the dual update function for the Lagrange multipliers/Uzawa and penalty formulations.
Section~\ref{s:acceleration} reviews two existing acceleration techniques for the Uzawa algorithm in contact mechanics, namely an accelerated projected gradient scheme~\cite{BECK2009} employed in~\cite{KANNO2020} and the well-known Anderson acceleration method~\cite{Anderson1965} applied in~\cite{HUANG2022}. This section also recalls the Crossed-Secant acceleration method~\cite{RAMIERE2015} and presents a numerical analysis highlighting its potential in the presence of divergent fixed-point iterations.
Section~\ref{s:numerical_exemples} is devoted to numerical experiments on three-dimensional frictionless contact problems between deformable bodies. Both academic and industrial test cases are considered, involving node-to-node and node-to-surface contact pairings. Few insights on multi-domain contact simulations, leveraging parallel solution strategies, are also provided.
Finally, conclusions and perspectives are drawn in the last section.

\section{Unified framework for iterative (accelerated) contact algorithms } \label{s:algo_unified}

We propose here a unified computational framework for solving frictionless contact problems that mainly lies on a
two-step iterative process, where each iteration $i$ is composed of:
\begin{subequations}\label{eq:algo_uni_2steps}
	\begin{align}
		\textbf{Solve step:} \quad 
		& \bK \bU^{i} = \bF_{\text{ext}} - \bB^\top \blambda^{i-1} , 
		\label{eq:unified_solve} \\[0.5em]
               \textbf{Update step:} \quad
                & \textbf{Evaluate $\blambda^i$ from up-to-date iterates}
		\label{eq:unified_update}
	\end{align}
\end{subequations}

At each iteration $i$, the primal variable $\bU^i$ (displacements) is
first obtained by solving the equilibrium equation with the current
estimate of the dual contact variable $\blambda^{i-1}$.
Without loss of generality, we assume that $\bK$ is non-singular in order to guarantee uniqueness. If this condition is not initially satisfied by the finite element model (e.g., due to rigid-body modes), standard procedures can be applied to render the stiffness matrix non-singular (see, for instance,~\cite[Remark 4.2]{NOUROMID1986}). After the solve step, the dual variable $\blambda^i$ is updated according to the current
displacement estimate $\bU^i$ and possibly using dual values from previous iterations ($\blambda^{i-1}, \dots$). The update step, which is the main focus of this study, will be presented and discussed in details below.

As mentioned in the introduction, this splitting strategy is numerically attractive, since the matrix of the linear system to be solved is restricted to the material stiffness coefficient matrix, which does not vary during the iterations, in contrast to some splitting formulation as the Uzawa-like solution for the augmented Lagrangian formulation, where an inner loop is required during the solve step to update the contact status (see, for example, \cite{ALART1991,SIMO1992}).

\subsection{Standard update step} \label{sec:standard-update-step}
The most standard formulations for updating the dual variable can be summarized by the following equation
\begin{equation} \label{eq:generic_F}
  {\blambda}^{i} = \mathcal{F}(\blambda^{i-1},\bU^{i}).
\end{equation}
%
A fundamental step in inequality-constraint problems, such as contact problems, is the projection step, which typically enforces the sign constraint on $\blambda$ at each iteration. 
In the generic formulation~\eqref{eq:generic_F}, this step is included in the definition of function $\mathcal{F}$. However, as this projection will play a crucial role in the proposed framework, we propose to rewrite the standard update step as

\begin{subequations}\label{eq:update_proj}
  \begin{align}[left={\textbf{Update step} \empheqlbrace}]
    \hat{\blambda}^{i} &= \mathcal{G}(\blambda^{i-1},\bU^{i}),
        \label{eq:update_G}\\
         {\blambda}^{i} &=\Pi_{\mathbb{R}^+} (\hat{\blambda}^{i}), \label{eq:update_Pi}
   \end{align}
\end{subequations}
where the function $\mathcal{G}$ generically denotes an update formula, whose expression depends on the chosen method. The detailed expressions of $\mathcal{G}$ considered in this paper are given in Section~\ref{sec:std_methods}.
The symbol $\Pi_{\mathbb{R}^+}$ corresponds to the projection operator onto the non-negative orthant $\mathbb{R}^+$.
It is defined as
\begin{equation*}
  \label{eq:def_pi}
   \begin{aligned}
     \Pi_{\mathbb{R}^+}(\mathbf{x}) &= \max({0}, \mathbf{x}), \quad \text{applied component-wise.}\\
     &\text{or equivalently}\\
     \left(\Pi_{\mathbb{R}^+}(\mathbf{x})\right)_j& = \max(0,x_j), \quad j=1,\dots,n;\; \mathbf{x} \in \mathbb{R}^n
   \end{aligned}
 \end{equation*}
 Hence, mathematically speaking,
\begin{equation}
	\label{eq:def_F}
	\mathcal{F} = \Pi_{\mathbb{R}^+} \circ \mathcal{G}.
\end{equation}

It is worth emphasizing that the projection step, applied to the dual contact variable $\blambda$, automatically determines a set of active constraints (cf. the non-smooth optimization problem) at each iteration, which converges as the algorithm progresses.
{Consequently, no additional loop (outer or active-set) is required to determine the active constraints set.}

\subsection{Accelerated update step} \label{sec:accel-update-step}
Looking at the update expression in Eq.~\eqref{eq:update_proj}, as $\bU^i$ depends
on $\blambda^{i-1}$ via the solve step, see Eq.~\eqref{eq:unified_solve}, we can formally write:
\begin{equation}
  \label{eq:point-fixe_lambda}
  \begin{aligned}
  \blambda^{i} &= \Pi_{\mathbb{R}^+} \circ \mathcal{G}(\blambda^{i-1},\bU^{i}) = \mathcal{F}(\blambda^{i-1},\bU^{i}) \\ 
  &=
  \mathcal{F}\left(\blambda^{i-1},\bK^{-1} (\bF_{\text{ext}} - \bB^\top
    \blambda^{i-1}) \right) \\
  &=\tilde{\mathcal F}(\blambda^{i-1})
  \end{aligned}
\end{equation}
Then, fixed-point acceleration algorithms can be applied to find $\blambda$, solution of the fixed-point equation $\blambda = \tilde{\mathcal F}(\blambda)$.
However, since standard acceleration methods construct an affine
(rather than convex) combination at each step, the accelerated iterate
may lie outside the admissible set.
We therefore propose to also apply the projector $\Pi_{\mathbb{R}^+}$
after the acceleration step in order to ensure that the constraints on
$\blambda$ are satisfied.
Although absent from the literature on acceleration method for contact problems~\cite{KANNO2020,HUANG2022}, this appears to be a highly appropriate step for constrained problems~\cite{Li_Accel_2020}.\\

We propose to also consider fixed-point acceleration strategies applied directly to the function $\mathcal{G}$ rather than to $\mathcal{F}$. This allows us to work with a smooth operator (instead of the composite non-smooth mapping $\mathcal{F}$).
From a mathematical standpoint, the fixed point of $\mathcal{F}$ coincides with that of $\mathcal{G}$ once the active constraints have been identified.
{This opens the way to acceleration strategies that can accommodate less restrictive numerical parameters while avoiding a strong dependence on the fixed-point residual definition. Indeed, prior to active-set identification, the iteration is governed by the projected mapping $\Pi_{\mathbb{R}^+} \circ \mathcal{G}$ rather than $\mathcal{G}$ itself. However, acceleration techniques capable of stabilizing iterations beyond the standard fixed-point convergence regime typically require continuous -- or locally smooth -- operators, which further motivates considering $\mathcal{G}$ instead of the composite non-smooth mapping $\mathcal{F}$.}

During the iterative process, when acceleration is applied to $\mathcal{G}$, the iterates may leave the admissible set, and a projection step is therefore required after each iteration.\\

In what follows, the function $\mathcal{A}$ generically denotes an acceleration procedure, where the parameters in square brackets are method-dependent and possibly not mandatory. Several choices will be detailed in Section~\ref{s:acceleration}.\\

Thus, depending on the considered fixed-point equation, two accelerated update step can be derived. For the acceleration on the non-smooth function $\mathcal{F}=\Pi_{\mathbb{R}^+} \circ \mathcal{G}$, we get:
\begin{subequations}\label{eq:update_accel_F}
  \begin{align}[left={\textbf{Update step} \empheqlbrace}]
   \hat{\blambda}^{i} &=\Pi_{\mathbb{R}^+} \circ \mathcal{G}(\blambda^{i-1},\bU^{i}) \\
     {\blambda}^{i} &=\Pi_{\mathbb{R}^+} \circ \mathcal{A}\left(\hat{\blambda}^{i},\hat{\blambda}^{i-1},[\hat{\blambda}^{i-2,}\dots,\hat{\blambda}^{0},{\blambda}^{i-1},\dots, {\blambda}^{0}] \right)
	\end{align}
      \end{subequations}

\noindent For the acceleration on the smooth function $\mathcal{G}$, we have:
\begin{subequations}\label{eq:update_accel_G}
  \begin{align}[left={\textbf{Update step} \empheqlbrace}]
   &  \hat{\blambda}^{i} = \mathcal{G}(\blambda^{i-1},\bU^{i}),\\
 & {\blambda}^{i} =\Pi_{\mathbb{R}^+} \circ \mathcal{A}\left(\hat{\blambda}^{i},\hat{\blambda}^{i-1},[\hat{\blambda}^{i-2,}\dots,\hat{\blambda}^{0},{\blambda}^{i-1},\dots, {\blambda}^{0}] \right) 
	\end{align}
      \end{subequations}

\subsection{Unified iterative solution algorithm}
\label{sec:unified_steps}
By gathering the different variants of the update steps, we propose the following unified pseudo-code for solving frictionless contact problems using an iterative splitting scheme, see Algorithm~\ref{alg:algo_uni}.

\renewcommand{\algorithmicensure}{\textbf{Initialization:}}
\begin{algorithm}
	\caption{Unified (accelerated) fixed-point framework for frictionless contact}
	\label{alg:algo_uni}
	\begin{algorithmic}[1]
          \REQUIRE  initial dual contact variable $\blambda^0$, boolean $\textit{acceleration}$, boolean $\textit{accel\_F}$, integer $\textit{minit\_accel}$
          \ENSURE  $\textit{convergence} \gets \text{false}$, $i \gets 1$, $ \hat{\blambda}^0 \gets \blambda^0$
		\WHILE{not \textit{convergence}}
		\STATEx 
		\quad \textbf{solve step:}
                \STATE \label{step:algo_uni_solve} \quad $\bK \bU^{i} = \bF_{\text{ext}} - \bB^\top\blambda^{i-1}$ 	
		\STATEx 
		\quad \textbf{update step:} 
		\STATE \label{step:algo_uni_updt} \quad $\hat{\blambda}^{i} = \mathcal{G}(\blambda^{i-1},\bU^i)$ 
                \STATE \quad \textbf{if} \textit{acceleration} \textbf{then} 
                \STATE \quad \quad \textbf{if} \textit{accel\_F} \textbf{then}
                \STATE \quad \quad \quad $\hat{\blambda}^{i} \gets \Pi_{\mathbb{R}^+}(\hat{\blambda}^{i})$
                \STATE \quad \quad \textbf{endif}
                \STATE \quad \quad \textbf{if} ($i \ge \textit{minit\_accel}$) \textbf{then} 
                \STATE \quad \quad \quad ${\blambda^{i}}=\mathcal{A}\left(\hat{\blambda}^{i},\hat{\blambda}^{i-1},[\hat{\blambda}^{i-2,}\dots,\hat{\blambda}^{0},{\blambda}^{i-1},\dots, {\blambda}^{0}] \right)$                
                \STATE \quad \quad \textbf{else}
                \STATE \quad \quad \quad $\blambda^{i} = \hat{\blambda}^{i}$
                \STATE \quad \quad \textbf{endif}
                \STATE \quad \textbf{else}
                \STATE \quad \quad $\blambda^{i} = \hat{\blambda}^{i}$
                \STATE \quad \textbf{endif} 
                \STATE \label{step:proj_final} \quad ${\blambda^{i}} \gets \Pi_{\mathbb{R}^+}({\blambda^{i}})$
		
		\STATEx \label{step:algo_uni_convergence}
		\quad \textbf{check convergence}:
                \STATE \quad \textbf{if} convergence criterion is respected \textbf{then}
                \STATE \quad \quad $\textit{convergence} \gets \text{true}$
                \STATE \quad \textbf{endif} 	
		\STATE $i \gets i + 1$
		\ENDWHILE
	\end{algorithmic}
\end{algorithm}

{In this Algorithm, $\blambda^0$ denotes the initial dual contact variable, the boolean variable $\textit{acceleration}$ indicates whether acceleration is enabled, while $\textit{accel\_F}$ activates the acceleration on the function $\mathcal{F}$ (otherwise, the acceleration is applied to $\mathcal{G}$). Finally, $\textit{minit\_accel}$ denotes the iteration number from which the acceleration procedure is activated.}

We can see that the core of the algorithm -- the solve step (see line~\ref{step:algo_uni_solve}), the evaluation of $\mathcal{G}$, the projection in the update step (see lines~\ref{step:algo_uni_updt} and~\ref{step:proj_final}), as well as the convergence check -- remains fixed. The remaining parts of the update step are modular, emphasizing the plug-and-play nature of the algorithm.

\section{Dual update function} \label{sec:std_methods}
In this section, we give the definition of the function $\mathcal{G}$ for both the Lagrange multipliers and penalty formulations. Both formulations require a numerical parameter, for which we recall the theoretical bounds ensuring the convergence of the standard update procedure, see Eq.~\eqref{eq:update_proj}.

\subsection{Uzawa algorithm for Lagrange multipliers formulation}  \label{s:LM} 
The expression of $\mathcal{G}$ for the Lagrange multiplier formulation returns to the classical Uzawa algorithm~\cite{UZAWA1958}, see Eq.~\eqref{eq:uzawa_std_update}.
\begin{equation}
  \label{eq:uzawa_std_update}
  \mathcal{G}(\blambda,\bU) = \blambda + \rho \left(\bB \bU - \bD\right), \quad \rho > 0
\end{equation}
As already mentioned in the introduction, this update corresponds gradient ascent step~\cite{CIARLET1989} involving the augmentation (or relaxation) parameter $\rho$.\\

A sufficient condition on $\rho$ ensuring the convergence of the Uzawa scheme for  problems with inequality constraints of the form~\eqref{eq:saddle_point_system} is given by:
\begin{equation}
	0 < \rho < \frac{2 \mu_{\text{min}}(\bK)}{\left \| \bB \right \|_2}
	\label{eq:uzawa_rho_bounds}
\end{equation}
where $\mu_{\text{min}}(\bK)$ denotes the smallest eigenvalue of the stiffness matrix $\bK$, and $\left \| \bB \right \|_2$ refers to the Euclidian matrix norm of the contact pairing matrix $\bB$, see~\cite{ALART1991} for example. \\
Although the theoretical bounds~\eqref{eq:uzawa_rho_bounds} on the Uzawa parameter guarantee convergence, they are often overly restrictive and difficult to evaluate in practice, making the selection of an efficient parameter value nontrivial. Even when these bounds are satisfied, the Uzawa algorithm exhibits only slow linear convergence, due to a large contraction factor~\cite{Nochetto2004}, with the convergence rate further deteriorating as $\rho$ decreases. Moreover, choosing a large value of $\rho$ may result in divergence. Consequently, substantial effort has been devoted over the years to identifying optimal choices of this parameter, see for instance~\cite{ALART1991,KANNO2020,Nochetto2004,Bacuta2006}.

One of the main contributions of this work is to show that an appropriate acceleration strategy, namely the Crossed-Secant method~\cite{RAMIERE2015}, {applied to the smooth function $\mathcal{G}$ (cf. update step in Eq.~\eqref{eq:update_accel_G})}, renders the Uzawa scheme both efficient and essentially insensitive to the choice of the augmentation parameter. 

\subsection{Penalty-based fixed-point formulation} \label{s:penalty}
Relying of the proposed expression of the equivalent Lagrange multiplier, see Eq.~\eqref{eq:lambda_penalty}, the definition of $\mathcal{G}$ for the penalty formulation reads:
\begin{equation}
  \label{eq:penalty_fp_update}
\mathcal{G}({\blambda}, \bU) = k_N  \left({\bB} \bU - {\bD}\right), \quad k_N \gg 1
\end{equation}

To establish theoretical bounds on the penalty parameter ensuring convergence of the penalty-based splitting method, we exploit the mathematical equivalence between the penalty formulation and a so-called regularized saddle-point system~\cite{Luenberger2008,Benzi2005}. This regularized saddle-point system is characterized by the following $2 \times 2$ block matrix:
\begin{equation}
  \label{eq:regularized_SP}
\begin{bmatrix}
      \bK & \bB^\top \\
      \bB & -\frac{1}{k_N} \mathbf{Id}
\end{bmatrix}
\end{equation}
{with $\mathbf{Id}$ denoting the identity matrix.}
Applying an Uzawa scheme to this problem and choosing the augmentation parameter as $k_N$ yields:
\begin{equation*}
  \mathcal{G}(\blambda,\bU) = \blambda + k_N \left(\bB \bU - \frac{1}{k_N} \blambda - \bD\right).
\end{equation*}
This formally corresponds to~\eqref{eq:penalty_fp_update}. Consequently, the penalty-based splitting procedure can be interpreted as an Uzawa algorithm applied to the regularized saddle point system~\eqref{eq:regularized_SP}, and the theoretical bounds~\eqref{eq:uzawa_rho_bounds} then hold for~$k_N$.

At first glance, this algorithm may appear of limited interest, since the theoretical convergence bounds impose a strict upper limit on the admissible values of $k_N$, while achieving an accurate penalized solution requires $k_N$ to be as large as possible. 

However, within the unified framework introduced in Section~\ref{s:algo_unified}, this formulation can be combined with the advanced Crossed-Secant acceleration technique~\cite{RAMIERE2015}, which, {when applied to $\mathcal{G}$}, effectively prevents the divergence of highly penalized fixed-point iterations. As a result, penalty coefficients can be selected well beyond the theoretical convergence limits, enabling the computation of accurate penalized solutions.

\section{Acceleration methods} \label{s:acceleration}
In this work, we focus on three acceleration techniques that are tested and compared on numerical examples in Section~\ref{s:numerical_exemples}. First, we consider the FISTA method with adaptive restart and the Anderson acceleration method, which have recently been applied in the literature to accelerate the Uzawa algorithm for frictionless contact problems, see respectively~\cite{KANNO2020} and~\cite{HUANG2022}. 
Although these strategies have been previously introduced, they have not yet been systematically compared. 
Moreover, motivated by the adaptive restart strategy employed in~\cite{KANNO2020}, which was shown to improve efficiency, we propose to incorporate a similar adaptive restart mechanism into the Anderson acceleration scheme. 
While these methods enhance convergence rates, their performance remains sensitive to the choice of the Uzawa augmentation parameter, often requiring careful, problem-dependent tuning.

To overcome this limitation, we investigate the Crossed-Secant method~\cite{RAMIERE2015}, an acceleration strategy that has demonstrated strong numerical performance in situations involving divergent fixed-point iterations~\cite{Ruiz_2022,INTROINI2024}. In addition, the Crossed-Secant method enables the penalty-based splitting strategy to achieve efficient convergence.

The detailed definitions of the acceleration operators (corresponding to the function $\mathcal{A}$ in Algorithm~\ref{alg:algo_uni}) are provided in the following subsections. Some theoretical insights on the Crossed Secant method are also provided.


\subsection{FISTA} \label{ss:fista}

As mentioned in the introduction, the Uzawa method can be interpreted
as a projected gradient scheme. Consequently, the Fast Iterative
Shrinkage/Thresholding Algorithm (FISTA) originally developed for
proximal gradient operators~\cite{BECK2009} was successfully applied
in~\cite{KANNO2020} to improve the convergence of the Uzawa method in
frictionless contact problems.
By incorporating Nesterov-type extrapolation~\cite{NESTEROV1983},
FISTA achieves the optimal convergence rate $O\!\left(1 / k^{2}
\right)$ in objective values.
FISTA is considered as one of the most
important algorithm of the past decades.\\
Following~\cite{KANNO2020}, we combine FISTA with an adaptive restart scheme~\cite{Odonoghue2015} based on a generalized gradient definition. 
As justified by the authors of~\cite{Odonoghue2015}, the restart is triggered when the momentum term and the negative generalized gradient form an obtuse angle, indicating that the momentum seems to be taking the new iterate in an unfavorable direction.
This acceleration scheme is outlined in Algorithm~\ref{alg:FISTA}.
Furthermore, in accordance with its original design for projected gradient methods, FISTA has been applied in~\cite{KANNO2020} to the operator $\mathcal{F} = \Pi_{\mathbb{R}^+}\circ \mathcal{G}$ (see Section~\ref{sec:accel-update-step}). Applying FISTA directly to $\mathcal{G}$ does not appear to be consistent with its intended design.


\begin{algorithm}
	\caption{FISTA with adaptive restart scheme: $\mathcal{A}\left(\hat{\blambda}^{i},\hat{\blambda}^{i-1}\right) $, $\textit{minit\_accel}$ = 2}
	\label{alg:FISTA}
	\begin{algorithmic}[1]
          \ENSURE $\tau^1 = 1$
		\IF{$ \left((\bB \bU^{i} - \bD)\cdot (\hat{\blambda}^{i} - \hat{\blambda}^{i-1}) \right) \ge 0$}  \label{step:accel_kanno_descent} 
		\STATE \quad $\tau^{i} = \frac{1}{2} \left(1 + \sqrt{1+4 (\tau^{i-1})^2} \right)$
		\STATE \quad $\beta^{i} = \dfrac{\tau^{i-1} - 1}{\tau^{i}}$
		\STATE \label{step:fista_extrapolation} \quad ${\blambda}^{i} = \hat{\blambda}^{i} + {\beta^{i}}(\hat{\blambda}^{i} - \hat{\blambda}^{i-1}) $
		\ELSE \label{step:accel_kanno_restart} 
		\STATE \quad $\tau^{i} = 1$
		\STATE \quad ${\blambda}^{i} = \hat{\blambda}^{i}$
		\ENDIF
	\end{algorithmic}
\end{algorithm}

\subsection{Anderson acceleration (AA) method} \label{ss:accel_anderson}

The well-known Anderson acceleration (AA) strategy \cite{Anderson1965, Walker2011} has been successfully applied to improve the convergence of Uzawa iterations in various contexts \cite{HO2017,HUANG2022}. However, to the best of our knowledge, these works consider the fixed-point system on $(\bU,\blambda)$ (cf. Eqs.~\eqref{eq:unified_solve} and \eqref{eq:generic_F}) rather than Eq.~\eqref{eq:point-fixe_lambda}, which defines a fixed-point equation solely for the dual variable. Accelerating a single physical variable typically enhances both stability and efficiency of the acceleration procedure in multiphysics or multivariable couplings, as illustrated for example in \cite{erbts-duester2012,INTROINI2024}.

The memory parameter (often denoted by $m$ in AA) is difficult to choose, as it requires balancing convergence performance against memory footprint. Following the results reported in~\cite{HUANG2022} for contact problems, small values of $m$ appear to be a good compromise and already provide significant acceleration. {This is consistent with the recent conclusions of~\cite{SCHAPIRA2023}, drawn from a different setting, 
that the choice of $m=1$ is the best with  regard to the aspects of robustness and simplicity of implementation.} Algorithm \ref{alg:accelAnderson} outlines the one-step AA applied to the dual variable. This approach, denoted here as Anderson-1, is also referred to in the literature as the Alternate Secant method~\cite{RAMIERE2015}, as it coincides with a vector secant method.

\begin{algorithm}
	\caption{Anderson-1: $\mathcal{A}\left(\hat{\blambda}^{i},\hat{\blambda}^{i-1},\blambda^{i-1},\blambda^{i-2} \right)$, $\textit{minit\_accel}$ = 2}
	\label{alg:accelAnderson}
	\begin{algorithmic}[1]
          \ENSURE $\bdelta^1 = \hat{\blambda}^{1} - \blambda^{0} $ 
		\STATE $\bdelta^{i} = \hat{\blambda}^{i} - \blambda^{i-1}$
		\STATE $\beta^i = \dfrac{( \bdelta^{i-1}-\bdelta^{i}) \cdot \bdelta^{i} }{\left\| \bdelta^{i} - \bdelta^{i-1} \right\|_2^2}$
		\STATE ${\blambda}^{i} = \hat{\blambda}^{i} + \beta^i \left(\hat{\blambda}^{i} -\hat{\blambda}^{i-1}\right) $
	\end{algorithmic}
\end{algorithm}
 
Thanks to the unified formalism, we can see that Anderson-1 and FISTA are similar momentum-based methods, differing mainly in their choice of the inertial parameter $\beta^i$. Although combinations of Anderson acceleration with restarting or windowing strategies have been proposed in the literature, no clear strategy has emerged~\cite{Ouyang_AA_2024,Krzysik_AA_2025}. 
In this context, combining Anderson-1 with an adaptive restart mechanism originally developed for projected gradient methods, as proposed in~\cite{Odonoghue2015}, appears to be a promising enhancement. This approach is detailed in Algorithm~\ref{alg:accelAnderson_restart}. Nevertheless, it should be emphasized that the global convergence of Anderson acceleration is not guaranteed in the general case~\cite{Krzysik_AA_2025}.

\begin{algorithm}
	\caption{Anderson-1 with adaptive restart scheme: $\mathcal{A}\left(\hat{\blambda}^{i},\hat{\blambda}^{i-1},\blambda^{i-1},\blambda^{i-2} \right)$, $\textit{minit\_accel}$ = 2}
	\label{alg:accelAnderson_restart}
	\begin{algorithmic}[1]
          \ENSURE $\bdelta^1 = \hat{\blambda}^{1} - \blambda^{0} $ 
		\IF{$ \left((\bB \bU^{i} - \bD) \cdot (\hat{\blambda}^{i} - \hat{\blambda}^{i-1}) \right) \ge 0$}  
		\STATE \quad $\bdelta^{i} = \hat{\blambda}^{i} - \blambda^{i-1}$
		\STATE \quad $\beta^i = \dfrac{( \bdelta^{i-1}-\bdelta^{i}) \cdot \bdelta^{i} }{\left\| \bdelta^{i} - \bdelta^{i-1} \right\|_2^2}$
		\STATE \quad ${\blambda}^{i} = \hat{\blambda}^{i} + \beta^i \left(\hat{\blambda}^{i} -\hat{\blambda}^{i-1}\right) $
		\ELSE 
		\STATE \quad ${\blambda}^{i} = \hat{\blambda}^{i}$
		\ENDIF
	\end{algorithmic}
\end{algorithm}

{Finally, Anderson acceleration fundamentally relies on the residual structure associated with a well-defined fixed-point mapping and is therefore naturally suited to accelerating the operator $\mathcal{F}$ (see Section~\ref{sec:accel-update-step}), but less appropriate for $\mathcal{G}$ prior to active-set identification.}

\subsection{Crossed-Secant (CS) acceleration method} \label{ss:accel_cs}

The Crossed-Secant acceleration method is a one-step acceleration technique that can be interpreted as a vector secant method. It differs from the one-step Anderson algorithm in the definition of the vector inverse, see~\cite{RAMIERE2015} for further details.
As illustrated in Algorithm~\ref{alg:accelCS}, the method can also be interpreted as a dynamic relaxation scheme, which explains why it is often referred to in the literature as dynamic Aitken relaxation~\cite{kuttler-wall2008,erbts-duester2012}. Indeed, introducing $\omega^i = 1 - \beta^i$, we obtain:

\begin{equation}
  \label{eq:CS_relax}
{\blambda}^{i} =  \omega^i \hat{\blambda}^{i} + (1 - \omega^i) \blambda^{i-1}.
\end{equation}

\begin{algorithm}
	\caption{Crossed-Secant acceleration: $\mathcal{A}\left(\hat{\blambda}^{i},\hat{\blambda}^{i-1},\blambda^{i-1},\blambda^{i-2} \right)$, $\textit{minit\_accel}$ = 2}
	\label{alg:accelCS}
	\begin{algorithmic}[1]
          \ENSURE $\bdelta^1 = \hat{\blambda}^{1} - \blambda^{0} $ 
          \STATE $\bdelta^{i} = \hat{\blambda}^{i} - \blambda^{i-1}$
          \STATE $\beta^i = \dfrac{(\hat{\blambda}^{i} -\hat{\blambda}^{i-1}) \cdot (\bdelta^{i} - \bdelta^{i-1})}{\left\| \bdelta^{i} - \bdelta^{i-1} \right\|_2^2}$
	\STATE ${\blambda}^{i} =  \hat{\blambda}^{i} - \beta^i \bdelta^{i} $
	\end{algorithmic}
\end{algorithm}

Thanks to its relaxation analogy, see Eq.~\eqref{eq:CS_relax}, we have
\begin{equation} \label{eq:CS_residual}
    {\blambda}^{i} - {\blambda}^{i-1} 
    = \omega^i (\hat{\blambda}^{i} -  \blambda^{i-1})
\end{equation}
where $(\hat{\boldsymbol{\lambda}}^{i} - \boldsymbol{\lambda}^{i-1})$ corresponds to the fixed-point residual $\boldsymbol{\delta}^{i}$ at iteration $i$. Assuming $\omega^i \neq 0$ (otherwise the sequence has already converged, since $\boldsymbol{\lambda}^{i-1} = \boldsymbol{\lambda}^{i-2}$), the convergence of the Crossed-Secant--accelerated sequence implies convergence toward the fixed-point solution.

We now give some insights on the behavior of the Crossed-Secant--accelerated sequence.
\begin{equation}
  \label{eq:demo}
  \begin{aligned}
    \left\|{\blambda}^{i} - {\blambda}^{i-1} \right \|_2
    &= \left|\omega^i\right| \left \| (\hat{\blambda}^{i} -  \blambda^{i-1})\right\|_2\\
    &=\dfrac{\left|({\blambda}^{i-2} -{\blambda}^{i-1}) \cdot (\bdelta^{i} - \bdelta^{i-1})\right|}{\left\| \bdelta^{i} - \bdelta^{i-1} \right\|_2^2} \left\|{\bdelta}^{i}\right \|_2\\
    &\le \dfrac{\left \|{\bdelta}^{i}\right \|_2}{\left\|\bdelta^{i} - \bdelta^{i-1}\right \|_2}\left\|{\blambda}^{i-2} -{\blambda}^{i-1} \right\|_2
  \end{aligned}
\end{equation}
Then, the accelerated residual decreases whenever
$\dfrac{\left\|\bdelta^{i}\right\|_2}{\left\|\bdelta^{i} - \bdelta^{i-1}\right\|_2} < 1,
$
which is equivalent to
\begin{equation}
  \label{eq:inf_vect}
  \left\|\bdelta^{i}\right\|_2 < \left\|\bdelta^{i} - \bdelta^{i-1}\right\|_2.
\end{equation}

\noindent Squaring both sides and using the expansion of the squared norm, we have
\begin{equation*}
    \left\|{\bdelta}^{i}\right\|_2^2 < \left\|\bdelta^{i}\right\|_2^2-2 \bdelta^{i}\cdot \bdelta^{i-1}  + \left\|\bdelta^{i-1}\right\|_2^2
\end{equation*}
so that Condition~\eqref{eq:inf_vect} can be rewritten as
\begin{equation}
  \label{eq:cond_finale}
    2 \bdelta^{i}\cdot \bdelta^{i-1} < \left\|\bdelta^{i-1}\right\|_2^2 
  \end{equation}
  
Several cases can be considered:
\begin{itemize}
\item If $\bdelta^{i}\cdot \bdelta^{i-1} \le 0$, then Condition~\eqref{eq:cond_finale} is always satisfied. This corresponds to the fixed-point residuals forming an obtuse angle, indicating that successive iterates approach the fixed point (possibly without convergence) from alternating directions.\\
So even if the fixed point diverges locally, the Crossed-Secant method can reduce the norm of the accelerated residual sequence, thereby helping to guide the iterations toward the fixed point.
In contrast, Anderson-1 generates a negative \(\beta^i\) in this case, more specifically $\beta^i \in ]-1,0[$, which stabilizes the step but does not guarantee a decrease in the residual norm. Hence, Anderson-1 mitigates overshoot
without actively correcting local divergence, unlike Crossed-Secant.

\item If $\bdelta^{i} \cdot \bdelta^{i-1} > 0$ (acute angle), then $\cos(\bdelta^{i},\bdelta^{i-1}) > 0$, and Condition~\eqref{eq:cond_finale} is satisfied if
  \begin{equation}
    \label{eq:angle_aigu}
    \left\|\bdelta^{i}\right\|_2 < \frac{1}{2 \cos(\bdelta^{i},\bdelta^{i-1})} \, \left\|\bdelta^{i-1}\right\|_2.
  \end{equation}
Particularly, if the angle between $\bdelta^{i}$ and $\bdelta^{i-1}$ is less than $\pi/3$, then the fixed-point residual must decrease. In the limiting case of colinear vectors, the new fixed-point residual must be less than half of the previous one.
\end{itemize}

Hence, the convergence of the Crossed-Secant iterates is generally non-monotone, and the non-verification of Condition~\eqref{eq:cond_finale} acts as an annealing criterion. This mechanism may temporarily cause overshoot or increase the distance between successive iterates, helping the sequence escape local stagnation or unfavorable search directions, and guiding it toward the fixed-point solution.
While convergence is not guaranteed in general, these overshoots can improve the robustness of the iterations and help avoid stagnation near suboptimal iterates.

\begin{remark} Applied to the unprojected Uzawa scheme, the Crossed-Secant method reduces to the well-known Barzilai-Borwein~(BB) method~\cite{Barzilai_Borwein_1988}, an iterative gradient approach for unconstrained optimization.\\
Indeed, Eq.~\eqref{eq:CS_relax} can be reformulated as
\begin{equation} \label{eq:CS_relax_2}
  {\blambda}^{i} = \blambda^{i-1}+ \omega^i {\bdelta}^{i}.
\end{equation}
Using the Uzawa definition~\eqref{eq:uzawa_std_update}, we have
\begin{equation}\label{eq:delta_Uzawa}
  \bdelta^i = \rho (\bB \bU^i - \bD) = \rho \bg^i
\end{equation}
with the definition of the generalized gradient $\bg^i = (\bB \bU^i - \bD)$. Then
\begin{equation}\label{eq:CS_Uzawa_wp}
  \begin{aligned}
    \omega^i \rho &= \dfrac{({\blambda}^{i-2} -{\blambda}^{i-1}) \cdot (\bdelta^{i} - \bdelta^{i-1})}{\left\| \bdelta^{i} - \bdelta^{i-1} \right\|_2^2}\rho\\
    &=  \dfrac{({\blambda}^{i-2} -{\blambda}^{i-1}) \cdot (\bg^{i} - \bg^{i-1})}{\left\| \bg^{i} - \bg^{i-1} \right\|_2^2}\\
    &= - \alpha^i \quad \text{of Barzilai-Borwein~\cite[Eq.(5)]{Barzilai_Borwein_1988}}
  \end{aligned}
\end{equation}
Finally, we get
\begin{equation} \label{eq:CS_BB}
  {\blambda}^{i} = \blambda^{i-1}- \alpha^i {\bg}^{i}.
\end{equation}
which coincides exactly with the BB method. Hence, $\omega^i = (1 - \beta^i)$ can be interpreted as a correction factor for the Uzawa augmentation parameter $\rho$, leading to a dynamic optimal coefficient given by $\omega^i \rho = (1 - \beta^i)\rho$. \hfill $\square$
\end{remark}

\section{Numerical results} \label{s:numerical_exemples}

In this section, we illustrate the results obtained using the various proposed strategies  through a series of three-dimensional mechanical contact problems without friction, including an academic Hertzian benchmark (Section~\ref{ss:Hertz}) and an industrial thermal-induced contact case where contact is induced by thermal expansion (Section~\ref{ss:IPG}). We also investigate the performance of the best retained strategies as the number of contacting domains increases, highlighting key observations and trends in case of parallel computations (see Section~\ref{ss:multi_spheres}).

Numerical experiments presented in this section are performed using the \textsc{Cast3M}~\cite{CEA_CAST3M_2023} finite element solver developed at CEA (French Alternative Energies and Atomic Energy Commission). In all simulations, the resulting linear systems were solved using the direct linear solver (based here on a Crout factorization) available in \textsc{Cast3M}. 
The computations were run on a dual-socket Intel Xeon Silver 4114 workstation clocked at 2.20 GHz, equipped with 64 GB of RAM. 

For all numerical examples and tested methods, we define the convergence criterion by

\begin{equation}\label{eq:relative_residual}
	{r}^i = \frac{\| \blambda^{i} - \blambda^{i-1} \|_2}{ \| \blambda^{i} \|_2 }.
\end{equation}
The convergence is verified if ${r}^i  \le \epsilon$ with the convergence tolerance set to $\epsilon~=~10^{-12}$, see Algorithm~\ref{alg:algo_uni}.

Moreover, the theoretical sufficient upper bound given in Eq.~\eqref{eq:uzawa_rho_bounds} for the augmentation parameter $\rho$ or the penalty coefficient $k_N$ has been evaluated using $\left\| \bB \right \|_2 = 1$, as suggested in the literature~\cite{ALART1991,HUANG2022}.
For the Hertzian contact problem, we found an upper bound of approximately $3.9 \cdot 10^3$, for the industrial test case the upper bound is about $2 \cdot 10^3$, while for the multi-domain contact case it is around $3.2 \cdot 10^4$. From a numerical perspective, the standard Uzawa and penalty-splitting algorithms can sometimes handle higher values before divergence is observed. Accordingly, for numerical simulations, the threshold guaranteeing convergence is of $10^4$ for the Hertzian contact problem, $10^3$ for the industrial example, and $10^5$ for the multi-domain contact cases.

For both Uzawa and penalty-splitted fixed-point formulations, we compare the different acceleration techniques introduced in Section~\ref{s:acceleration}:
\begin{itemize}
	\item \textbf{FISTA + AR}: FISTA scheme with adaptive restart (AR), see Algorithm~\ref{alg:FISTA};
	\item \textbf{Anderson-1}: Anderson-1 acceleration, see Algorithm~\ref{alg:accelAnderson};
	\item \textbf{Anderson-1 + AR}: Anderson-1 acceleration with adaptive restart, see Algorithm~\ref{alg:accelAnderson_restart};
	\item \textbf{Crossed-Secant}: Crossed-Secant acceleration, see Algorithm~\ref{alg:accelCS}.
\end{itemize} 

\noindent \textbf{Accuracy assessement.}  

\begin{itemize}
	\item To assess accuracy of the proposed schemes in terms of Signorini contact conditions satisfaction, the following quantities are evaluate: \\
	\emph{(i)} \textbf{effective gap}
	\begin{equation}
	{\bg}^i_{\text{max}} = \max_{\mathcal{Z}^{i-1}} |\bB \bU^i - \bD|
      \end{equation}
      on a so-called active contact zone, identified at the beginning of iteration $i$ as $\mathcal{Z}^{i-1} = \{ j \in [\![1; N_{\lambda}]\!] ; \lambda_j^{i-1} > 0 \}$ \\	 
	 \emph{(ii)} \textbf{complementarity condition} 
	 \begin{equation}
		\max | \blambda^i \odot (\bB \bU^i - \bD) |
         \end{equation}
	
	\item As additional indicators of accuracy, the iterative converged results, marked as $\cdot^{\mathrm{comp}}$, are systematically compared with the reference solution obtained using the standard saddle-point Lagrange multiplier method (implying the saddle-point matrix inversion), marked as~$\cdot^{\mathrm{SP}}$.
	Denoting by $\mathcal{Q}$ a quantity of interest, the relative error is defined as
	\begin{equation}
		\boldsymbol{e}_{\mathcal{Q}} =
		\frac{\| \mathcal{Q}^{\mathrm{SP}} - \mathcal{Q}^{\mathrm{comp}} \|_2}
		{\| \mathcal{Q}^{\mathrm{SP}} \|_2}.
	\end{equation}
The quantities of interest are typically the contact forces at convergence $\bF_C = - \bB^\top \blambda^\star$, and the final displacements $\bU^\star$, where the superscript~$\star$ denotes the converged values. 
\end{itemize}

Before presenting the numerical results for each test case, we first discuss in Section~\ref{sec:influence_proj} the importance of the fixed-point function selected for acceleration.

\subsection{Impact of the projection-step before acceleration} \label{sec:influence_proj}

We investigate here the effect of the chosen fixed-point function ($\mathcal{F}$ or $\mathcal{G}$, see Section~\ref{sec:accel-update-step}) on the different acceleration schemes. As highlighted in Algorithm~\ref{alg:algo_uni}, this involves whether a projection step is applied prior to acceleration.

By a slight abuse of notation, we denote the two possible strategies as follows:
\begin{itemize}
\item $\Pi \circ \mathcal{A} \circ \Pi$: acceleration applied to $\mathcal{F}$, see Eq.\eqref{eq:update_accel_F},
\item $\Pi \circ \mathcal{A}$: acceleration applied to $\mathcal{G}$, see  Eq.\eqref{eq:update_accel_G}.
\end{itemize}
{Table~\ref{tab:projection} summarizes the results obtained for all the test cases described in the following sections.} We focus on the impact of the projection-step position on the four acceleration strategies considered here.
We distinguish between so-called in-range and out-of-range values of the augmentation/penalty parameter. In-range values correspond to the range of parameter values for which the convergence of standard (non-accelerated) strategies is numerically observed. The upper bound of in-range values depends of the test case and has been given above (typically $10^4$ or $10^5$).

\begin{table}[h!]
	\centering
	\renewcommand{\arraystretch}{1.4}
	\setlength{\tabcolsep}{12pt}
	\scalebox{0.95}{
	\begin{tabular}{c|c|c||c|c|}
		\cline{2-5}
		& \multicolumn{2}{c||}{\textbf{In-range value}} 
		& \multicolumn{2}{c|}{\textbf{Out-of-range value}} \\
		\cline{2-5}
		& $\Pi \circ \mathcal{A} \circ \Pi $ &  $\Pi \circ \mathcal{A} $& $\Pi \circ \mathcal{A} \circ \Pi $ &  $\Pi \circ \mathcal{A}$ \\
		\hline
		\multicolumn{1}{|c|}{\textbf{FISTA + AR}}  
		& \cellcolor{mygray!50}\cmark & \cmark  & \xmark & \xmark \\
		\hline
		\multicolumn{1}{|c|}{\textbf{Anderson-1}}  
		& \cellcolor{mygray!50}\cmark & \xmark  & \xmark & \xmark \\
		\hline
		\multicolumn{1}{|c|}{\textbf{Anderson-1 + AR}}  
		& \cellcolor{mygray!50}\cmark & \xmark  & \xmark & \xmark \\
		\hline
		\multicolumn{1}{|c|}{\textbf{Crossed-Secant}}  
		& \cmark & \cellcolor{mygray!50}\cmark & \xmark & \cellcolor{mygray!50}\cmark \\
		\hline
	\end{tabular}
	}
	\caption{Influence of the projection step before acceleration on the performance of acceleration methods. The recommended strategy is highlighted in grey.}
	\label{tab:projection}
\end{table}

The reported behaviors are in good agreement with the theoretical insights given in Sections~\ref{sec:accel-update-step} and~\ref{s:acceleration}. For conciseness, the subsequent results focus on the best-performing fixed-point variants for each acceleration strategy, highlighted in grey in the Table~\ref{tab:projection}.
Hence, we retain the double-projection strategy $\Pi \circ \mathcal{A} \circ \Pi$ (i.e., acceleration on~$\mathcal{F}$) for the FISTA with adaptive restart and Anderson schemes, with or without adaptive restart, and the $\Pi \circ \mathcal{A}$ variant (i.e., acceleration on $\mathcal{G}$) for the Crossed-Secant scheme.  For the FISTA with AR, this leads to the best strategy in terms of the number of iterations, {cf. $i$ in Algorithm \ref{alg:algo_uni}}. For the Crossed-Secant method, the chosen approach ensures convergence for all augmentation/penalty parameter values. Note, that since the Crossed-Secant is a dynamic relaxation method, applying it to $\mathcal{F}$ (which includes a projection before acceleration) alters its intrinsic behavior, whereas applying it to $\mathcal{G}$ preserves the original dynamics.

\subsection{3D academic benchmark: Hertzian contact} \label{ss:Hertz}

\begin{figure}[htbp]
	\centering
	\begin{subfigure}{0.48\textwidth}
		\captionsetup[subfigure]{labelformat=empty}
		\centering
		\includegraphics[page=1, width=1\linewidth]{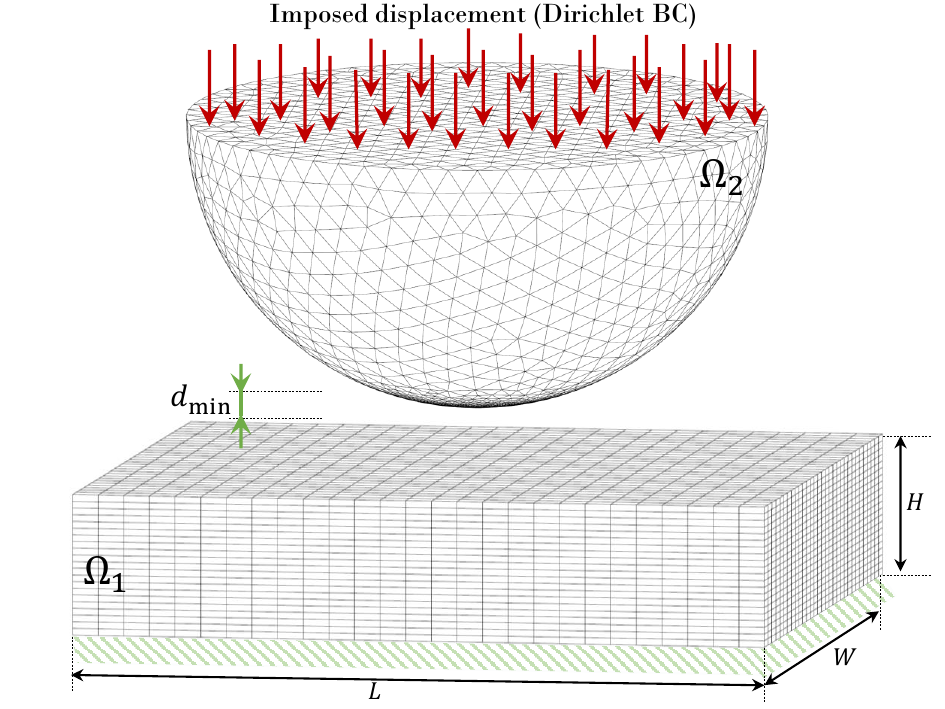}
		\caption{Problem setup}
		\label{fig:test_Hertz_pb}
	\end{subfigure}
	\hfill
	\begin{subfigure}{0.48\textwidth}
		\centering
		\includegraphics[page=2,width=\linewidth]{figures/test_sphere_plan_sketch_tet.pdf}
		\caption{Mesh of the half-sphere $\Omega_2$}
		\label{fig:test_Hertz_mesh}
	\end{subfigure}
	\caption{Hertzian contact problem -- computational domains $\Omega_1$ and $\Omega_2$ with boundary conditions.}
	\label{fig:test_Hertz}
\end{figure}

As the first numerical validation, we consider the classical Hertzian contact problem between a half-sphere and a plane, a well-established and widely used benchmark in computational contact mechanics. 
The test consists of a deformable elastic block $\Omega_1$ and a deformable elastic half-sphere $\Omega_2$, as shown in Figure~\ref{fig:test_Hertz}. 
The block $\Omega_1$ has length and width $L = W = 5 \cdot 10^{-2}\,\text{m}$ and thickness $H = 1 \cdot 10^{-2}\,\text{m}$, while the half-sphere $\Omega_2$ has radius $R = 2 \cdot 10^{-2}\,\text{m}$.
The minimal initial gap between the two bodies is set to $d_{\text{min}}= \min_{\Gamma^C} d = 0\,\text{m}$.
Both domains are modeled as isotropic, linear elastic materials. The Young’s modulus of the half-sphere is $E_{2} = 2.1 \cdot 10^{9}\,\text{Pa}$, whereas the block has a Young’s modulus of $E_{1} = 2.1 \cdot 10^{11}\,\text{Pa}$ (rigidity contrast of $100$). Identical Poisson’s ratios are assumed for both domains, with $\nu_{1} = \nu_{2} = 0.3$.
Contact is induced by prescribing a vertical displacement of $u_{\text{D}} = 3 \cdot 10^{-4}\,\text{m}$ on the equatorial plane of the half-sphere, thereby enforcing indentation into the block, which is fixed on its bottom surface.
Throughout the simulation, small deformations, linear elasticity, and frictionless normal contact are assumed.

In this test case, the block (domain $\Omega_1$) is discretized using trilinear hexahedral finite elements (or $\mathcal{Q}_1$-elements with 8 discretization nodes), while the sphere (domain $\Omega_2$) is discretized with linear tetrahedral elements ($\mathcal{P}_1$-elements with 4 nodes) with local refinement in the contact zone (see Figure~\ref{fig:test_Hertz_mesh}). The mesh of the block domain contains $19{,}712$ nodes ($17{,}577$ elements) while the sphere has $13{,}009$ nodes ($71{,}491$ elements).
As shown in Figure~\ref{fig:test_Hertz}, the meshes in domains are non-matching, leading to a node-to-surface pairing contact formulation.

\subsubsection{Application of accelerated Uzawa/Lagrange Multipliers formulation} \label{ss:sphere_plan_uzawa}

\noindent We first investigate the performance of the accelerated Uzawa iterative schemes. The solution process is carried out using Algorithm~\ref{alg:algo_uni}, with the dual contact variable update step from Eq.~\eqref{eq:uzawa_std_update}. 
\newline

\noindent \textbf{Standard Uzawa resolution and comparison of acceleration schemes} \\

Table~\ref{tab:hertz_uzawa_combined} summarizes the performance of the Uzawa schema and of the various acceleration strategies applied to it for augmentation parameters $\rho \in \{10^2,10^3, 10^4\}$. For these values of $\rho$, which are below the numerical convergence upper bound, all iterative strategies converge while satisfying the non-penetration and complementarity conditions with high precision. Moreover, the accuracy of the displacement and contact force solutions are in very good agreement with the reference saddle-point Lagrange multiplier method solutions.

\begin{table}[h!]
	\centering
	\renewcommand{\arraystretch}{1.5}
	\setlength{\tabcolsep}{12pt}
	\scalebox{0.7}{
		\begin{tabular}{|c|c|c|c|c|c|}
			\hline
			$\rho$ & \makecell{Nb\\iterations} 
			& \makecell{Relative contact\\force error} 
			& \makecell{Relative\\displacement error} 
			& \makecell{Gap\\function (m)} 
			& \makecell{Complementarity\\condition ($\mathrm{N\,m}$)} \\
			\hline
			\multicolumn{6}{|c|}{\textbf{Standard Uzawa method}} \\
			\hline
			$\rho = 10^{2}$ & 258745 &$ 3.19 \cdot 10^{-9}$  & $ 6.76 \cdot 10^{-11}$  & $ 5.20 \cdot 10^{-13}$ & $ 4.81 \cdot 10^{-12}$  \\ 
			$\rho = 10^{3}$ & 29995 & $3.20\cdot 10^{-10}$ & $6.78\cdot 10^{-12}$ & $5.53 \cdot 10^{-14}$ & $4.89\cdot 10^{-13}$   \\ 
			$\rho = 10^{4}$ & 3406 & $3.19\cdot 10^{-11}$ & $6.78\cdot 10^{-13}$ & $5.82\cdot 10^{-15}$ & $4.94\cdot 10^{-14}$ \\  
			\hline
			\multicolumn{6}{|c|}{\textbf{FISTA acceleration  scheme with adaptive restart}} \\
			\hline
			$\rho = 10^2$ & 3715 & $1.93 \cdot 10^{-9}$ & $6.36 \cdot 10^{-11}$ & $3.13 \cdot 10^{-13}$ & $6.09 \cdot 10^{-12}$ \\
			$\rho = 10^3$ & 1332  & $4.42 \cdot 10^{-11}$ & $1.82 \cdot 10^{-12}$ & $1.06 \cdot 10^{-14}$ & $1.09 \cdot 10^{-13}$ \\
			$\rho = 10^4$ & 398  & $3.55 \cdot 10^{-12}$ & $1.43 \cdot 10^{-13}$ & $7.81 \cdot 10^{-16}$ & $7.09 \cdot 10^{-15}$ \\
			\hline
			\multicolumn{6}{|c|}{\textbf{Anderson-1 acceleration scheme}} \\
			\hline
			$\rho = 10^2$ & 41528 & $1.30 \cdot 10^{-09}$ & $3.02 \cdot 10^{-11}$ & $2.40 \cdot 10^{-13}$ & $3.99 \cdot 10^{-12}$ \\
			$\rho = 10^3$ & 11695  & $5.60 \cdot 10^{-11}$ & $1.29 \cdot 10^{-12}$ & $1.13 \cdot 10^{-14}$ & $1.80\cdot 10^{-13}$ \\
			$\rho = 10^4$ & 1962  & $5.73 \cdot 10^{-12}$ & $5.73 \cdot 10^{-14}$ & $3.60 \cdot 10^{-16}$ & $9.12 \cdot 10^{-15}$ \\
			\hline
			\multicolumn{6}{|c|}{\textbf{Anderson-1 acceleration scheme with adaptive restart}} \\
			\hline
			$\rho = 10^{2}$ & 502 &$ 5.23 \cdot 10^{-10}$  & $6.87 \cdot 10^{-11}$  & $ 6.58 \cdot 10^{-14}$ & $ 1.54 \cdot 10^{-12}$  \\ 
			$\rho = 10^{3}$ & 342 & $4.80\cdot 10^{-11}$ & $5.60\cdot 10^{-12}$ & $1.36 \cdot 10^{-14}$ & $1.50\cdot 10^{-13}$   \\ 
			$\rho = 10^{4}$ & 227 & $3.98\cdot 10^{-12}$ & $3.14\cdot 10^{-13}$ & $1.11\cdot 10^{-15}$ & $1.34\cdot 10^{-14}$ \\  
			\hline
			\multicolumn{6}{|c|}{\textbf{Crossed-Secant acceleration scheme}} \\
			\hline
			$\rho = 10^{2}$ & 136 & $2.28\cdot 10^{-12}$ & $2.22 \cdot 10^{-13}$ & $1.27\cdot 10^{-15}$ & $1.00 \cdot 10^{-14}$   \\ 
			$\rho = 10^{3}$ & 132 & $3.11\cdot 10^{-12}$ & $9.97\cdot 10^{-14}$ & $7.92\cdot 10^{-16}$ & $5.03\cdot 10^{-15}$   \\ 
			$\rho = 10^{4}$ & 136 & $6.20\cdot 10^{-12}$ & $1.79\cdot 10^{-13}$ & $1.56 \cdot 10^{-15}$ & $1.00\cdot 10^{-14}$ \\  
			\hline
		\end{tabular}
	}
	\caption{Hertzian contact problem -- Uzawa/LM formulation -- performance of different acceleration schemes for the augmentation parameter $\rho \in \{10^2, 10^3, 10^4\}$}
	\label{tab:hertz_uzawa_combined}
\end{table}

A consistent trend is observed for the standard Uzawa method, as well as for its FISTA and Anderson-accelerated variants: across all error measures (displacement and contact force errors, as well as complementarity and gap values), each error decreases inversely proportionally to the parameter~$\rho$. For the standard Uzawa method, this observation is consistent with the chosen convergence criterion, defined in Equation~\eqref{eq:relative_residual}.
In practice it leads to the relation $\dfrac{\left\| \mathbf{g}^i \right\|_2}{\left\| \boldsymbol{\lambda}^i \right\|_2} \le \dfrac{\epsilon}{\rho}$ when the active zone does not vary anymore.
Consequently, the precision of the converged values is directly
related to the choice of the augmentation parameter~$\rho$ (for a
given convergence threshold $\epsilon$). Even though the acceleration
methods aim to modify the expression of $\blambda^i - \blambda^{i-1}$,
a dependence of the error measures on $\rho$ is observed for
the momentum-based methods (FISTA+AR and Anderson-1), but not for the
Crossed-Secant method.
Equation~\eqref{eq:CS_residual} shows that, for the Crossed-Secant method, the convergence relation becomes $\dfrac{\left\| \mathbf{g}^i \right\|_2}{\left\| \boldsymbol{\lambda}^i \right\|_2} \le \dfrac{\epsilon}{\left|\omega^i\right|\rho},$
with an automatically corrected augmentation parameter $\omega^i\rho$.

For the momentum-based acceleration schemes, the generic formulation
$ {\blambda}^{i} = \hat{\blambda}^{i} + {\beta^{i}}(\hat{\blambda}^{i} - \hat{\blambda}^{i-1})$
implies that the expression for ${\blambda}^{i} - {\blambda}^{i-1}$ consists of terms proportional to $\beta^i\rho$ together with the Uzawa term $\rho \mathbf{g}^i$. 
Consequently, the gap precision and other error measures retain a dependence on $\rho$.

\begin{figure}[htbp]
	\centering
	\begin{subfigure}{0.48\textwidth}
		\centering
		\includegraphics[width=\linewidth]{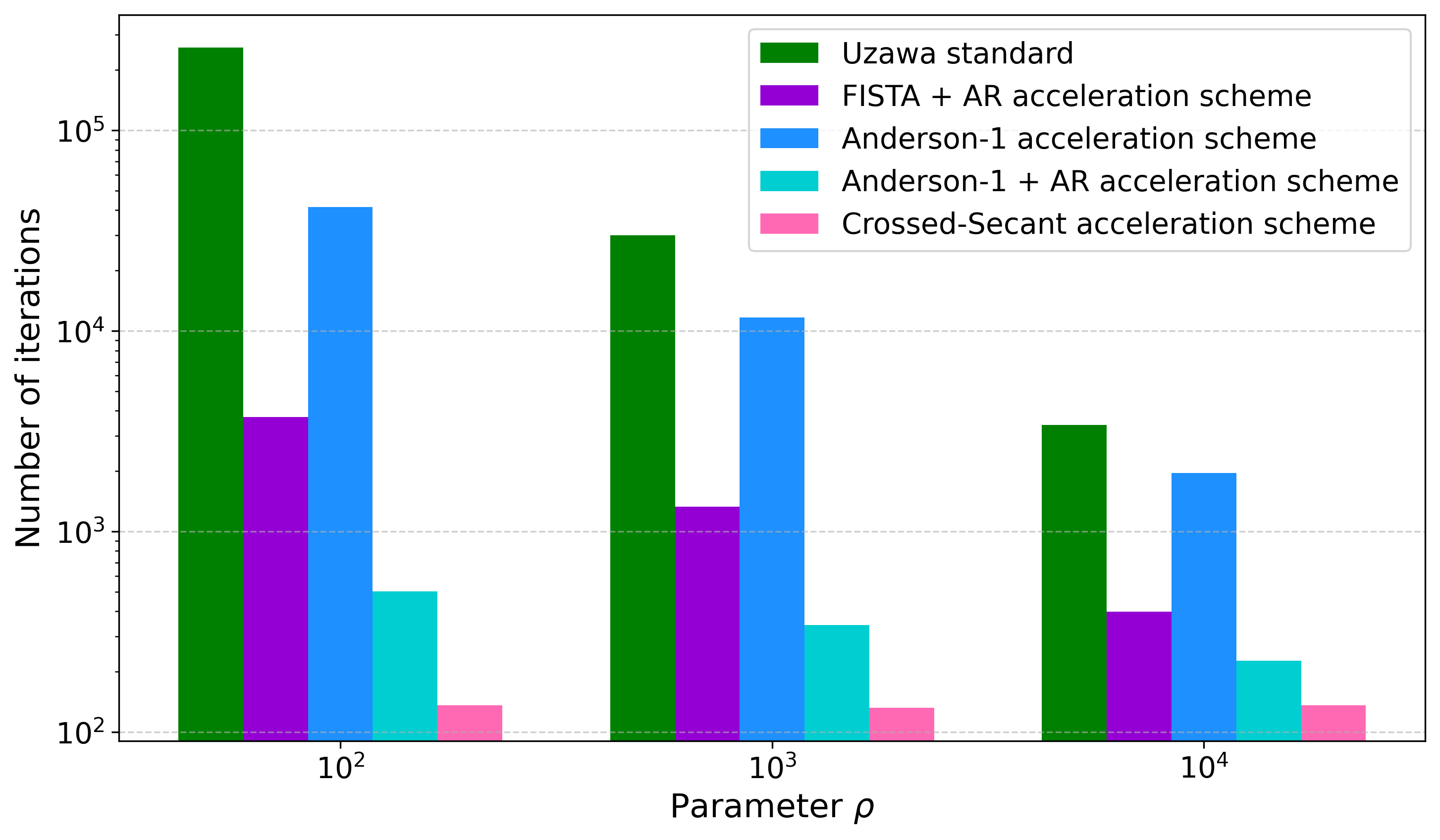}
		\caption{ Number of iterations for $\rho \in \{10^2, 10^3, 10^4\}$}
		\label{fig:sphere_plan_uzawa_iterations_diffrho}
	\end{subfigure}
	\hfill
	\begin{subfigure}{0.48\textwidth}
		\centering
		\includegraphics[width=\linewidth]{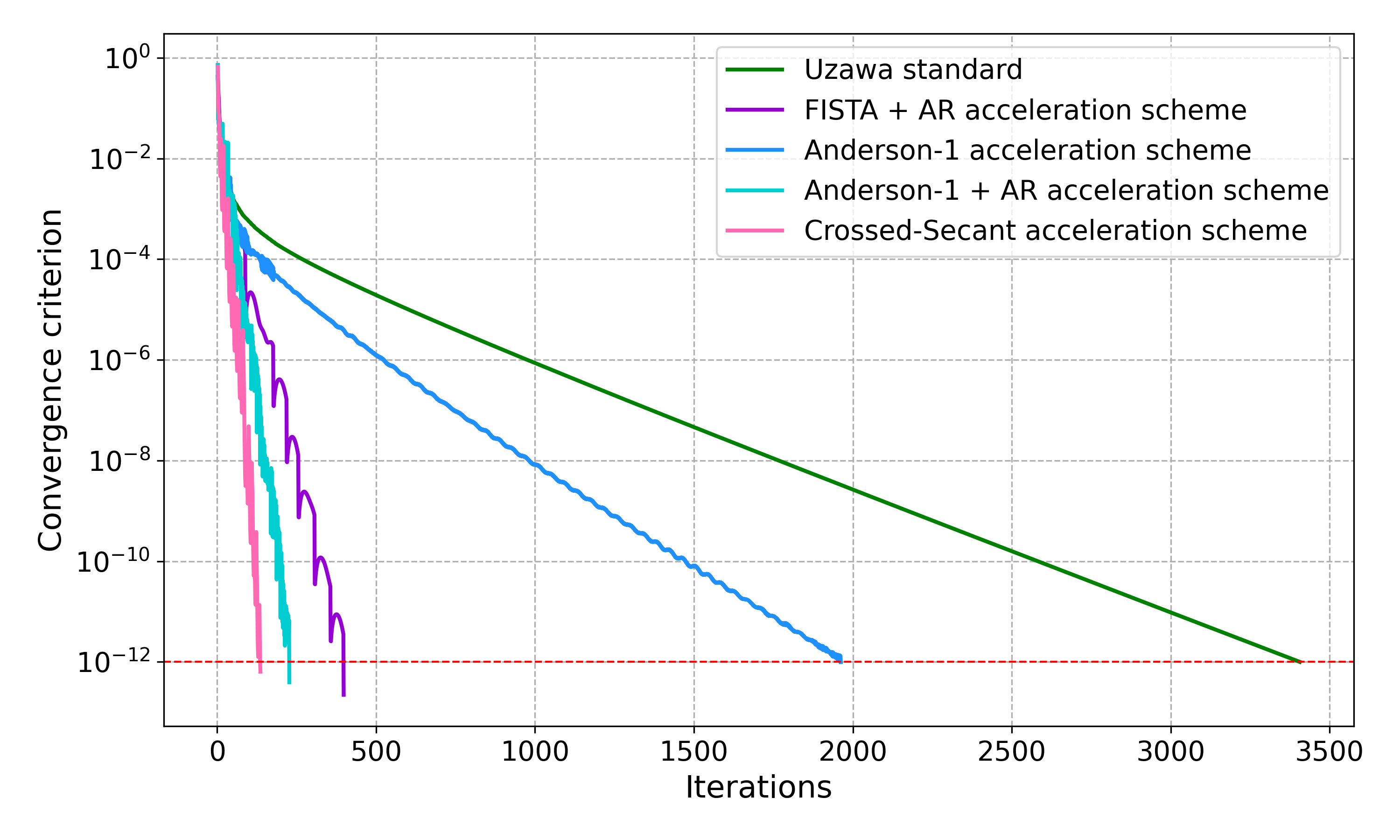}
		\caption{Conv. criterion for $\rho = 10^4$}
		\label{fig:sphere_plan_uzawa_residual_rho104}
	\end{subfigure}
	\caption{Hertzian contact problem -- Uzawa/LM formulation -- performance of different acceleration schemes }
	\label{fig:hertz_compar_uza}
\end{figure}

Figure~\ref{fig:sphere_plan_uzawa_iterations_diffrho} compares the number of iterations required by each acceleration scheme for the chosen values of $\rho$.
To compare the performance of the methods, we consider it sufficient to report only the number of iterations. The first iteration is the most computationally expensive, and its cost is nearly identical for all methods. In contrast, the subsequent iterations are relatively inexpensive, with a similar computational cost per iteration across all methods. Consequently, the total computational time can be approximated as the sum of the costs of the individual iterations, making the number of iterations a reliable indicator of the overall computational cost.

As it can be observed, the standard Uzawa method requires a very high number of iterations, particularly for small values of the augmentation parameter~$\rho$.
Introducing FISTA acceleration with adaptive restart significantly reduces the number of iterations, especially for large~$\rho$, while maintaining very low relative errors. The Anderson-1 acceleration scheme also improves convergence, and adding adaptive restart to Anderson-1 further decreases iterations count. 
 Finally, the Crossed-Secant acceleration scheme achieves the lowest number of iterations, with a number of iterations that is almost independent of $\rho$, while preserving excellent numerical accuracy. 
Hence, in practice, the Anderson-1 method with AR and the Crossed-Secant method enable the efficient use of the Uzawa method with low values of $\rho$.

Figure~\ref{fig:sphere_plan_uzawa_residual_rho104} reports the evolution of the convergence criterion ${r}^i$ (see Eq.~\eqref{eq:relative_residual}) as a function of the iteration count $i$ for $\rho = 10^4$. 
One can observe, that in contrast to the smooth and monotonically decreasing convergence profile of the standard Uzawa method, the accelerated variants exhibit oscillatory behavior with varying amplitudes.

To characterize the convergence behavior in a quantitative manner, we rely on an local indicator $p^i$ of the convergence order  
\begin{equation}
	p^i = \frac{\ln(\frac{{r}^{i+1}}{{r}^i})}{\ln(\frac{{r}^{i}}{{r}^{i-1}})}
	\label{eq:conv_rate}
\end{equation} 
In order to mitigate the influence of oscillations, the reported value of the convergence order $p$ in Table~\ref{tab:conv-order}, should be understood as a representative average trend rather than a strict asymptotic order.
\begin{table}[!ht]
		\centering
\scalebox{0.85}{
  \begin{tabular}{c|ccccc}
    & Uzawa standard & FISTA + AR & Anderson-1 & Anderson-1 + AR& Crossed-Secant\\
    \hline
    $p$  & 1.0 & 1.21 & 1.12 & 1.27 & 1.41
  \end{tabular}
}
  \caption{Approximate convergence order for the different strategies} \label{tab:conv-order}
\end{table}

As expected, the standard Uzawa method exhibits linear convergence ($p \simeq 1.0$). 
{All considered acceleration strategies exhibit superlinear convergence behavior, although Anderson-1 acceleration (without AR) only shows a weakly superlinear regime, with an observed convergence order of approximately 1.12.}
The Crossed-Secant strategy achieves the best convergence rate, with an estimated order of $1.41$, which is however lower than the expected theoretical order of $1.6$ for a Secant method.\\


\noindent \textbf{Performance of Crossed-Secant-accelerated Uzawa method} \\

This section further investigates the performance of the Crossed-Secant acceleration scheme and its robustness with respect to the augmentation parameter.

Table~\ref{tab:hertz_uzawa_cs_allrho} summarizes the results obtained for different values of the augmentation parameter spanning $\rho \in [10^1,10^{19}]$. 
We highlight here that the Crossed-Secant acceleration allows the use of augmentation parameter values far exceeding the restricted bounds typically acceptable for the standard Uzawa method and other acceleration techniques, as discussed previously.
Moreover, in contrast to the standard Uzawa method, where small values of~$\rho$ lead to slow convergence and a slight progressive degradation in accuracy, the proposed approach achieves rapid convergence even for very small values of~$\rho$, without any noticeable loss of precision.

\begin{table}[h!]
	\centering
	\renewcommand{\arraystretch}{1.3} 
	\scalebox{0.75}{
		\begin{tabular}{|c|c|c|c|c|c|}
			\hline
			$\rho$ 
			& \makecell{Nb\\iterations} 
			& \makecell{Relative contact\\force error} 
			& \makecell{Relative\\displacement error} 
			& \makecell{Gap\\function (m)} 
			& \makecell{Complementarity\\condition ($\mathrm{N\,m}$)} \\
			\hline
			$\rho = 10^{1}$ & 155 & $3.76 \cdot 10^{-12}$ & $2.95 \cdot 10^{-13}$ & $6.46 \cdot 10^{-16}$ & $7.90\cdot 10^{-15}$  \\ \hline
			$\rho = 10^{3}$ & 132 & $3.11\cdot 10^{-12}$ & $9.97\cdot 10^{-14}$ & $7.92\cdot 10^{-16}$ & $5.03\cdot 10^{-15}$ \\ \hline
			$\rho = 10^{5}$ & 145 & $3.46\cdot 10^{-12}$ & $2.66\cdot 10^{-13}$ & $7.22 \cdot 10^{-16}$ & $9.57\cdot 10^{-15}$ \\ \hline
			$\rho = 10^{7}$ & 156 & $3.99\cdot 10^{-12}$ & $2.46 \cdot 10^{-13}$ & $8.96 \cdot 10^{-16}$ & $1.07\cdot 10^{-14}$  \\ \hline
			$\rho = 10^{9}$ & 172 & $4.05 \cdot 10^{-13}$ & $2.10 \cdot 10^{-13}$ & $8.80 \cdot 10^{-17}$ & $2.80 \cdot 10^{-15}$  \\ \hline
			$\rho = 10^{11}$ & 179 & $2.80 \cdot 10^{-12}$ & $9.32 \cdot 10^{-14}$ & $6.05 \cdot 10^{-16}$ & $5.46 \cdot 10^{-15}$\\ \hline
			$\rho = 10^{13}$ & 187 & $7.35 \cdot 10^{-12}$ & $1.66\cdot 10^{-13}$ & $1.96 \cdot 10^{-15}$ & $2.64\cdot 10^{-14}$ \\ \hline
			$\rho = 10^{15}$ & 233 & $1.20 \cdot 10^{-12}$ & $6.83 \cdot 10^{-14}$ & $2.80\cdot 10^{-15}$ & $2.79 \cdot 10^{-16}$  \\ \hline
			$\rho = 10^{17}$ & 211 & $4.50\cdot 10^{-12}$ & $2.54 \cdot 10^{-13}$ & $1.45 \cdot 10^{-15}$ & $9.22 \cdot 10^{-15}$\\ \hline
			$\rho = 10^{19}$ & 274 & $2.47 \cdot 10^{-13}$ & $8.82 \cdot 10^{-14}$ & $1.13 \cdot 10^{-16}$ & $1.67 \cdot 10^{-15}$ \\ \hline 
		\end{tabular}
	}
	\caption{Hertzian contact problem -- Uzawa/Lagrange Multipliers formulation -- performance of the Crossed-Secant acceleration scheme for the augmentation parameter $\rho \in [10^1, 10^{19}]$}
	\label{tab:hertz_uzawa_cs_allrho}
\end{table}
      
\begin{figure}[htbp]
	\centering
	\begin{subfigure}{0.48\textwidth}
		\centering
		\includegraphics[width=\linewidth]{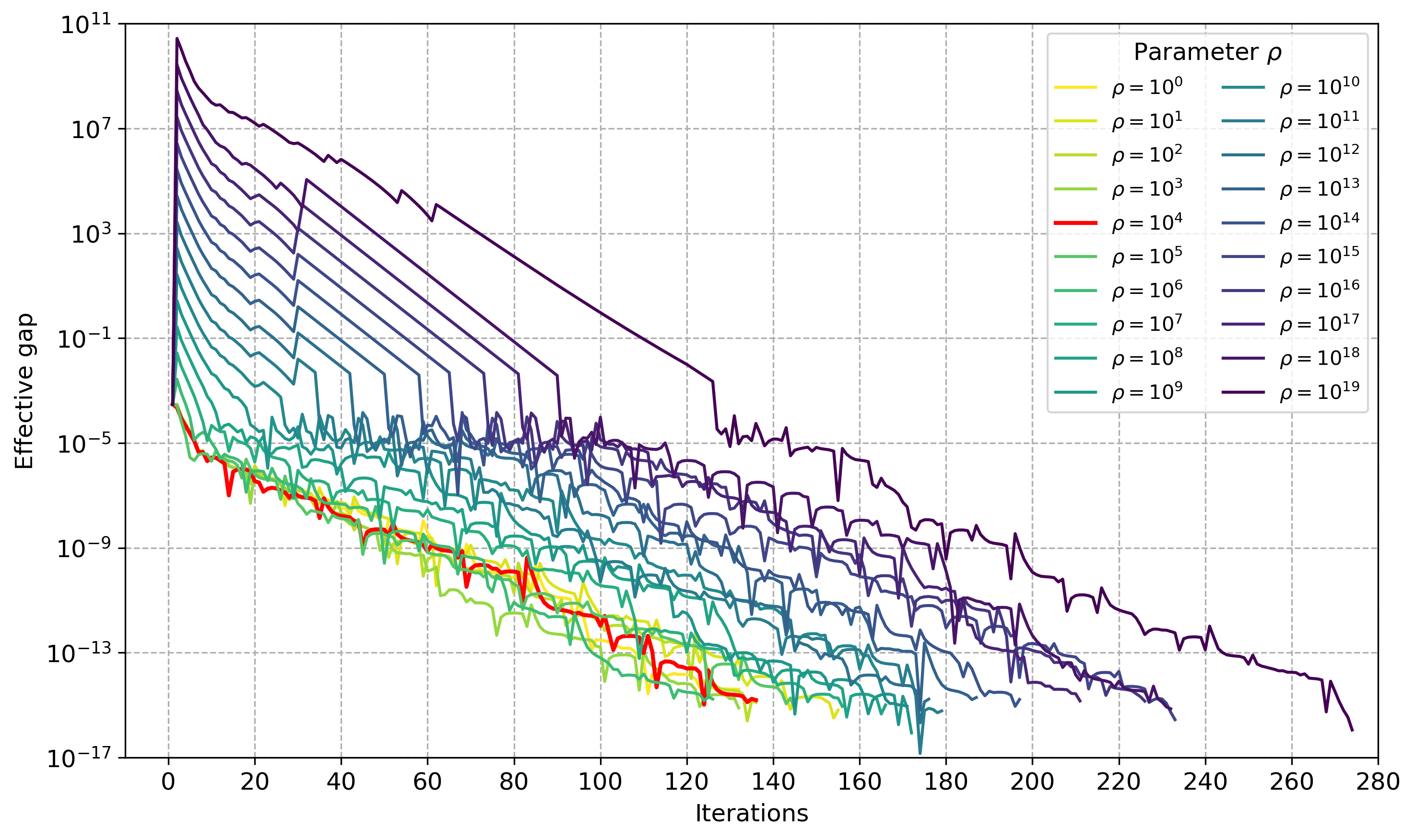}
		\caption{Effective gap for  $\rho \in [1, 10^{19}]$}
		\label{fig:sphere_plan_uzawa_accelCS_DN}
	\end{subfigure}
	\hfill
	\begin{subfigure}{0.48\textwidth}
		\centering
		\includegraphics[width=\linewidth]{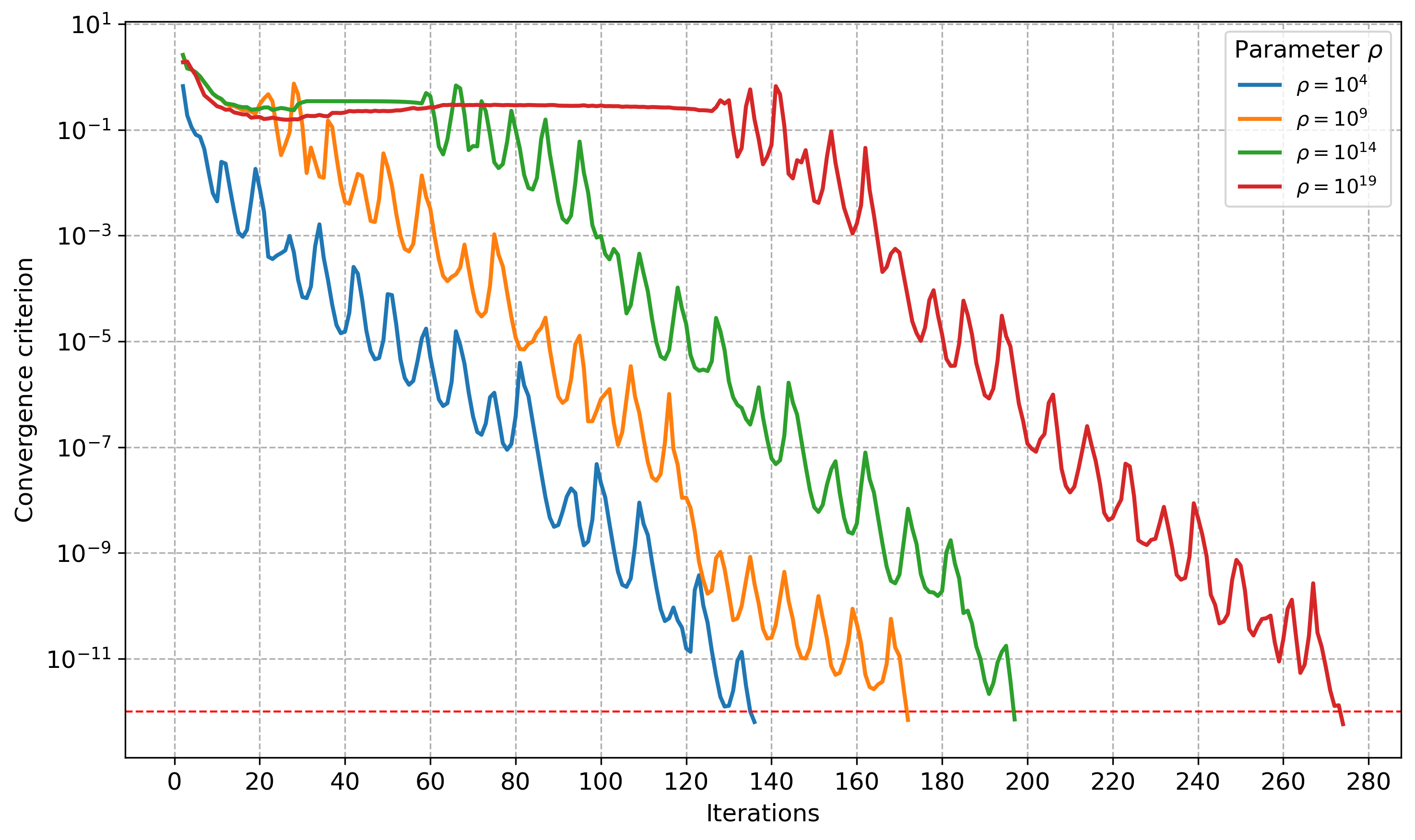}
		\caption{Conv. criterion for $\rho \in \{10^4, 10^{9}, 10^{14}, 10^{19}\}$}
		\label{fig:sphere_plan_uzawa_accelCS_RSDnrmEu}
	\end{subfigure}
	\caption{Hertzian contact problem -- Uzawa/Lagrange Multipliers formulation -- performance of the Crossed-Secant acceleration scheme for different augmentation parameters}
	\label{fig:hertz_accel_cs_uza}
\end{figure}

Figure~\ref{fig:sphere_plan_uzawa_accelCS_DN} shows the evolution of the effective gap for all parameter values. The initial gap ${d}_{\text{min}} = 0$ m (iteration $0$) is omitted from the figure because of the logarithmic scale.
One can observe that for augmentation parameters $\rho \leq 10^4$ (the maximum admissible value for the standard Uzawa method, highlighted in red in Figure~\ref{fig:sphere_plan_uzawa_accelCS_DN}), the effective gap decreases from the first iterations.
By contrast, for larger values of $\rho$, the effective gap initially
increases due to the large $\blambda$ values computed with these
$\rho$ values (see Eq.~\eqref{eq:uzawa_std_update}). Consequently,
large $\rho$ values may lead to iterations where the solids are no
longer in contact and are even significantly separated. The standard
Uzawa method and momentum-based acceleration techniques are
ineffective in this case and may even diverge. 
{In practice, these methods tend only to sharpen the contact zone, starting from a large initial set of nodes.}
{The Crossed-Secant acceleration scheme, in contrast, can handle scenarios in which some iterates exhibit temporary separation of the contacting domains.}
This feature allows the Crossed-Secant scheme to achieve convergence even for large, out-of-range parameter values.

Figure~\ref{fig:sphere_plan_uzawa_accelCS_RSDnrmEu} illustrates the
evolution of the relative iterate convergence criterion computed using
Eq.~\eqref{eq:relative_residual}.
For values of $\rho$ exceeding the guaranteed convergence upper bound, a characteristic feature is the presence of an initial plateau, whose duration increases with~$\rho$. 
This plateau corresponds to iterations where the active contact zone varies significantly, due to over-estimated contact forces (through the large parameter value), which may temporarily induces empty contact zones. Once the active contact zone has converged, the algorithm naturally transitions to the effective convergence regime (the projection reduces in this case to the identity operator). To the best of our knowledge, the algorithm manages to handle this scenario thanks to the particular properties of the Crossed-Secant method, demonstrated in Section~\ref{ss:accel_cs}.

As a consequence, the total number of iterations required to reach convergence also increases with~$\rho$.
One observes that the convergence behavior is qualitatively similar in the decreasing portion of the curves and consistent with the rates reported above for~$\rho=10^4$. 

Overall, as our results suggest, coupling the Uzawa method with the Crossed-Secant acceleration scheme applied to the fixed-point $\mathcal{G}$ (without projection prior to acceleration) greatly enhances Uzawa’s practical flexibility, allowing the use of a broad range of~$\rho$ values that were previously inaccessible, while reaching great levels of accuracy across all error measures.

\subsubsection{Application of penalty-splitted formulation}
\label{ss:sphere_plan_penalty} 

\begin{figure}[htbp]
	\centering
	\begin{subfigure}{0.48\textwidth}
		\centering
		\includegraphics[width=\linewidth]{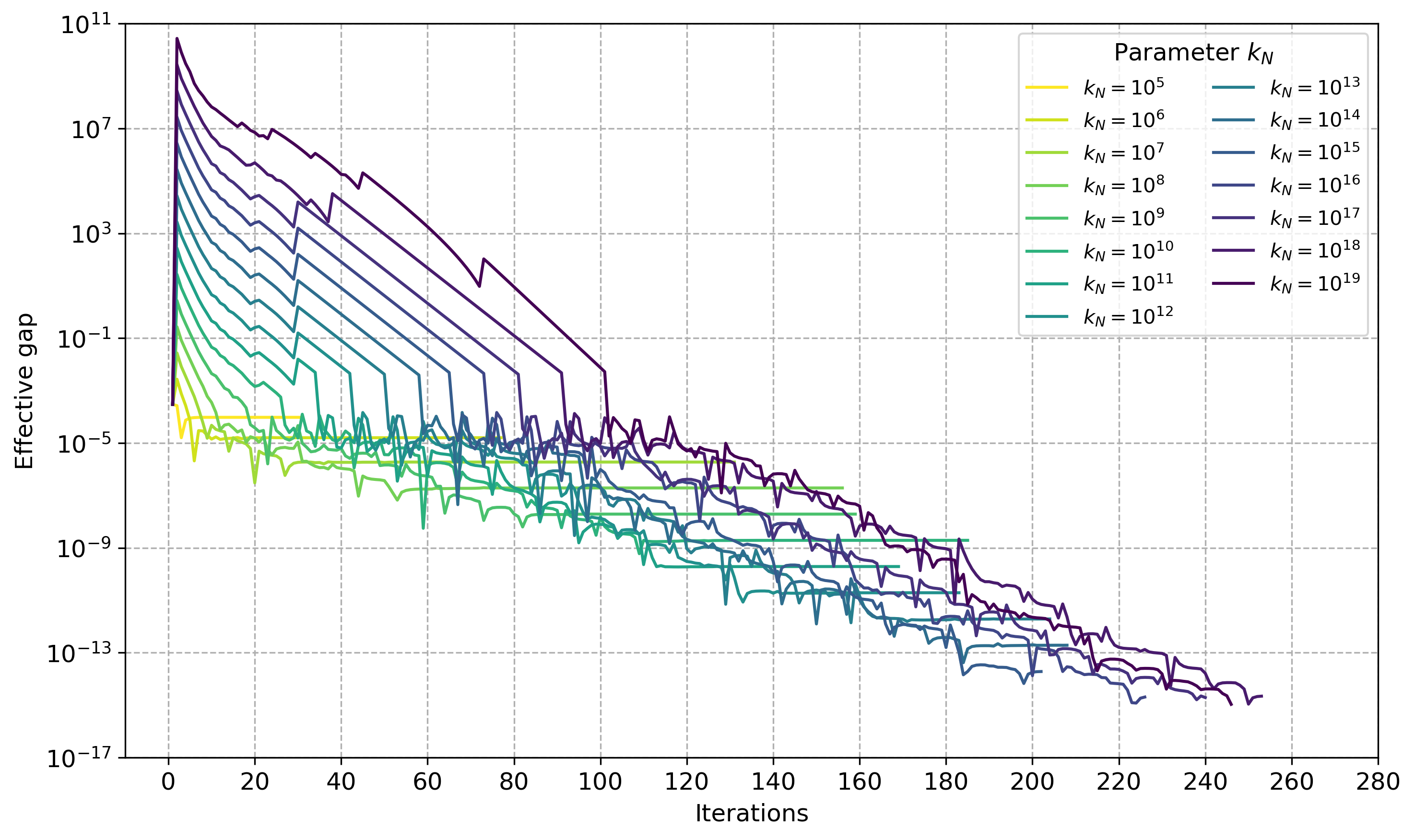}
		\caption{Effective gap for $k_N \in [10^5, 10^{19}]$}
		\label{fig:sphere_plan_penalty_accelCS_DN}
	\end{subfigure}
	\hfill
	\begin{subfigure}{0.48\textwidth}
		\centering
		\includegraphics[width=\linewidth]{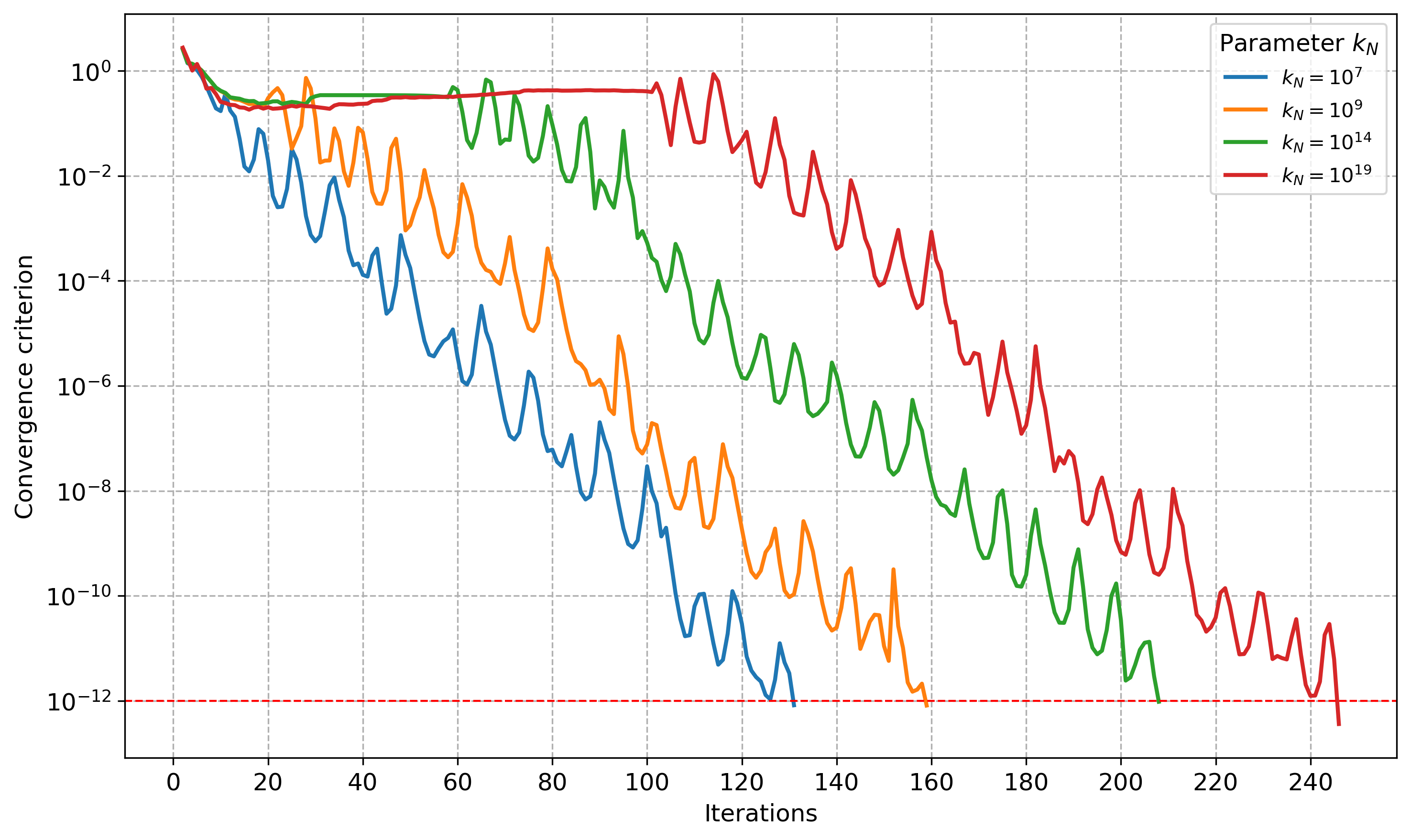}
		\caption{Conv. criterion for $k_N \in \{10^7, 10^{9}, 10^{14}, 10^{19}\}$}
		\label{fig:sphere_plan_penalty_accelCS_RSDnrmEu}
	\end{subfigure}
	\caption{Hertzian contact problem -- penalty-splitted fixed-point  formulation -- performance of the Crossed-Secant acceleration scheme for different penalty parameters }
	\label{fig:hertz_accel_cs_penalty}
\end{figure}

This section is dedicated to the evaluation of the proposed penalty-splitted fixed-point iterative strategy with the dual contact variable update step defined in Eq.~\eqref{eq:penalty_fp_update}.
For consistency, we also applied the standard non-accelerated penalty fixed-point scheme, as well as the FISTA+AR and Anderson acceleration techniques.
While these methods do converge, their applicability is very limited due to the constraints on the penalty parameter, which, combined with the approximate nature of the penalty formulation, lead to low-accuracy results. Only the Crossed-Secant acceleration scheme enables us to achieve the desired level of accuracy by exceeding the classical convergence upper bounds on the penalty parameter.
Therefore, only the results obtained with the Crossed-Secant acceleration scheme are reported here.  Table~\ref{tab:hertz_penalty_cs} demonstrates the results in terms of the accuracy measure over a wide range of penalty coefficients \(k_N \in [10^{5}, 10^{19}]\), while  Figure~\ref{fig:hertz_accel_cs_penalty} reports the evolution of the effective gap and the convergence criterion.

\begin{table}[h!]
	\centering
	\renewcommand{\arraystretch}{1.3} 
	\scalebox{0.75}{
		\begin{tabular}{|c|c|c|c|c|c|}
			\hline
			$k_N$ 
			& \makecell{Nb\\iterations} 
			& \makecell{Relative contact\\force error} 
			& \makecell{Relative\\displacement error} 
			& \makecell{Gap\\function (m)} 
			& \makecell{Complementarity\\condition ($\mathrm{N\,m}$)} \\
			\hline
			$k_N  = 10^{5}$ & 31 & $3.24\cdot 10^{-1}$ & $9.69 \cdot 10^{-2}$ & $9.58\cdot 10^{-5}$ & $9.18 \cdot 10^{-4} $ \\ \hline
			$k_N  = 10^{7}$ & 131 & $5.87\cdot 10^{-3}$ & $1.40 \cdot 10^{-3}$ & $1.88\cdot 10^{-6}$ & $3.55 \cdot 10^{-5} $ \\ \hline
			$k_N  = 10^{9}$ & 159 & $6.02\cdot 10^{-5}$ & $1.43 \cdot 10^{-5}$ & $1.94\cdot 10^{-8}$ & $3.79\cdot 10^{-7}$  \\ \hline
			$k_N  = 10^{11}$ & 169 & $6.01 \cdot 10^{-7}$ & $1.41 \cdot 10^{-7}$ & $1.94 \cdot 10^{-10}$ & $3.79\cdot 10^{-9}$ \\ \hline
			$k_N  = 10^{13}$ & 204 & $6.03\cdot 10^{-9}$ & $1.40 \cdot 10^{-9}$ & $1.94 \cdot 10^{-12}$ & $3.79\cdot 10^{-11}$ \\ \hline
			$k_N  = 10^{15}$ & 202 & $6.08\cdot 10^{-11}$ & $1.39 \cdot 10^{-11}$ & $1.92 \cdot 10^{-14}$ & $3.74 \cdot 10^{-13}$ \\ \hline
			$k_N  = 10^{17}$ & 240 & $7.38 \cdot 10^{-12}$ & $3.77 \cdot 10^{-13}$ & $1.95 \cdot 10^{-15}$ & $1.84 \cdot 10^{-14}$ \\ \hline
			$k_N  = 10^{19}$ & 246 & $2.87\cdot 10^{-12}$ & $1.11 \cdot 10^{-13}$ & $1.07 \cdot 10^{-15}$ & $4.12 \cdot 10^{-15}$ \\ \hline
			
		\end{tabular}
	}
	\caption{Hertzian contact problem -- penalty-splitted fixed-point formulation -- performance of the Crossed-Secant acceleration scheme for the penalty parameter $k_N \in [10^5, 10^{19}]$}
	\label{tab:hertz_penalty_cs}
\end{table}

In agreement with the theoretical properties of penalty formulation, it can be observed, that the final accuracy is strongly related to the value of the penalty parameter. For smal~$k_N$ values, the method converges more rapidly, but the associated error measures remain relatively large. As $k_N$ increases, all error measures decrease by several orders of magnitude, illustrating the asymptotic consistency of the penalty approach, for which exact enforcement of the contact constraint is achieved only in the limit $k_N \to \infty$, see for instance~\cite{Epalle2025,Chouly2013Penalty,Nour_1987}. As illustrated in Figure~\ref{fig:sphere_plan_penalty_accelCS_DN}, for smal penalty parameters, the effective gap exhibits a saturation behavior at the end of the iterative process, corresponding to the intrinsic nature of the penalty formulation. In this regime, the contact constraint is enforced only approximately, and the final gap value is governed by the penalty parameter, remaining proportional to $1/k_N$.

Figure~\ref{fig:sphere_plan_penalty_accelCS_RSDnrmEu} illustrates the convergence behavior of the iterative penalty-splitted approach, as measured by the relative iterate convergence criterion defined in Eq.~\eqref{eq:relative_residual}. A convergence pattern very similar to that of the Uzawa formulation is observed, confirming our theoretical insight in Section~\ref{s:penalty}. In particular, the decreasing portion of the convergence curves is essentially independent of the penalty parameter, while an initial plateau with a duration related to $k_N$ is observed, in close analogy with the behavior reported for the Uzawa formulation (see Figure~\ref{fig:sphere_plan_uzawa_accelCS_RSDnrmEu}).

With classical penalty solvers, the ill-conditioning of the system matrix typically restricts the usable range of~$k_N$, even when using efficient preconditioners (see~\cite{Epalle2025}), making the use of large values infeasible in practice. In contrast, coupling our operator-splitting iterative fixed-point formulation with the Crossed-Secant acceleration scheme enables the robust solution of highly penalized problems, achieving near machine-precision accuracy comparable to that obtained with a Lagrange multiplier formulation. No specialized preconditioner is required. Overall, the proposed method preserves the simplicity of the penalty framework while substantially improving robustness and accuracy, allowing the use of very large~$k_N$ values and achieving precision levels that classical penalty solvers generally cannot reach.

\subsubsection{Validation against the Hertz's analytical solution} \label{ss:sphere_plan_sol_analyt}

To validate the numerical results obtained with the Crossed-Secant method, a comparison is made with the Hertz's contact theory.
The original Hertz problem, in which a sphere is loaded with a concentrated normal force on its top, has been replaced here by the equivalent scenario, in which a displacement is imposed in the equatorial plane of the sphere (cf. Section \ref{ss:Hertz}) and the corresponding resultant normal reaction $F$ is used for the analytical solution. This approach has
been widely adopted in previous studies and has been shown to yield sufficiently accurate results for the evaluation of the analytical contact stresses. \\
The analytical Hertz solution gives a circular contact area of radius $a = \sqrt[3]{\dfrac{3 F R}{4 E^\ast}}$, and a maximum contact pressure $P_{\rm max} = \dfrac{3 F}{2 \pi a^2}$,
where $R$ is the sphere radius and $E^\ast$ is the effective Young's modulus defined as $
\dfrac{1}{E^\ast} = \dfrac{1-\nu_1^2}{E_1} + \dfrac{1-\nu_2^2}{E_2}$
with $E_1$, $E_2$ and $\nu_1$, $\nu_2$ being the Young's moduli and
Poisson ratios of the plane and sphere, respectively. In the present
case, the computed contact area radius and maximum pressure are $a \simeq 2.5401 \cdot 10^{-3} $ m and $P_{\rm max} \simeq 0.1847 $ GPa, respectively.
The Hertz contact pressure is given by $P = P_{\rm max} \, \sqrt{1 - \left( \dfrac{r}{a} \right)^2 }, \quad 0 \le r \le a$, where $r$ is the distance from the center of the contact area.

\begin{figure}[htbp]
	\centering
	\begin{subfigure}{0.48\textwidth}
		\centering
		\includegraphics[width=\linewidth]{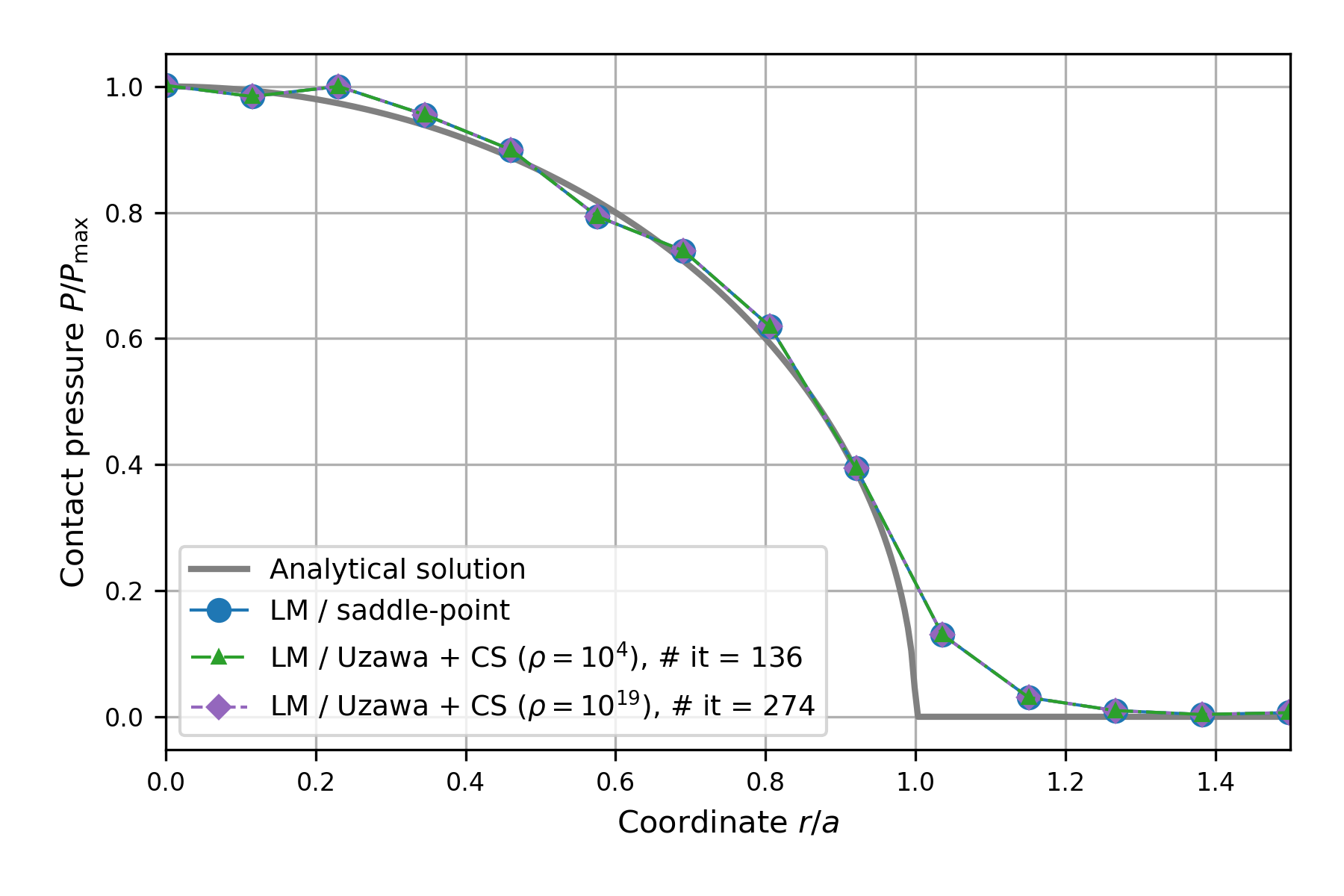}
		\caption{LM/Uzawa formulation + Crossed-Secant  (CS)  acceleration method for the augmentation parameter $\rho \in \{10^4, 10^{19}\}$}
		\label{fig:sphere_plan_sol_analyt_uza}
	\end{subfigure}
	\hfill
	\begin{subfigure}{0.48\textwidth}
		\centering
		\includegraphics[width=\linewidth]{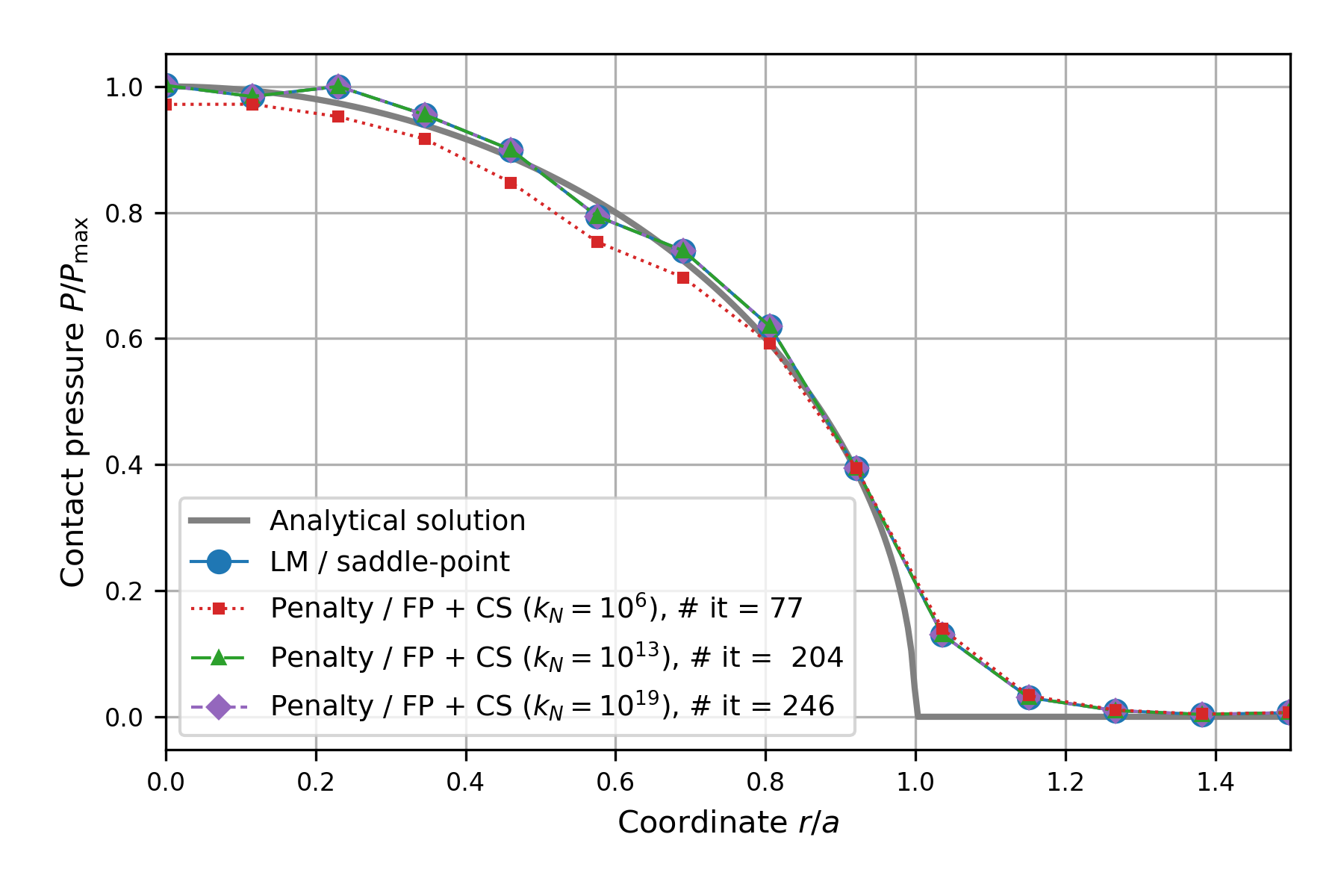}
		\caption{Penalty/fixed-point (FP) formulation + Crossed-Secant (CS) acceleration method for the penalty parameter $k_N \in \{10^6, 10^{13}, 10^{19}\}$}
		\label{fig:sphere_plan_sol_analyt_pen}
	\end{subfigure}
	\caption{Hertzian contact problem -- validation against the Hertzian analytical solution }
	\label{fig:sphere_plan_sol_analyt}
\end{figure}
	
Figure~\ref{fig:sphere_plan_sol_analyt} shows the dimensionless
contact pressure $P/ P_{\rm max}$ plotted against the dimensionless
radial coordinate $r/a$. Figures~\ref{fig:sphere_plan_sol_analyt_uza}
and~\ref{fig:sphere_plan_sol_analyt_pen} report the results obtained
using the Crossed-Secant--accelerated Uzawa and the penalty fixed-point iterative methods, detailed in Sections~\ref{ss:sphere_plan_uzawa} and~\ref{ss:sphere_plan_penalty}, respectively. In both cases the LM/saddle-point solution, considered
as the numerical reference solution, is also plotted. A
two-dimensional representation is presented, taking advantage of the
symmetry of the problem. Overall, the agreement between the numerical
and analytical solutions is very good, although small discrepancies
are observed near the boundary of the contact region due to the mesh
resolution. Moreover, the iterative solutions are virtually indistinguishable from the reference Lagrange multiplier solutions, except when small penalty parameters are used.

\subsection{3D industrial problem: thermal induced contact} \label{ss:IPG}

The second example is an industrial problem arising from nuclear fuel rod behavior simulation. The modeled phenomenon, known as Pellet-Cladding Interaction (PCI), refers to a set of physical phenomena occurring in Pressurized Water Reactors (PWRs) during irradiation, see for example~\cite{INTROINI2024,MICHEL2008,sercombe2020}. 

\begin{figure}[htbp]
	\centering
		\includegraphics[width=1\linewidth]{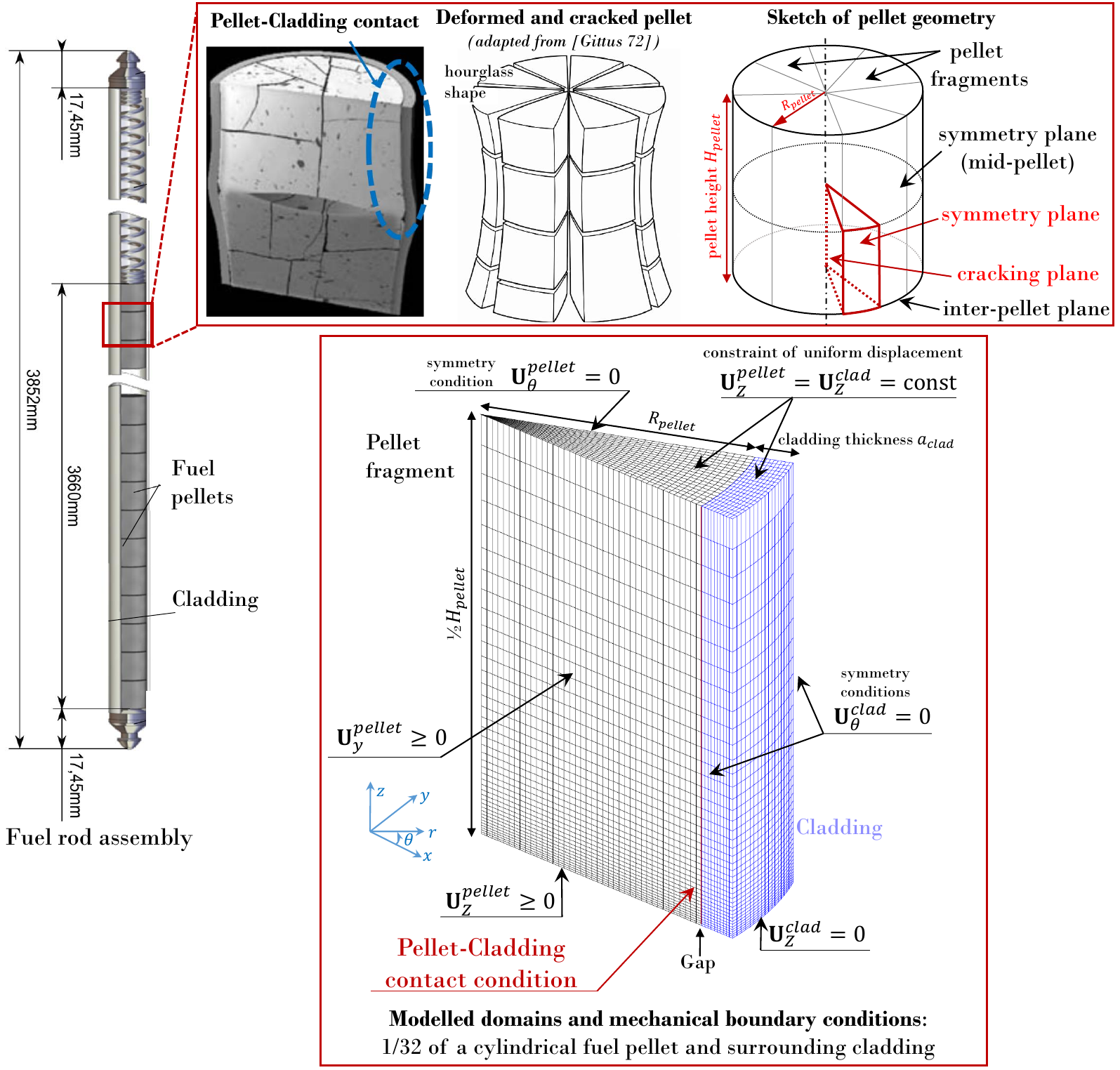}
		\caption{3D industrial problem: thermal-induced Pellet-Cladding Interaction -- fuel rod assembly, example of pellet-cladding contact, deformed and cracked pellet (adapted from [Gittus 1972, \cite{Gittus72}]), computational domain with mechanical boundary conditions}
	\label{fig:test_IPG_setup}
\end{figure}

Figure \ref{fig:test_IPG_setup}-left illustrates the whole fuel rod composed of cylindrical pellets enclosed in a metallic cladding tube. 
The behavior of the fuel rod is primarily driven by high temperature and irradiation. Heat generation induces fuel pellets thermal expansion and cracking. 
The strong thermal gradient between the center and the periphery of the pellet causes differential thermal expansion, leading the pellet fragments to adopt a so-called hourglass shape, see Figure~\ref{fig:test_IPG_setup}-top. Combined with other irradiation phenomena reducing the initial pellet--cladding gap, the thermal strain results in the closure of the gap and in the contact between the pellet fragments and the cladding. 
The local strains applied to the cladding produce a characteristic “bamboo-like” shape at the end of irradiation.

{As done in legacy fuel performance codes (see for example~\cite{INTROINI2024}), assuming a pre-cracked fuel pellet and exploiting the resulting symmetries, the analysis can be restricted to a representative fragment corresponding to 1/32 of the fuel pellet, including the surrounding cladding. The upper boundary of the fragment coincides with the mid-pellet plane, whereas the lower boundary coincides with the inter-pellet plane.}
{To facilitate the reproducibility of this test case, the industrial geometry was simplified while preserving the key features governing pellet-cladding mechanical interaction. The fuel pellet is idealized as a perfect cylinder and assumed to be pre-cracked into eight identical radial fragments. Furthermore, the temperature field is provided as an input to the mechanical analysis (see below).}

The geometry and mechanical boundary conditions are shown in Figure~\ref{fig:test_IPG_setup}-bottom. The cylindrical fuel pellet has a radius of $R_{\mathrm{pellet}} = 4.1$ mm and a height of $H_{\mathrm{pellet}} = 13.8$ mm. Owing to axial symmetry, only the lower half of the pellet, corresponding to a height of $H_{\mathrm{pellet}}/2$, is modeled.
In addition, exploiting the symmetry associated with the pre-cracked geometry, only a 22.5\textdegree{} sector is considered. The surrounding cladding segment has a thickness of $a_{\mathrm{clad}} = 0.57$ mm. A uniform initial pellet--cladding gap $d = 2\,\mu\mathrm{m}$ is considered in this static test case.

In this study, both the fuel and the cladding are modeled as isotropic linear elastic solids. Mechanical properties are defined as follows: for the fuel, Young's modulus $E_{\mathrm{pellet}} = 1.9 \cdot 10^{11}$ Pa and Poisson's ratio $\nu_{\mathrm{pellet}} = 0.3$; for the cladding, $E_{\mathrm{clad}} = 7.8 \cdot 10^{10}$ Pa and $\nu_{\mathrm{clad}} = 0.34$.

\begin{figure}[htbp]
	\centering
	
	\begin{subfigure}{0.48\textwidth}
		\centering
		\includegraphics[width=\linewidth]{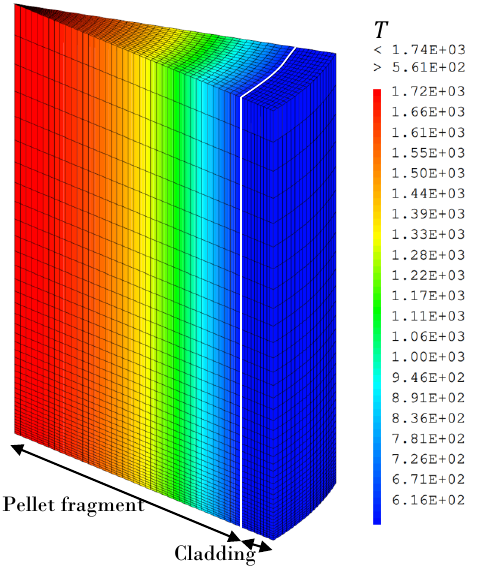}
		\caption{Imposed temperature field}
		\label{fig:IPG_thermalBC}
	\end{subfigure}
	\hfill
	\begin{subfigure}{0.48\textwidth}
		\centering
		\includegraphics[width=\linewidth]{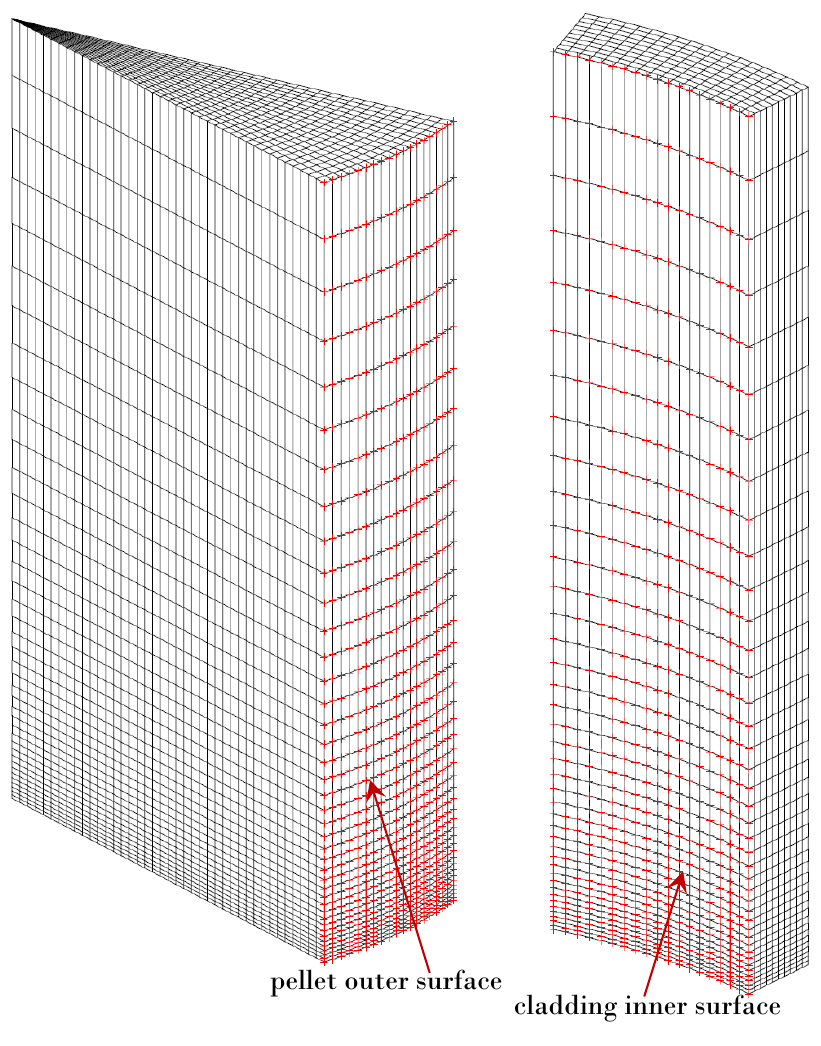}
		\caption{Active contact nodes}
		\label{fig:test_IPG_contact_zones}
	\end{subfigure}
	
	\caption{Thermal-induced contact configuration.}
	\label{fig:test_IPG_contact}
\end{figure}

{As in a previous 2D benchmark~\cite{LRL_CMAME_2017}, a radial temperature field is prescribed in the pellet, while a uniform temperature field is prescribed in the cladding, see Eq. \eqref{eq:champ-temp}. This temperature distribution is representative of typical operating conditions and illustrated in Figure~\ref{fig:IPG_thermalBC}.}
	\begin{equation} \label{eq:champ-temp}
		\begin{aligned}
			T(r) &= 7.0 \times 10^7 \left(R_{\mathrm{pellet}}^2-r^2\right)+561 \ \mathrm{K}
			&& \text{in the pellet}, \\
			T &= 561 \ \mathrm{K}
			&& \text{in the cladding}.
                      \end{aligned}
	\end{equation}

{In the mechanical analysis, thermal expansion is accounted for through the thermal strain $\boldsymbol{\varepsilon}_{\mathrm{th}} = \alpha_T (T - T_{ref})\,\mathbf{Id}$, where the coefficient of thermal expansion is taken as $\alpha_T = 10^{-5}\,\mathrm{K}^{-1}$ for the fuel and $\alpha_T = 6 \times 10^{-6}\,\mathrm{K}^{-1}$ for the cladding and the reference temperature is chosen to be $T_{ref}= 293 \,\mathrm{K}$.
The resulting differential thermal expansion induces contact between the pellet and the cladding.}

{To accurately capture the characteristic hourglass-shapes deformation of the fuel pellet, which predominantly affects its lower part, an \textit{a priori} refined mesh is employed. Thanks to adaptive mesh refinement preliminary studies~\cite{KRL_CMAME_2022}, the mesh is refined in the lower regions of both the pellet and the cladding. Since the meshes are matching in the (potential) contact region, a node-to-node contact discretization is adopted. The pellet fragment is discretized with $31,\!201$ nodes ($28,\!800$ elements), whereas the corresponding cladding portion is discretized with $8,\!569$ nodes ($7,\!200$ elements). The active contact surfaces comprise the entire outer surface of the pellet and the entire inner surface of the cladding (see Figure~\ref{fig:test_IPG_contact_zones}), resulting in $779$ contact nodes on each surface.}
{The strong mesh grading ($h_{min}\simeq 100 \,\mu\mathrm{m}$ and $h_{max}\simeq 511 \,\mu\mathrm{m}$) also makes this benchmark challenging for Uzawa-type contact algorithms. The large element-size ratio adversely affects the smallest eigenvalue of the stiffness matrix, thereby limiting the admissible augmentation parameter and potentially leading to a large number of Uzawa iterations.}

\subsubsection{Performance of standard and accelerated Uzawa methods} \label{ss:IPG_uzawa}

As for the Hertzian contact case, we first compare the performance of
the standard Uzawa formulation with various acceleration schemes.
Three augmentation parameter values within the admissible range,
$\rho \in \{10^{1},\,10^{2},\,10^{3}\}$, have been tested.

Based on the results of the previous academic example, a maximum number of iterations, $i_{\max}=5{,}000$, is imposed for all strategies, in addition to the classical convergence criterion ${r}^i \leq \epsilon$, with $\epsilon = 10^{-12}$ (see Eq.~\eqref{eq:relative_residual}).
This additional stopping criterion is introduced to highlight the ability of the proposed acceleration schemes to achieve sufficient accuracy within an acceptable number of iterations, and consequently within reasonable computational times.

\begin{figure}[htbp]
	\centering
	\begin{subfigure}{0.48\textwidth}
		\centering
		\includegraphics[width=\linewidth]{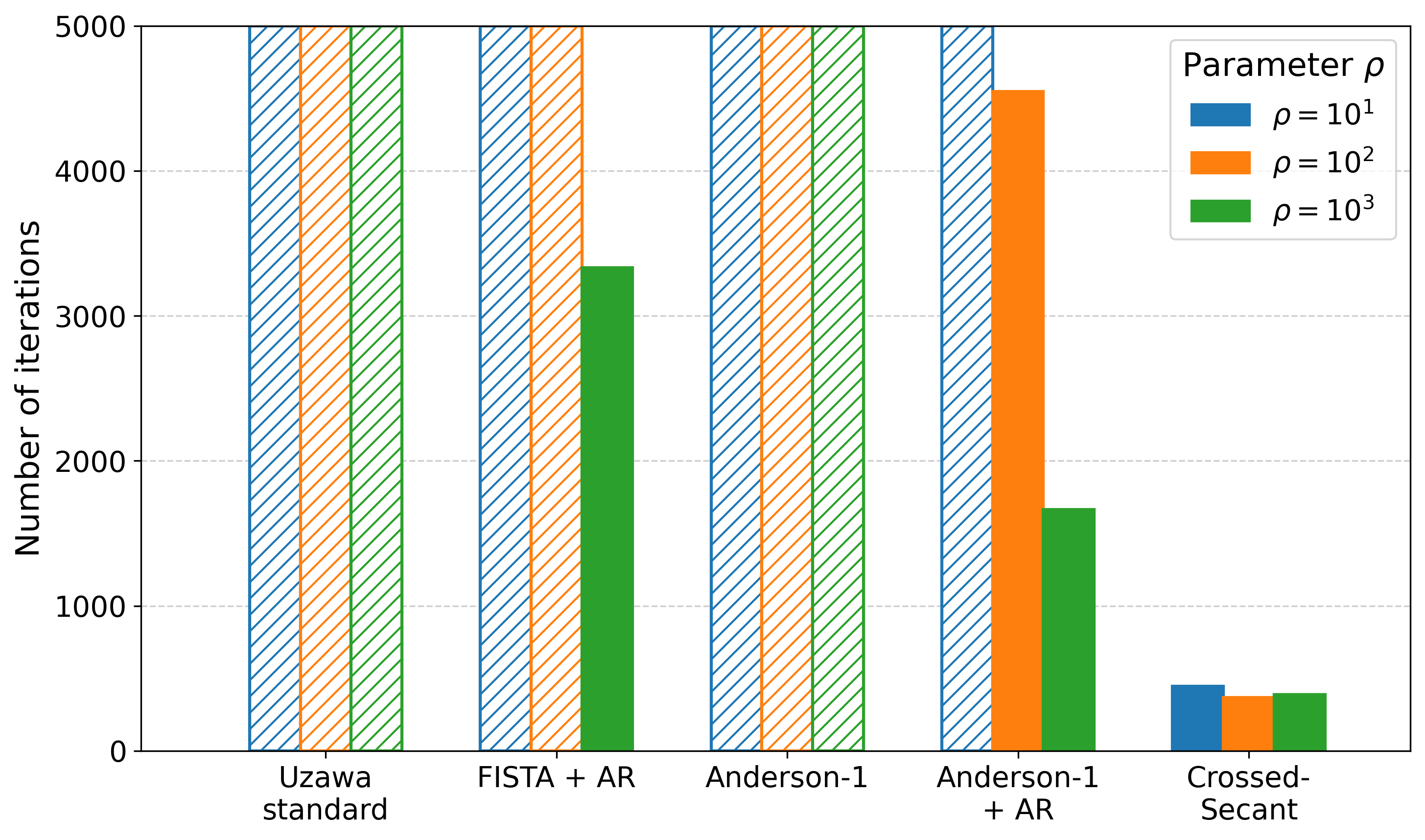}
		\caption{Number of iterations}
		\label{fig:IPG_uzawa_iterations_methods_diffrho}
	\end{subfigure}
	\hfill
	\begin{subfigure}{0.48\textwidth}
		\centering
		\includegraphics[width=\linewidth]{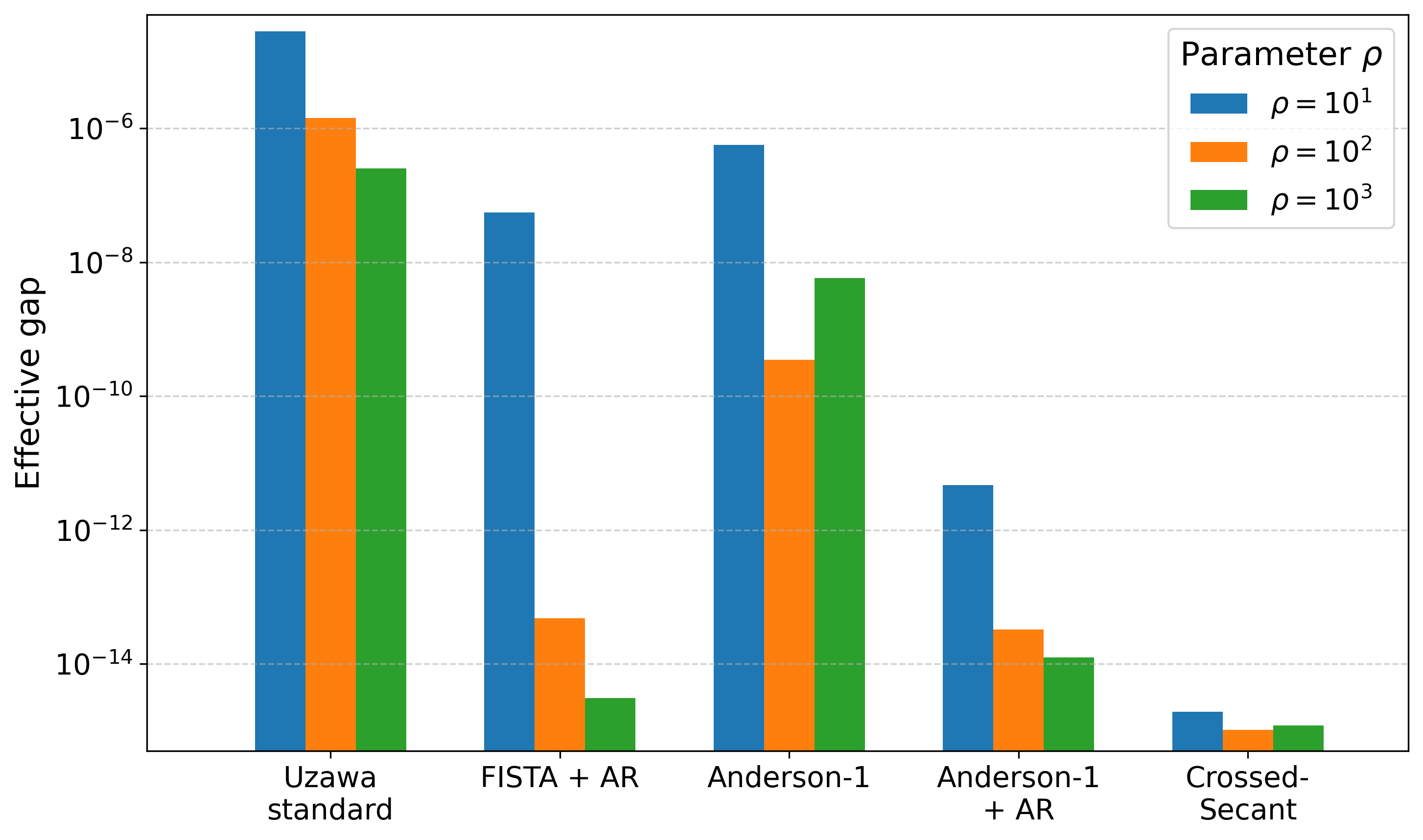}
		\caption{Effective gap}
		\label{fig:IPG_uzawa_gap_methods_diffrho}
	\end{subfigure}
	\caption{Industrial thermal-induced contact -- Uzawa/LM formulation -- performance of different acceleration schemes for $\rho \in {10^1, 10^2, 10^3}$ -- Stopping criterion: $r^i \leq \epsilon$ or $i = i_{\max}=5{,}000$. (a) Number of iterations required to satisfy the stopping criterion; methods that do not satisfy the convergence criterion ($r^i \leq \epsilon$) before reaching $i_{\max}$ are represented with hatched bars; (b) effective gap at the end of the iterative process.}
	\label{fig:IPG_compar_uza}
\end{figure}

Figures~\ref{fig:IPG_uzawa_iterations_methods_diffrho} and~\ref{fig:IPG_uzawa_gap_methods_diffrho} compare the convergence behavior and the final solution accuracy of the different acceleration strategies.

Figure~\ref{fig:IPG_uzawa_iterations_methods_diffrho} reports the number of iterations at the end of the algorithm. Strategies that fail to satisfy the convergence criterion before reaching the prescribed maximum number of iterations, $i_{\max}=5{,}000$, are displayed in hatched bars. The standard Uzawa method and Anderson-1 do not converge within the prescribed iteration budget for any tested value of $\rho$, whereas FISTA with adaptive restart converges only for $\rho=10^3$, and Anderson-1 with adaptive restart only for $\rho\ge10^2$. In contrast, the Crossed-Secant acceleration scheme is the only method that satisfies the convergence criterion for all tested values of $\rho$, requiring fewer than 450 iterations in every case. To further assess the quality of the computed solutions, Figure~\ref{fig:IPG_uzawa_gap_methods_diffrho} presents the effective gap at the end of the iterative process, i.e., either when the convergence criterion is satisfied or when the prescribed maximum number of iterations is reached. The Crossed-Secant scheme consistently yields the smallest effective gap, close to machine precision, for all tested values of $\rho$, whereas the remaining strategies exhibit significantly larger residual gaps especially when convergence is not achieved within the prescribed maximum number of iterations. In the latter case, for small $\rho$, the effective gaps remain comparable to the initial gap.
These results confirm the robustness of the proposed Crossed-Secant acceleration scheme, which combines rapid convergence with consistently high solution accuracy for all tested augmentation parameters within the admissible range.

\begin{figure}[h]
	\centering

	\includegraphics[trim=0cm 1cm 0cm 1cm,
	clip,width=\linewidth]{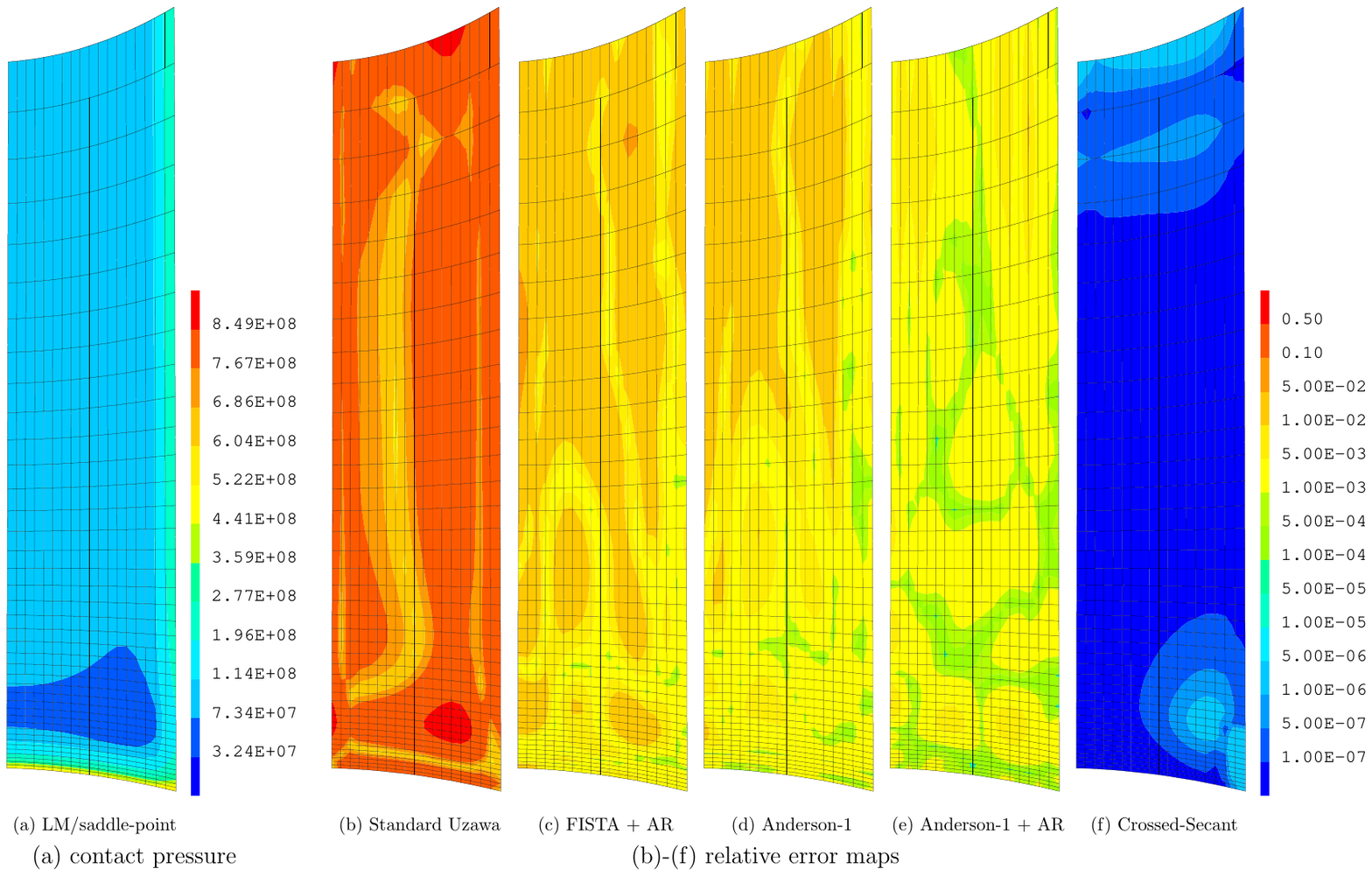}

	\caption{Industrial thermal-induced contact -- Uzawa/LM formulation -- performance of different acceleration schemes for $\rho = 10^3$, with all methods stopped at exactly $i_\text{max} = 388$; contact pressure and associated error maps on the cladding inner surface -- (a) contact pressure obtained with the reference LM/saddle-point solution; (b)--(f) relative error maps for the standard and accelerated Uzawa variants with respect to the reference LM/saddle-point solution}
	\label{fig:IPG_force_error_map}
\end{figure}

{
As the contact pressure is of particular interest from an industrial perspective, the accuracy of this latter ($p_c= - \sigma_{rr}$ here) is investigated for all (accelerated) Uzawa contact strategies setting the augmentation parameter $\rho=10^3$.
All results are reported after a fixed number of iterations, $i_{\max}=388$, corresponding to the number of iterations required for convergence of the Crossed-Secant acceleration scheme with this value of $\rho$.
This fixed iteration count allows a direct comparison of the contact pressure accuracy across the different strategies.
Figure~\ref{fig:IPG_force_error_map}(a) displays the contact pressure on the cladding inner surface obtained with the reference LM/saddle-point solution, while Figures~\ref{fig:IPG_force_error_map}(b)--(f) report the relative error maps with respect to the reference saddle-point solution.
As expected, large contact pressures are localized at the bottom of the cladding inner surface, where the hourglass-shape phenomenon occurs.
The error maps show that only the Crossed-Secant acceleration scheme reaches a very satisfactory error level with respect to the reference solution (relative error below $10^{-6}$). In contrast, the standard Uzawa, FISTA with adaptive restart, and Anderson-based strategies still exhibit significant errors after the same fixed number of iterations. These error maps are closely correlated with the convergence behavior of these methods.}



\subsubsection{Performance of Crossed-Secant-accelerated Uzawa and penalty-splitted methods}

{In this section, we perform a detailed analysis of the performance of the Crossed-Secant acceleration scheme, especially when the augmentation/penalty parameter is chosen outside the admissible range. The investigated range spans from $10^0$ to $10^{18}$.}

\begin{figure}[!htbp]
	\centering
	\begin{subfigure}{0.8\textwidth}
		\centering
		\includegraphics[width=\linewidth]{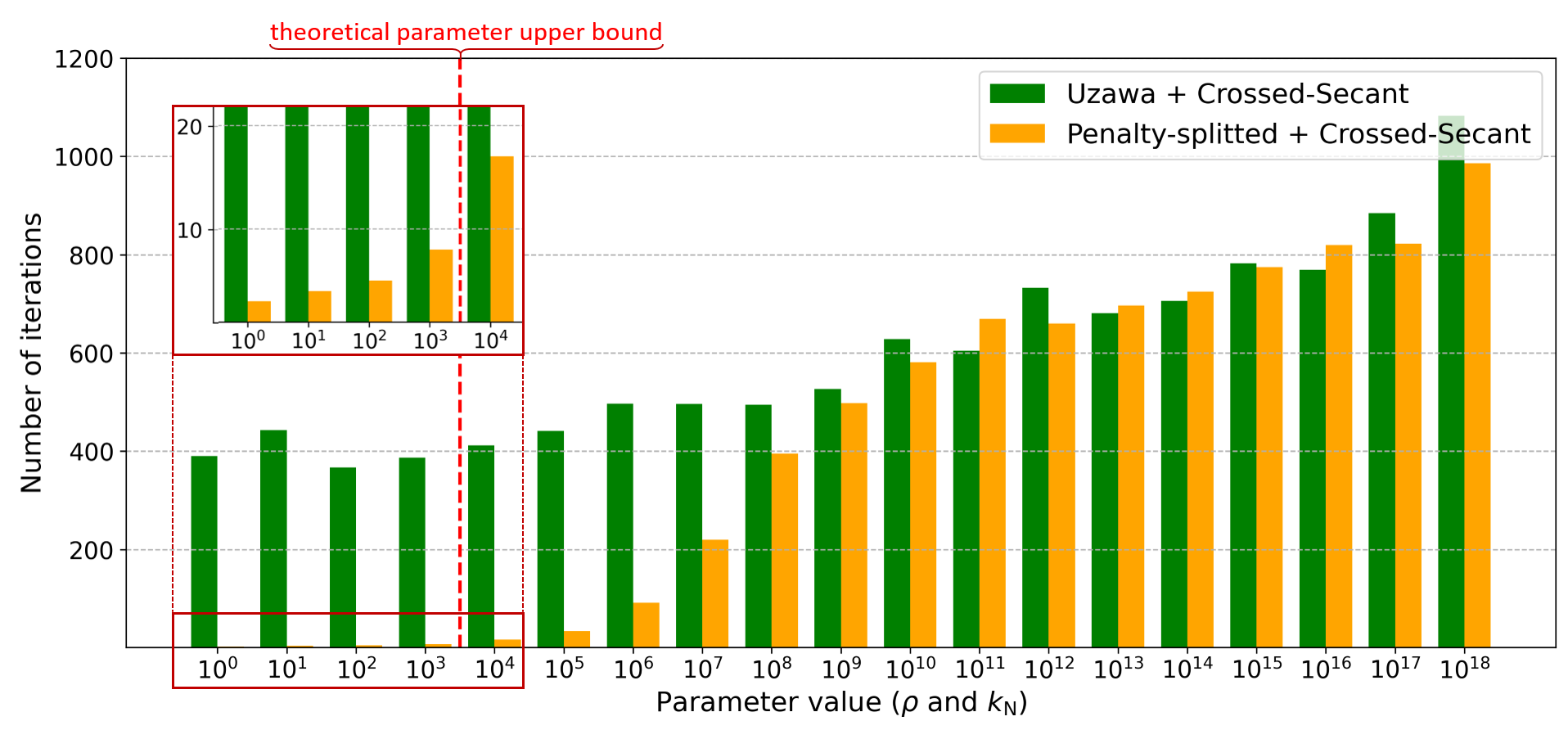}
		\caption{Number of iterations}
		\label{fig:IPG_uzawa_vs_penalty_ITER}
	\end{subfigure}

	\begin{subfigure}{0.8\textwidth}
		\centering
		\includegraphics[width=\linewidth]{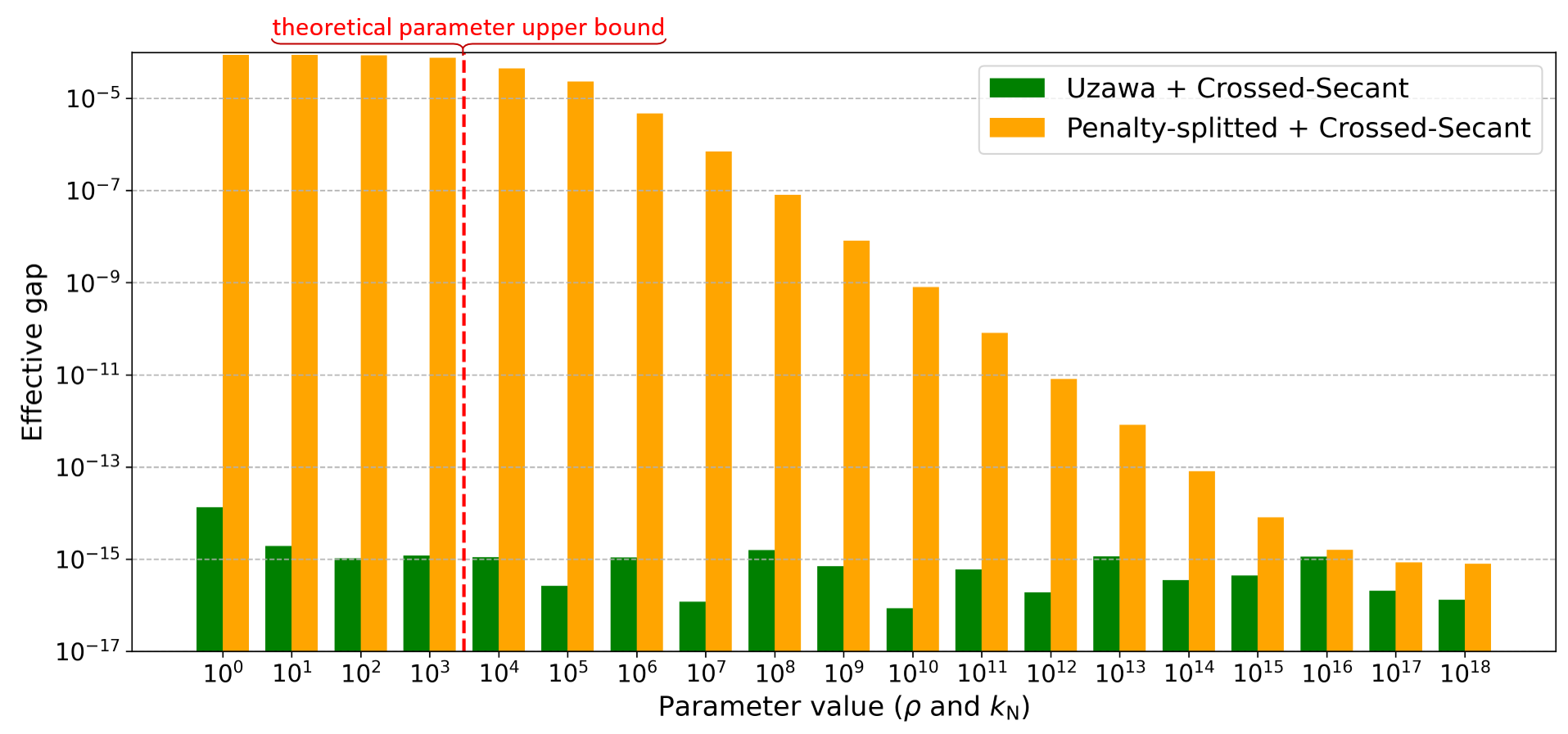}
		\caption{Effective gap}
		\label{fig:IPG_uzawa_vs_penalty_DN}
	\end{subfigure}
	\caption{Industrial thermal-induced contact -- performance of the Crossed-Secant accelerated Uzawa/LM and penalty-splitted methods for $\rho, k_N \in [10^0, 10^{18}]$ }
	\label{fig:IPG_accelCS_DN}
\end{figure}

{Figures~\ref{fig:IPG_uzawa_vs_penalty_ITER} and~\ref{fig:IPG_uzawa_vs_penalty_DN} compare the influence of the parameters, $\rho$ and $k_N$, on the convergence behavior and the final effective gap of the Crossed-Secant accelerated Uzawa and penalty-splitting formulations. 
  For the Uzawa formulation, the number of iterations remains relatively stable within the admissible range and increases slightly with the augmentation parameter outside this range. In contrast, for the penalty formulation, the number of iterations increases progressively with the penalty parameter $k_N$, particularly for small values of the penalty coefficient. For high penalty coefficients, the number of iterations becomes comparable to that obtained with the Uzawa method for the same value of the augmentation parameter.
  }
{This behavior can be explained by the strength of the coupling between the displacement solve and the contact force update within the iterative scheme. For small $k_N$, the contact force responds weakly to the gap, so the mechanical and contact updates are only loosely coupled: the iterative process converges quickly, but toward a solution with a significant residual gap. As $k_N$ increases, the contact force becomes increasingly sensitive to the gap, strengthening this coupling and making the problem numerically stiffer, approaching the behavior of the Uzawa formulation, which enforces the contact constraint strongly by construction. This explains why the iteration count for the penalty method progressively increases with $k_N$ and converges toward that of the Uzawa scheme.
However, this similarity in iteration counts between the two formulations does not reflect an equivalent convergence in solution accuracy. The effective gap obtained with the penalty formulation systematically improves as $k_N$ increases,
and only reaches error levels comparable to the Uzawa formulation (around $10^{-15}$~m) for sufficiently large values, $k_N \geq 10^{16}$.
This confirms that the residual constraint violation is inherent to the penalty method, and is only progressively reduced as $k_N$ increases, whereas the Uzawa formulation enforces the contact constraint essentially exactly, independently of $\rho$.
These results confirm the efficiency of the proposed Crossed-Secant accelerated penalty-splitting formulation, which enables highly penalized yet stable computations that would otherwise be numerically infeasible with classical penalty matrix methods due to ill-conditioning.}



\begin{figure}[!htbp]
	\centering
	
	\begin{subfigure}{\textwidth}
		\centering
		\includegraphics[page=1,width=1.\textwidth]{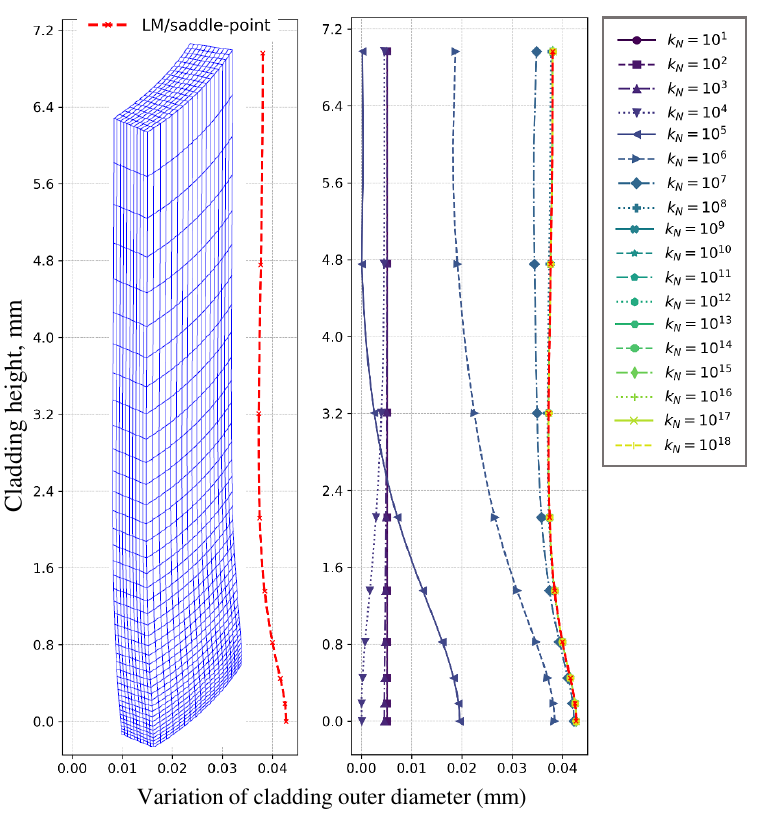}
	\end{subfigure}
	\caption{Industrial thermal-induced contact -- penalty-splitted formulation -- cladding profilometry ("bamboo-like" shape deformation)}
	\label{fig:IPG_penalty_evol}
\end{figure}

To complement the numerical analysis, Figure~\ref{fig:IPG_penalty_evol} compares the cladding profilometry obtained with the proposed penalty-splitted scheme (for different values of $k_N$) to the reference LM/saddle-point solution. 
Cladding outer profilometry provides direct evidence of the
pellet-cladding interaction and the associated residual strains. From our reference results (obtained with the saddle-point matrix solution), the simulated cladding's external diameter variations (mm), see
Figure~\ref{fig:IPG_penalty_evol}-left, successfully qualitatively reproduce the
characteristic primary ridges observed experimentally in
base post-irradiation profilometry (see~\cite{sercombe2020}, Figure (25)).
As shown in Figure~\ref{fig:IPG_penalty_evol}, small penalty parameters lead to cladding deformations that differ significantly from the reference solution, whereas for $k_N \geq 10^{9}$, the computed cladding profile is almost indistinguishable from the reference solution. This confirms that the use of large penalty parameters is crucial to obtain a reliable solution. Note that we focus here only on the penalty-splitted scheme, where the impact of the penalty parameter can be clearly observed. For the Crossed-Secant accelerated Uzawa method, the cladding profiles are superimposed on the reference solution even for small augmentation parameters, making it qualitatively impossible to see any deviations.

\subsection{Toward multi-domain contact problems and parallel performance} \label{ss:multi_spheres}

In order to further appreciate the potentiality of the proposed accelerated iterative methods, namely the Uzawa and penalty-splitted fixed-point strategies enhanced by the Crossed-Secant acceleration scheme, we investigate their computational performance as the number of interacting domains increases (Section~\ref{ss:multi_domain_contact}).
Then, in Section \ref{ss:hpc_insights}, for a fixed number of contacting domains, we give some insights on the parallel performance of these iterative methods by increasing the number of processors. The objective is to compare the computational efficiency of the proposed iterative accelerated strategies with the reference Lagrange multiplier saddle-point method. \\

\begin{figure}[!ht]
	\captionsetup[subfigure]{labelformat=empty}
	\centering
	\includegraphics[page=1, width=1.\linewidth]{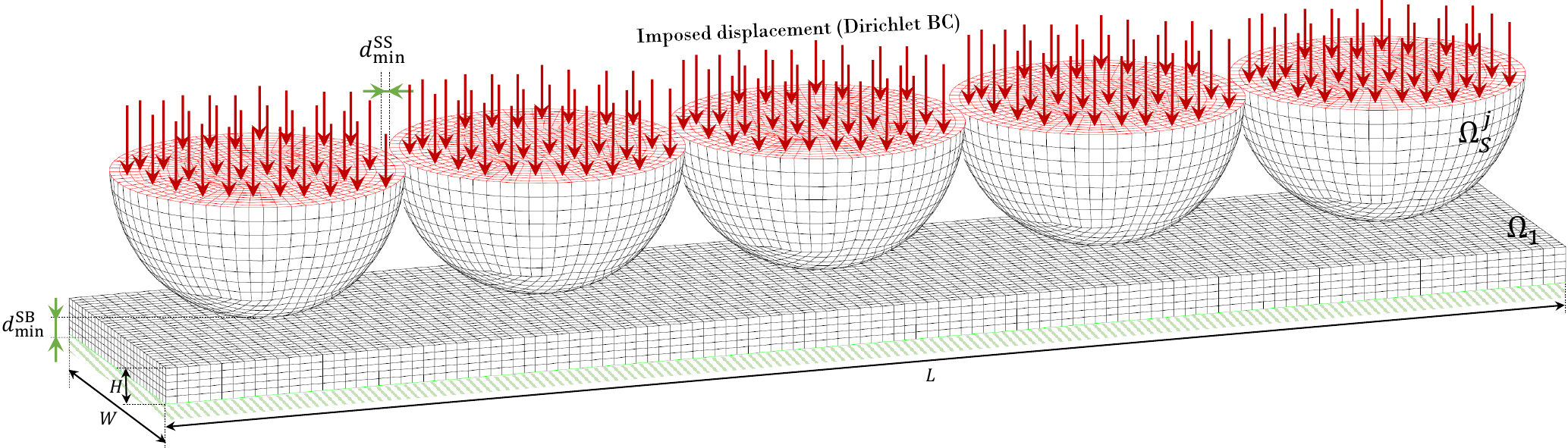}
	\caption{Multi-domain contact problem -- computational domains with boundary conditions for $N_S = 5$}
	\label{fig:multi_sph}
      \end{figure}

For the numerical experiments, we consider a basic setup similar to the Hertzian contact case described in Section~\ref{ss:Hertz}, but with an increasing number of deformable half-spheres $\{\Omega_S^j\}_{j=1}^{N_S}$ coming into contact with an elastic block $\Omega_1$. A sketch of the problem is shown in Figure~\ref{fig:multi_sph} for $N_S = 5$. All domains are modeled as deformable, isotropic, linear elastic solids. Each sphere has the same radius $R = 2\cdot 10^{-2}$ m and material properties as described in Section~\ref{ss:Hertz}.  The block domain $\Omega_1$ has a thickness of $H= 0.5 \cdot 10^{-2}$ m, width $W= 2 R$ m, and length which varies depending on the number of spheres $L = 2 R  N_S$ m.  The material properties of the block remain those specified in Section~\ref{ss:Hertz}.

A minimal initial gap of $d^{\text{SB}}_{\text{min}} = 1 \cdot 10^{-3}$~m is prescribed between the spheres and the block, while a zero initial gap is imposed between neighboring spheres, $d_{\text{min}}^{\text{SS}} = 0$~m.
Contact is enforced by applying a vertical displacement of $u_\text{D} = 4 \cdot 10^{-3}$~m at the top of each half-sphere, inducing indentation in the block. The initial zero gap between spheres results in mutual contact at their equatorial planes.
Between the spheres, a node-to-node contact is assumed since the meshes are matching, whereas node-to-surface pairing is used for the contact between the block and each sphere.

For the numerical application of the Crossed-Secant–accelerated fixed-point schemes, the following parameter choices are adopted. For the accelerated Uzawa method, the augmentation parameter is fixed at $\rho = 10^{4}$ (in-range value), which allows achieving the desired precision in a limited number of iterations.  
For the penalty-splitted accelerated method, the penalty parameter is set to $k_N = 10^{18}$ (out-of-range value) in order to reach a level of accuracy comparable to that obtained with the Uzawa approach.

\subsubsection{Multi-domain contact}
\label{ss:multi_domain_contact}

We first assess the performance of the sequential computation as the number of spheres increases, with $N_S = 1,\dots,20$, corresponding to a total number of domains of $N_S + 1$.

\begin{figure}[!ht]
	\captionsetup[subfigure]{labelformat=empty}
	\centering
	\includegraphics[width=1.\linewidth]{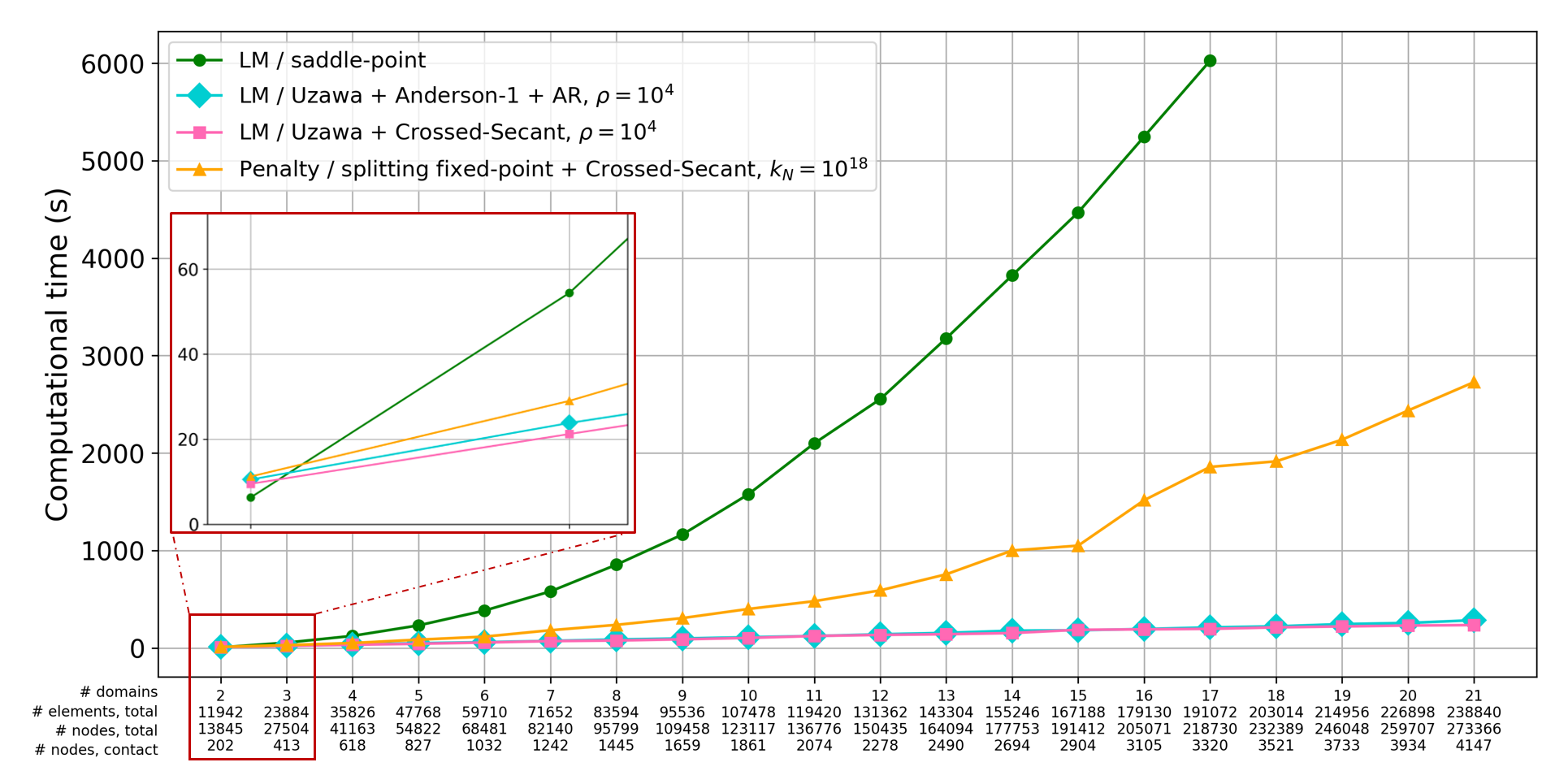}
	\caption{Multi-domain contact problem -- computational time as a function of the number of interacting domains (sequential implementation)}
	\label{fig:multi_sp_ncpu1}
\end{figure}

Figure~\ref{fig:multi_sp_ncpu1} reports the computational times (in seconds) as a function of the number of interacting domains, the total number of nodes, and the number of nodes in active contact. These computational times correspond to the solution process. For the accelerated iterative scheme, this consists in the application of Algorithm~\ref{alg:algo_uni}. For the standard Lagrange multiplier method, it corresponds to the solution of the underlying saddle-point system, see Eq.~\eqref{eq:saddle_point_system}, which includes an active-set outer loop to update the constraint set.

First, it can be observed that for a small number of contacting domains, the standard saddle-point matrix solution is only slightly more computationally efficient than the accelerated fixed-point splitting methods. However, its computational cost increases rapidly with the number of contacting domains and quickly becomes prohibitive. Consequently, for more than $4$ domains in contact (i.e., $3$ half-spheres in contact with the block), this approach is outperformed by the accelerated splitting fixed-point methods.  
It is worth noting that computations using the standard LM/saddle-point formulation were not feasible beyond $17$ contacting domains, corresponding to approximately $220\,000$ nodes, due to the increasing problem complexity. In contrast, the proposed accelerated iterative methods are able to handle larger configurations with ease.      

For the Crossed-Secant accelerated penalty-splitted strategy, it has been observed that as the number of interacting domains increases, a larger number of iterations is required for convergence, from $55$ iterations for $N_S = 1$ to $1\,117$ for $N_S = 20$. In particular, more iterations are needed to identify the final active constraint set, i.e., the correct contact zone. This behavior is typical of the out-of-range parameter effects observed in the previous examples with two contact bodies. The size of the initial plateau (see Figure~\ref{fig:sphere_plan_penalty_accelCS_RSDnrmEu}) is therefore also related to the number of contact domains.

For the Crossed-Secant accelerated Uzawa formulation with an in-range parameter, $\rho = 10^{4}$, it achieves the lowest computational times across the entire range of contact domain numbers. Due to the in-range parameter value, the number of iterations remains very limited and stable, between $32$ and $40$ for $1 \le N_S \le 20$, highlighting strong robustness with respect to both problem size and the number of contact constraints. 

{As the improved Anderson method proposed in this contribution, by combining it with adaptive restart, has shown promising results in the previous examples, its results on the multi-domain contact problem are also reported in Figure~\ref{fig:multi_sp_ncpu1}, for the same augmentation parameter $\rho = 10^{4}$. As expected, its performance in terms of computation time is similar, though slightly worse, to that of the Crossed-Secant/Uzawa algorithm, as it converges in a few more iterations ($58$ to $68$ iterations for $1 \leq N_S \leq 20$). However, let us underline that this algorithm is only applicable for augmentation parameters within the admissible range of the Uzawa method.}


In terms of accuracy, all accelerated approaches achieve a maximum effective gap of approximately $10^{-15}$~m and a complementarity condition measure of $10^{-12}~\mathrm{N\,m}$. In the case of the penalty-splitted approach, we once again emphasize that the proposed Crossed-Secant acceleration scheme makes it possible to handle very large penalty parameters efficiently, thereby achieving levels of accuracy that surpass those generally attainable with the standard penalty method.

\subsubsection{Insights into HPC potential and parallel performance}
\label{ss:hpc_insights}

We extend our analysis to parallel computations by increasing the number of processors for a fixed number of contacting domains, $N_S = 15$ (i.e., a total of $16$ domains). This allows us to provide insights into the parallel performance of the proposed accelerated splitting fixed-point methods.

\begin{figure}[!ht]
	\captionsetup[subfigure]{labelformat=empty}
	\centering
	\includegraphics[width=1.\linewidth]{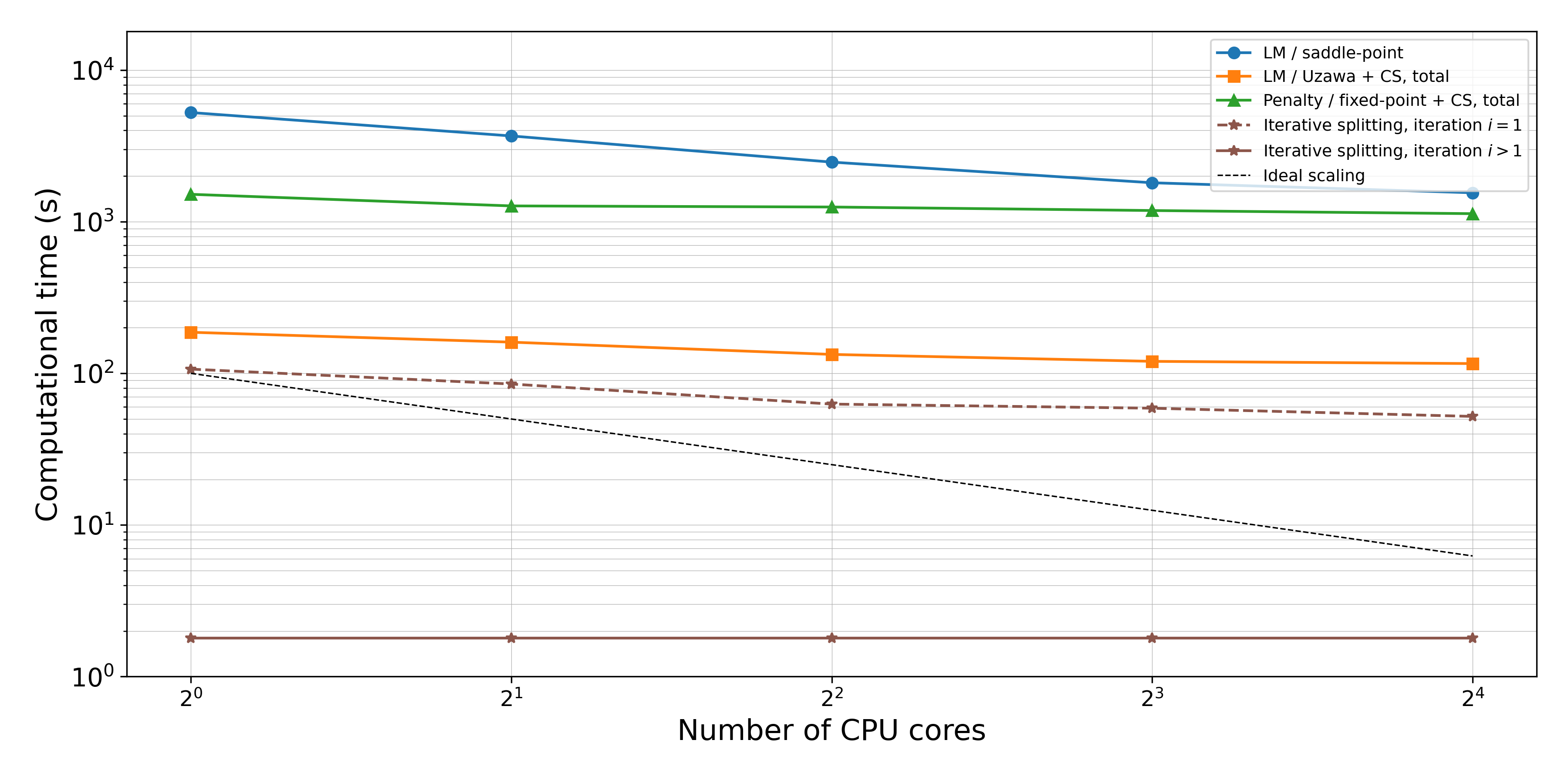}
	\caption{Multi-domain contact problem -- computational times for increasing number of contact domains}
	\label{fig:multi_sp_ncpu}
\end{figure}

Figure~\ref{fig:multi_sp_ncpu} presents the computational times obtained with the three solving strategies considered in this example. For the accelerated fixed-point methods, three measures are reported: the total computational time, the computational time of the first iteration, and the characteristic time of a single iteration after the first one.
Let us emphasize that the last two measures are identical for the accelerated Uzawa and penalty-splitted strategies, implying that the difference in total computational time depends solely on the difference in the number of iterations required to reach convergence.
Moreover, we recall that the matrix system is factorized during the first iteration and reused in subsequent iterations, which explains the large difference in computational time between the first iteration (around $10^2$~s) and the following ones (around $1.7$~s).

Finally, only the standard Lagrange multiplier saddle-point strategy and the standard stiffness direct solution benefit from parallel implementation in the Cast3M finite element software~\cite{CEA_CAST3M_2023} used for this demonstration.
For instance, matrix-vector products and scalar products are not yet parallelized. As a result, the cost of iterations beyond the first remains constant with respect to the number of processors. Addressing this issue could improve the strong scaling performance of the accelerated fixed-point strategies.

The first conclusion to be drawn is that, even though the standard LM/saddle-point solution strategy exhibits some computational time savings with an increasing number of processors, the strong scaling is far from ideal. However, as already mentioned in the introduction, the efficient parallelization of saddle-point systems remains an ongoing research topic.
The same conclusion applies to the strong scaling behavior of the standard stiffness system solution (cf. the first iteration of the iterative splitting schemes). However, for the latter, better parallelization is available in some modern finite element software.

Globally, the computational time of the saddle-point matrix solution is roughly two orders of magnitude larger than that of the stiffness matrix solution. This confirms the interest of iterative splitting algorithms, provided that the number of iterations remains limited. The difference in computational time between the two accelerated fixed-point strategies depends on the number of iterations required for convergence: $40$ for the Uzawa method and $720$ for the penalty-term–treated approach. Consequently, the Uzawa method reduces the saddle-point solution time by a factor of approximately $20$, while the penalty-splitted approach achieves a speed-up of about $2$ to $3$.
Note that using an out-of-range parameter for the Uzawa method (e.g., $\rho = 10^{18}$, equivalent to the penalty parameter used here) leads to a similar number of iterations and computational time as the out-of-range penalty method. 

In our study, the ultimate parallel efficiency of the methods is primarily constrained by solver-level design choices rather than by the algorithms themselves. Future work will therefore focus on implementing the proposed methods within a high-performance finite element solver to fully exploit their scalability potential.
{Indeed, the proposed accelerated operator splitting schemes seem well suited for parallel environments due to their intrinsic algorithmic structure. In particular, the left-hand side matrix is a regular rigidity matrix, for which  existing scalable linear solvers could be used. Moreover, it remains unchanged across iterations, allowing for efficient reuse of matrix factorizations or iterative solver initialization. The computational overhead associated with the proposed acceleration strategies remains negligible at the iteration level, as they only involve a limited number of vector operations with linear complexity. However, since these operations are repeated at each iteration, their efficient parallelization is crucial to preserve scalability when convergence requires a large number of iterations. Their simple structure makes them naturally suitable for large-scale distributed-memory implementations.}

\section{Conclusion} \label{s:conclusions}

This work introduces a unified iterative framework for the efficient solution of frictionless mechanical contact problems using an operator splitting strategy. Each iteration involves solving a standard stiffness system with a constant matrix, and updating the contact forces through an expression depending on the contact numerical formulation. For the Lagrange multiplier formulation, the Uzawa scheme is recovered, while a penalty-splitted scheme is obtained for the penalty formulation.

A central contribution of this study is the integration of the Crossed-Secant acceleration strategy within this unified framework, applied before the projection onto the admissible set. When applied to both Uzawa and penalty-based splitting schemes, the Crossed-Secant accelerated scheme provides faster convergence than other acceleration strategies reported in the literature, including the enhanced Anderson acceleration proposed here with an adaptive restart procedure. Moreover, the Crossed-Secant scheme also demonstrates remarkable robustness. In particular, a key and distinctive feature of the proposed approach is its ability to reach convergence using parameter values, namely the augmentation parameter and the penalty parameter, beyond the restrictive theoretical and numerical bounds. This improvement removes one of the main bottlenecks of classical Uzawa and penalty formulations, which typically require careful and problem-dependent parameter tuning.

For the Uzawa algorithm, the Crossed-Secant acceleration achieves the machine-level precision independently of the parameter value. 
For the penalty-splitted formulation, the Crossed-Secant acceleration enables stable and efficient convergence for very large penalty parameters, making it possible to reach levels of accuracy comparable to those of the Uzawa and saddle-point formulations, that is generally unattainable with classical penalty approaches due to matrix conditioning and convergence issues. {Importantly, this improvement is obtained while solving only standard stiffness systems, without the need for specialized preconditioning.}

Beyond its convergence and robustness, the proposed framework exhibits increasing efficiency compared with the classical Lagrange multiplier saddle-point formulation as the number of contacting domains grows. The proposed methods also appear well suited for large-scale simulations. Numerical experiments on multi-domain contact configurations already reveal encouraging trends, highlighting their potential for large (multi-)contact systems. Future work will further investigate scalability aspects, with particular emphasis on massively parallel high-performance computing implementations.

Our strategy is intended to be extended to nonlinear material behaviors and frictional contact formulations to further demonstrate the versatility and practical relevance of the proposed framework for large-scale, industrial-grade simulations.

\section{Acknowledgements}

The authors are grateful to the PLEIADES project, financially supported by the CEA that funded this research work.



        
\bibliographystyle{unsrt}
\bibliography{bibl_contact}

@article{ALART1991,
title = {A mixed formulation for frictional contact problems prone to {N}ewton like solution methods},
journal = {Computer Methods in Applied Mechanics and Engineering},
volume = {92},
number = {3},
pages = {353-375},
year = {1991},
issn = {0045-7825},
doi = {https://doi.org/10.1016/0045-7825(91)90022-X},
url = {https://www.sciencedirect.com/science/article/pii/004578259190022X},
author = {P. Alart and A. Curnier}
}

@article{Gittus72,
	title = {Theoretical analysis of the strains produced in nuclear fuel cladding tubes by the expansion of cracked cylindrical fuel pellets},
	journal = {Nuclear Engineering and Design},
	volume = {18},
	number = {1},
	pages = {69-82},
	year = {1972},
	issn = {0029-5493},
	doi = {https://doi.org/10.1016/0029-5493(72)90037-4},
	url = {https://www.sciencedirect.com/science/article/pii/0029549372900374},
	author = {J.H. Gittus}
}

@book{Nocedal2006,
	address = {New York, NY},
	author = {J. Nocedal and S.J Wright},
	edition = {2. ed.},
	isbn = {978-0-387-30303-1},
	pagetotal = {XXII, 664},
	ppn_gvk = {502988711},
	publisher = {Springer},
	series = {Springer series in operations research and financial engineering},
	title = {Numerical optimization},
	url = {http://gso.gbv.de/DB=2.1/CMD?ACT=SRCHA&SRT=YOP&IKT=1016&TRM=ppn+502988711&sourceid=fbw_bibsonomy},
	year = 2006
}

@article{SCHAPIRA2023,
	title = {Performance of acceleration techniques for staggered phase-field solutions},
	journal = {Computer Methods in Applied Mechanics and Engineering},
	volume = {410},
	pages = {116029},
	year = {2023},
	issn = {0045-7825},
	doi = {https://doi.org/10.1016/j.cma.2023.116029},
	url = {https://www.sciencedirect.com/science/article/pii/S0045782523001536},
	author = {Yaron Schapira and Lars Radtke and Stefan Kollmannsberger and Alexander Düster}
}

@article{Dafermos1992,
	author = {Dafermos, S. C. and McKelvey, S. C.},
	title = {Partitionable variational inequalities with applications to network and economic equilibria},
	year = {1992},
	issue_date = {May       1992},
	publisher = {Plenum Press},
	address = {USA},
	volume = {73},
	number = {2},
	issn = {0022-3239},
	journal = {J. Optim. Theory Appl.},
	month = may,
	pages = {243–268},
	numpages = {26}
}

@book{Kinderlehrer2000,
	author = {Kinderlehrer, D. and Stampacchia, G.},
	title = {An Introduction to Variational Inequalities and Their Applications},
	publisher = {Society for Industrial and Applied Mathematics},
	year = {2000},
	doi = {10.1137/1.9780898719451},
	address = {},
	edition   = {},
	URL = {https://epubs.siam.org/doi/abs/10.1137/1.9780898719451},
	eprint = {https://epubs.siam.org/doi/pdf/10.1137/1.9780898719451}
}

@article{ODEN1980,
	title = {Theory of variational inequalities with applications to problems of flow through porous media},
	journal = {International Journal of Engineering Science},
	volume = {18},
	number = {10},
	pages = {1173-1284},
	year = {1980},
	issn = {0020-7225},
	doi = {https://doi.org/10.1016/0020-7225(80)90111-1},
	url = {https://www.sciencedirect.com/science/article/pii/0020722580901111},
	author = {J.T. Oden and N. Kikuchi}
}

@article{Lions1967,
	author = {Lions, J. L. and Stampacchia, G.},
	title = {Variational inequalities},
	journal = {Communications on Pure and Applied Mathematics},
	volume = {20},
	number = {3},
	pages = {493-519},
	doi = {https://doi.org/10.1002/cpa.3160200302},
	url = {https://onlinelibrary.wiley.com/doi/abs/10.1002/cpa.3160200302},
	eprint = {https://onlinelibrary.wiley.com/doi/pdf/10.1002/cpa.3160200302},
	year = {1967}
}

@BOOK{Duvaut1976,
   title = {Inequalities in mechanics and physics.},
   author = {Duvaut, G. and Lions, J.-L.},
   publisher = {Springer Science \& Business Media},
   year = {1976},
   doi={10.1007/978-3-642-66165-5}
}

@BOOK{Dautray1988,
   title = {Mathematical Analysis and Numerical Methods for Science and Technology. Volume 2 : Functional And Variational Methods},
   author = {Dautray, R. and Lions, J.-L.},
   publisher = {Springer},
   year = {1988}
}

@article{Karatzas2000,
	author = {Karatzas, I. and Shreve, S.},
	year = {2000},
	month = {06},
	pages = {},
	title = {Methods of Mathematical Finance},
	volume = {95},
	journal = {Journal of the American Statistical Association},
	doi = {10.2307/2669426}
}

@book{Migorski2013,
	TITLE = {{Nonlinear Inclusions and Hemivariational Inequalities. Models and Analysis of
	Contact Problems}},
	AUTHOR = {Migorski, S. and Ochal, A. and Sofonea, M.},
	URL = {https://hal.science/hal-01898484},
	PUBLISHER = {{Springer-Verlag New York}},
	SERIES = {Advances in Mechanics and Mathematics},
	VOLUME = {26},
	YEAR = {2013},
	DOI = {10.1007/978-1-4614-4232-5},
	HAL_ID = {hal-01898484},
	HAL_VERSION = {v1},
}

@article{Han2019,
	TITLE = {{Numerical analysis of hemivariational inequalities in contact mechanics}},
	AUTHOR = {Han, W. and Sofonea, M.},
	URL = {https://hal.science/hal-03544749},
	JOURNAL = {{Acta Numerica}},
	PUBLISHER = {{Cambridge University Press (CUP)}},
	VOLUME = {28},
	PAGES = {175-286},
	YEAR = {2019},
	MONTH = May,
	DOI = {10.1017/S0962492919000023},
	HAL_ID = {hal-03544749},
	HAL_VERSION = {v1},
}

@book{Sofonea2018,
	TITLE = {{Variational-Hemivariational Inequalities with Applications}},
	AUTHOR = {Sofonea, M. and Migorski, S.},
	URL = {https://hal.science/hal-01627037},
	PUBLISHER = {{Chapman \& Hall/CRC}},
	SERIES = {Monographs and Research Notes in Mathematics},
	YEAR = {2018},
	HAL_ID = {hal-01627037},
	HAL_VERSION = {v1},
}

@Book{Panagiotopoulos1993,
  author = 	 {Panagiotopoulos, P. D.},
  title = 	 {Hemivariational Inequalities: Applications in Mechanics and Engineering},
  publisher = 	 {Springer},
  year = 	 {1993},
}

@article{HUANG2022,
title = {An accelerated method of {U}zawa algorithm in contact problems},
journal = {Computers \& Mathematics with Applications},
volume = {127},
pages = {97-104},
year = {2022},
issn = {0898-1221},
doi = {https://doi.org/10.1016/j.camwa.2022.09.030},
url = {https://www.sciencedirect.com/science/article/pii/S0898122122004163},
author = {Z. Huang and X. Cheng}
}

@article{KANNO2020,
author = {Y. Kanno},
year = {2020},
month = {10},
pages = {},
title = {An Accelerated {U}zawa Method for Application to Frictionless Contact Problem},
volume = {14},
journal = {Optimization Letters},
doi = {10.1007/s11590-019-01481-2},
url = {https://doi.org/10.1007/s11590-019-01481-2}
}

@article{RAMIERE2015,
	title = {Iterative residual-based vector methods to accelerate fixed point iterations},
	journal = {Computers \& Mathematics with Applications},
	volume = {70},
	number = {9},
	pages = {2210-2226},
	year = {2015},
	issn = {0898-1221},
	doi = {https://doi.org/10.1016/j.camwa.2015.08.025},
	url = {https://www.sciencedirect.com/science/article/pii/S0898122115004046},
	author = {I. Ramière and T. Helfer}
}

@article{MOTA2025,
	author = {A. Mota and D. Koliesnikova and I. Tezaur and J. Hoy},
	title = {A Fundamentally New Coupled Approach to Contact Mechanics via the {D}irichlet-{N}eumann {S}chwarz Alternating Method},
	journal = {International Journal for Numerical Methods in Engineering},
	volume = {126},
	number = {9},
	pages = {e70039},
	doi = {https://doi.org/10.1002/nme.70039},
	url = {https://onlinelibrary.wiley.com/doi/abs/10.1002/nme.70039},
	eprint = {https://onlinelibrary.wiley.com/doi/pdf/10.1002/nme.70039},
	year = {2025}
}

@article{Signorini1959,
	author =       {Signorini, A.},
	title =        {Questioni di elasticit{\`a} non linearizzata e
	semilinearizzata},
	journal =      {Rendiconti di Matematica e delle sue Applicazioni,
	Quinta Serie},
	ISSN =         {0034-4427},
	volume =       18,
	pages =        {95--139},
	year =         1959,
	language =     {Italian},
}

@incollection{Fichera1972,
	author =       {G. Fichera},
	editor =       {Siegfried Fl{\"u}gge and Clifford A. Truesdell},
	title =        {Boundary value problems of elasticity with
	unilateral constraints},
	booktitle =    {Festk{\"o}rpermechanik/Mechanics of Solids},
	series =       {Handbuch der Physik (Encyclopedia of Physics)},
	volume =       {VIa/2},
	year =         1972,
	edition =      {paperback 1984},
	pages =        {391--424},
	publisher =    {Springer-Verlag},
	address =      {Berlin--Heidelberg--New York},
	isbn =         {0-387-13161-2},
	zbl =          {0277.73001}
}

@incollection{Epalle2025,
title = {Parallel simulation and adaptive mesh refinement for {3D} elastostatic contact mechanics problems between deformable bodies},
editor = {F. Chouly and S. Bordas and R. Becker and P. Omnes},
series = {Advances in Applied Mechanics},
publisher = {Elsevier},
volume = {61},
pages = {287-345},
year = {2025},
booktitle = {Error Control, Adaptive Discretizations, and Applications, Part 4},
issn = {0065-2156},
doi = {https://doi.org/10.1016/bs.aams.2025.08.003},
url = {https://www.sciencedirect.com/science/article/pii/S0065215625000225},
author = {A. Epalle and I. Ramière and G. Latu and F. Lebon}
}

@article{HO2017,
	author = {Ho, N. and Olson, S. D. and Walker, H. F.},
	title = {Accelerating the {U}zawa {A}lgorithm},
	journal = {SIAM Journal on Scientific Computing},
	volume = {39},
	number = {5},
	pages = {S461-S476},
	year = {2017},
	doi = {10.1137/16M1076770},	
	URL = {https://doi.org/10.1137/16M1076770},
	eprint = {https://doi.org/10.1137/16M1076770}
}

@article{Walker2011,
	author = {Walker, H. F. and Ni, P.},
	title = {Anderson Acceleration for Fixed-Point Iterations},
	journal = {SIAM Journal on Numerical Analysis},
	volume = {49},
	number = {4},
	pages = {1715-1735},
	year = {2011},
	doi = {10.1137/10078356X},
	URL = {https://doi.org/10.1137/10078356X},
	eprint = {https://doi.org/10.1137/10078356X}
}

@article{Nochetto2004,
	author = {Nochetto, R. H. and Pyo, J.},
	title = {Optimal relaxation parameter for the {U}zawa Method},
	year = {2004},
	issue_date = {October 2004},
	publisher = {Springer-Verlag},
	address = {Berlin, Heidelberg},
	volume = {98},
	number = {4},
	issn = {0029-599X},
	url = {https://doi.org/10.1007/s00211-004-0522-0},
	doi = {10.1007/s00211-004-0522-0},
	journal = {Numer. Math.},
	month = oct,
	pages = {695–702},
	numpages = {8}
}

@article{Bacuta2006,
	author = {Bacuta, C.},
	title = {A Unified Approach for {U}zawa Algorithms},
	journal = {SIAM Journal on Numerical Analysis},
	volume = {44},
	number = {6},
	pages = {2633-2649},
	year = {2006},
	doi = {10.1137/050630714},
	URL = {https://doi.org/10.1137/050630714},
	eprint = {https://doi.org/10.1137/050630714}
}

@article{Chouly2013Penalty,
	author = {Chouly, F. and Hild, P.},
	title = {On convergence of the penalty method for unilateral contact problems},
	year = {2013},
	issue_date = {March, 2013},
	publisher = {Elsevier Science Publishers B. V.},
	address = {NLD},
	volume = {65},
	issn = {0168-9274},
        url = {https://doi.org/10.1016/j.apnum.2012.10.003},
	doi = {10.1016/j.apnum.2012.10.003},
	journal = {Appl. Numer. Math.},
	month = mar,
	pages = {27–40},
	numpages = {14}
}

@article{Franceschini2022,
title = {A reverse augmented constraint preconditioner for {L}agrange multiplier methods in contact mechanics},
journal = {Computer Methods in Applied Mechanics and Engineering},
volume = {392},
pages = {114632},
year = {2022},
issn = {0045-7825},
doi = {https://doi.org/10.1016/j.cma.2022.114632},
url = {https://www.sciencedirect.com/science/article/pii/S0045782522000421},
author = {A. Franceschini and M. Ferronato and M. Frigo and C. Janna}
}

@article{BECK2009,
	author = {A. Beck and M. Teboulle},
	title = {A Fast Iterative Shrinkage-Thresholding Algorithm for Linear Inverse Problems},
	journal = {SIAM Journal on Imaging Sciences},
	volume = {2},
	number = {1},
	pages = {183-202},
	year = {2009},
	doi = {10.1137/080716542},
	URL = {https://doi.org/10.1137/080716542}
}

@book{Laursen2003,
	author = {Laursen, T.},
	year = {2003},
	month = {01},
	pages = {476},
	title = {Computational Contact and Impact Mechanics: Fundamentals of Modeling Interfacial Phenomena in Nonlinear Finite Element Analysis},
	publisher= {Springer},
	isbn = {9783540429067},
	doi = {10.1007/978-3-662-04864-1}
}

@article{Chouly2013Nitsche,
	author = {Chouly, F. and Hild, P.},
	title = {A {N}itsche-Based Method for Unilateral Contact Problems: Numerical Analysis},
	journal = {SIAM Journal on Numerical Analysis},
	volume = {51},
	number = {2},
	pages = {1295-1307},
	year = {2013},
	doi = {10.1137/12088344X},
	URL = {https://doi.org/10.1137/12088344X},
	eprint = {https://doi.org/10.1137/12088344X}
}

@article{Chouly2019,
	author =       {F. Chouly and M. Fabre and P. Hild and R. Mlika and
	J. Pousin and Y. Renard},
	journal =      {UCL Workshop 2016, UCL (University College London),
	Jan 2016, London, United Kingdom},
	pages =        {94--141},
	title =        {An overview of recent results on {N}itsche's method
	for contact problems},
	year =         2019
}

@book{Glowinski1989,
  title        = {{Augmented Lagrangian and Operator-Splitting Methods in Nonlinear Mechanics}},
  author       = {Glowinski, R. and Le Tallec, P.},
  year         = {1989},
  series       = {Studies in Applied and Numerical Mathematics},
  volume       = {9},
  publisher    = {Society for Industrial \& Applied Mathematics (SIAM)},
  address      = {Philadelphia, PA, USA},
  isbn         = {9780898712308},
  note         = {SIAM Studies in Applied Mathematics},
}

@article{Krause2002,
	author = {Krause, R. and Wohlmuth, B.},
	year = {2002},
	month = {12},
	pages = {139-148},
	title = {A {Dirichlet–Neumann} type algorithm for contact problems with friction},
	volume = {5},
	journal = {Computing and Visualization in Science},
	doi = {10.1007/s00791-002-0096-2}
}

@article{Puso2004_mortar3D,
	author =       {Puso, M. A.},
	doi =          {https://doi.org/10.1002/nme.865},
	journal =      {International Journal for Numerical Methods in
	Engineering},
	number =       3,
	pages =        {315-336},
	title =        {A {3D} mortar method for solid mechanics},
	url =
	{https://onlinelibrary.wiley.com/doi/abs/10.1002/nme.865},
	volume =       59,
	year =         2004
}

@book{WriggersZavarise2004,
	author =       {P. Wriggers and G. Zavarise},
	publisher =    {John Wiley \& Sons, Ltd. Edited by E. Stein, R. de
	Borst and T.J.R. Hughes},
	title =        {Encyclopedia of Computational Mechanics, Volume 2:
	Solids and Structures},
	year =         2004
}

@article{Wriggers1985,
	author =       {P. Wriggers and J. Simo and R. Taylor},
	journal =      {Proceedings of the NUMETA '85 Conference, Elsevier,
	Amsterdam},
	title =        {Penalty and augmented {L}agrangian formulations for
	contact problems},
	year =         1985
}

@article{SIMO1992,
	title = {An augmented lagrangian treatment of contact problems involving friction},
	journal = {Computers \& Structures},
	volume = {42},
	number = {1},
	pages = {97-116},
	year = {1992},
	issn = {0045-7949},
	doi = {https://doi.org/10.1016/0045-7949(92)90540-G},
	url = {https://www.sciencedirect.com/science/article/pii/004579499290540G},
	author = {J.C. Simo and T.A. Laursen}
}

@article{Popp2012,
	author = {Popp, A. and Wohlmuth, B. and Gee, M. and Wall, W.},
	year = {2012},
	month = {07},
	pages = {B421-B446},
	title = {Dual Quadratic Mortar Finite Element Methods for 3D Finite Deformation Contact},
	volume = {34},
	journal = {SIAM Journal on Scientific Computing},
	doi = {10.1137/110848190}
}

@article{Wriggers2008,
	author =       {P. Wriggers and G. Zavarise},
	journal =      {Computational Mechanics},
	pages =        {407--420},
	title =        {A formulation for frictionless contact problems
	using a weak form introduced by {N}itsche},
	volume =       41,
	year =         2008,
}

@article{Kikuchi1981,
	author = {Kikuchi, N. and Song, Y.},
	year = {1981},
	month = {04},
	pages = {},
	title = {Penalty/finite-element approximation of a class of unilateral problems in linear elasticity},
	volume = {39},
	journal = {Quarterly of Applied Mathematics},
	doi = {10.1090/qam/613950}
}

@book{Kikuchi1988,
	publisher = {SIAM},
	booktitle={Studies in Applied and Numerical Mathematics},
	year = {1988},
	title = {Contact problems in elasticity: a study of variational inequalities and finite element methods},
	author = {N. Kikuchi and J. T. Oden},
}

@book{Wriggers2006,
	author = {Wriggers, P.},
	year = {2006},
	month = {01},
	publisher = {Springer},
	title = {Computational Contact Mechanics},
	isbn = {978-3-540-32608-3},
	doi = {10.1007/978-3-540-32609-0}
}

@book{Yastrebov2013,
	author = {Yastrebov, V.},
	year = {2013},
	month = {02},
	publisher = {Wiley},
	title = {Numerical Methods in Contact Mechanics},
	isbn = {9781848215191},
	doi = {10.1002/9781118647974}
}

@article{DIONE2019,
	title = {Optimal convergence analysis of the unilateral contact problem with and without Tresca friction conditions by the penalty method},
	journal = {Journal of Mathematical Analysis and Applications},
	volume = {472},
	number = {1},
	pages = {266-284},
	year = {2019},
	issn = {0022-247X},
	doi = {https://doi.org/10.1016/j.jmaa.2018.11.023},
	url = {https://www.sciencedirect.com/science/article/pii/S0022247X18309612},
	author = {I. Dione}
}

@article{NESTEROV1983,
	author = {Nesterov, Y. E.},
	title = {A method of solving a convex programming problem with convergence rate {$O\!\left(1 / k^{2} \right)$}},
	journal = {Dokl. Akad. Nauk SSSR},
	volume = {269},
	number = {3},
	pages = {543-547},
	year = {1983},
	doi = {},
	URL = {},
}

@book{CIARLET1989, 
	place={Cambridge}, 
	series={Cambridge Texts in Applied Mathematics}, 
	title={Introduction to Numerical Linear Algebra and Optimisation}, publisher={Cambridge University Press}, 
	author={Ciarlet, P. G.}, 
	year={1989}, 
	collection={Cambridge Texts in Applied Mathematics}
}

@book{UZAWA1958,
	author    = {K. J. Arrow and L. Hurwicz and H. Uzawa},
	title     = {Studies in Linear and Non-Linear Programming},
	year      = {1958},
	publisher = {Stanford University Press},
	address   = {Redwood City}
}

@article{Odonoghue2015,
	author = {B. O'Donoghue and E. Cand\`{e}s},
	title = {{Adaptive Restart for Accelerated Gradient Schemes}},
	year = {2015},
	issue_date = {June 2015},
	publisher = {Springer-Verlag},
	address = {Berlin, Heidelberg},
	volume = {15},
	number = {3},
	issn = {1615-3375},
	url = {https://doi.org/10.1007/s10208-013-9150-3},
	doi = {10.1007/s10208-013-9150-3},
	journal = {Found. Comput. Math.},
	month = jun,
	pages = {715–732},
	numpages = {18}
}

@article{Anderson1965,
	author = {Anderson, D. G.},
	title = {Iterative Procedures for Nonlinear Integral Equations},
	year = {1965},
	issue_date = {Oct. 1965},
	publisher = {Association for Computing Machinery},
	address = {New York, NY, USA},
	volume = {12},
	number = {4},
	issn = {0004-5411},
	url = {https://doi.org/10.1145/321296.321305},
	doi = {10.1145/321296.321305},
	journal = {Journal of the ACM},
	month = oct,
	pages = {547–560},
	numpages = {14}
}

@article{Benzi2005, 
	title={Numerical solution of saddle point problems}, 
	volume={14}, 
	DOI={10.1017/S0962492904000212}, 
	journal={Acta Numerica}, 
	author={Benzi, M. and Golub, G. H. and Liesen, J.}, year={2005}, 
	pages={1–137}
}

@article{NOUROMID1986,
	title = {A two-level iteration method for solution of contact problems},
	journal = {Computer Methods in Applied Mechanics and Engineering},
	volume = {54},
	number = {2},
	pages = {131-144},
	year = {1986},
	issn = {0045-7825},
	doi = {https://doi.org/10.1016/0045-7825(86)90122-2},
	url = {https://www.sciencedirect.com/science/article/pii/0045782586901222},
	author = {B. Nour-Omid and P. Wriggers}
}

@article{Siefert2006,
	author = {Siefert, M. and de Sturler, E.},
	year = {2006},
	month = {07},
	pages = {1275-1296},
	title = {Preconditioners for Generalized Saddle-Point Problems},
	volume = {44},
	journal = {SIAM Journal on Numerical Analysis},
	doi = {10.1137/040610908}
}

@book{Luenberger2008,
	author = {Luenberger, D. and Ye, Y.},
	year = {2008},
	title = {Linear and Nonlinear Programming},
	publisher = {Springer},
}

@article{Acary2024,
	author = {V. Acary and P. Armand and H. M. Nguyen and M. Shpakovych},
	title = {Second-order cone programming for frictional contact mechanics using interior point algorithm},
	journal = {Optimization Methods and Software},
	volume = {39},
	number = {3},
	pages = {634--663},
	year = {2024},
	publisher = {Taylor \& Francis},
	doi = {10.1080/10556788.2023.2296438},
	URL = {https://doi.org/10.1080/10556788.2023.2296438	},
	eprint = { 	https://doi.org/10.1080/10556788.2023.2296438}
}

@article{TEMIZER2014,
	title = {An interior point method for isogeometric contact},
	journal = {Computer Methods in Applied Mechanics and Engineering},
	volume = {276},
	pages = {589-611},
	year = {2014},
	issn = {0045-7825},
	doi = {https://doi.org/10.1016/j.cma.2014.03.018},
	url = {https://www.sciencedirect.com/science/article/pii/S0045782514001042},
	author = {{\.I}. Temizer and M.M. Abdalla and Z. G{\"u}rdal}
}

@article{KOTHARI2022,
	title = {{A generalized multigrid method for solving contact problems in Lagrange multiplier based unfitted Finite Element method}},
	journal = {Computer Methods in Applied Mechanics and Engineering},
	volume = {392},
	pages = {114630},
	year = {2022},
	issn = {0045-7825},
	doi = {https://doi.org/10.1016/j.cma.2022.114630},
	url = {https://www.sciencedirect.com/science/article/pii/S004578252200041X},
	author = {H. Kothari and R. Krause}
}

@article{INTROINI2024,
	title = {{ALCYONE}: the fuel performance code of the {PLEIADES} platform dedicated to PWR fuel rods behavior},
	journal = {Annals of Nuclear Energy},
	volume = {207},
	pages = {110711},
	year = {2024},
	issn = {0306-4549},
	doi = {https://doi.org/10.1016/j.anucene.2024.110711},
	url = {https://www.sciencedirect.com/science/article/pii/S0306454924003748},
	author = {C. Introïni and I. Ramière and J. Sercombe and B. Michel and T. Helfer and J. Fauque}
}

@article{MICHEL2008,
	title = {{3D} fuel cracking modelling in pellet cladding mechanical interaction},
	journal = {Engineering Fracture Mechanics},
	volume = {75},
	number = {11},
	pages = {3581-3598},
	year = {2008},
	issn = {0013-7944},
	doi = {10.1016/j.engfracmech.2006.12.014},
	url = {https://www.sciencedirect.com/science/article/pii/S0013794406004759},
	author = {B. Michel and J. Sercombe and G. Thouvenin and R. Chatelet}
}

@incollection{sercombe2020,
	TITLE = {Modeling of pellet cladding interaction},
	AUTHOR = {J. Sercombe and B. Michel and C. Riglet-Martial and O. Fandeur},
	URL = {https://hal.science/hal-04825099},
	BOOKTITLE = {{Comprehensive Nuclear Materials - Second Edition. Volume 2 : Oxide Fuel Systems in Thermal and Fast Neutron Spectrum Reactors}},
	EDITOR = {Rudy J.M. Konings and Roger E. Stoller},
	PUBLISHER = {{Elsevier}},
	VOLUME = {2},
	YEAR = {2020},
	DOI = {10.1016/B978-0-12-803581-8.00715-3},
	HAL_ID = {hal-04825099},
	HAL_VERSION = {v1},
}

@article{Abide16,
  title={Analysis of two active set type methods to solve unilateral contact problems},
  author={Abide, S. and Barboteu, M. and Danan, D.},
  journal={Applied Mathematics and Computation},
  volume={284},
  pages={286--307},
  year={2016},
  publisher={Elsevier}
}

@Article{kuttler-wall2008,
  author = 	 {U. K\"uttler and W. A. Wall},
  title = 	 {Fixed-point fluid-structure interaction solvers with dynamic relaxation},
  journal = 	 {Computational Mechanics},
  year = 	 {2008},
  OPTkey = 	 {},
  volume = 	 {43},
  OPTnumber = 	 {},
  pages = 	 {61-72},
  OPTmonth = 	 {},
  OPTnote = 	 {},
  OPTannote = 	 {}
}

@Article{erbts-duester2012,
  author = 	 {P. Erbts and A. D\"uster},
  title = 	 {Accelerated staggered coupling schemes for problems of thermoelasticity at finite strains},
  journal = 	 {Computers and Mathematics with Applications},
  year = 	 {2012},
  OPTkey = 	 {},
  volume = 	 {64},
  OPTnumber = 	 {},
  pages = 	 {2408-2430},
  OPTmonth = 	 {},
  OPTnote = 	 {},
  OPTannote = 	 {}
}

@Article{Ruiz_2022,
  author = 	 {C. Ruiz and N. Lebaz and I. Ramière and D. Mangin and M. Bertrand},
  title = 	 {Fixed point convergence and acceleration
    for steady state population balance modelling of precipitation
    processes: Application to neodymium oxalate},
  journal = 	 {Chemical Engineering Research and Design},
  year = 	 {2022},
  OPTkey = 	 {},
  volume = 	 {177},
  pages = 	 {767-777},
  doi= {10.1016/j.cherd.2021.11.030},
}

@article{Nour_1987,
author = {Nour-Omid, B. and Wriggers, P.},
title = {A note on the optimum choice for penalty parameters},
journal = {Communications in Applied Numerical Methods},
volume = {3},
number = {6},
pages = {581-585},
doi = {https://doi.org/10.1002/cnm.1630030620},
url = {https://onlinelibrary.wiley.com/doi/abs/10.1002/cnm.1630030620},
eprint = {https://onlinelibrary.wiley.com/doi/pdf/10.1002/cnm.1630030620},
year = {1987}
}

@Article{Zulehner_2001,
  author = 	 {W. Zulehner},
  title = 	 {Analysis of iterative methods for saddle point problems: a unified approach},
  journal = 	 {Mathematics of Computation},
  year = 	 {2001},
  OPTkey = 	 {},
  volume = 	 {71},
  number = 	 {238},
  pages = 	 {479-505},
  doi  =         {10.1090/S0025-5718-01-01324-2},
  OPTmonth = 	 {},
  OPTnote = 	 {},
  OPTannote = 	 {}
}

@article{Ouyang_AA_2024,
author = {Ouyang, W. and Liu, Y. and Milzarek, A.},
title = {Descent Properties of an {A}nderson Accelerated Gradient Method with Restarting},
journal = {SIAM Journal on Optimization},
volume = 34,
number = 1,
pages = {336-365},
year = 2024,
doi = {10.1137/22M151460X}
}

@article{Krzysik_AA_2025,
author = {Krzysik, O. A. and De Sterck, H. and Smith, A.},
title = {Asymptotic Convergence of Restarted {A}nderson Acceleration for Certain Normal Linear Systems},
journal = {SIAM Journal on Scientific Computing},
OPTvolume = {0},
OPTnumber = {0},
pages = {S135-S160},
year = {2025},
doi = {10.1137/24M1672262},
    
}

@article{Barzilai_Borwein_1988,
    author = {Barzilai, J. and Borwein, J. M.},
    title = {Two-Point Step Size Gradient Methods},
    journal = {IMA Journal of Numerical Analysis},
    volume = {8},
    number = {1},
    pages = {141-148},
    year = {1988},
    month = {01},
    issn = {0272-4979},
    doi = {10.1093/imanum/8.1.141},
    url = {https://doi.org/10.1093/imanum/8.1.141},
    eprint = {https://academic.oup.com/imajna/article-pdf/8/1/141/2402762/8-1-141.pdf},
}

@article{Li_Accel_2020,
title = {Accelerated fixed-point formulation of topology optimization: Application to compliance minimization problems},
journal = {Mechanics Research Communications},
volume = {103},
pages = {103469},
year = {2020},
issn = {0093-6413},
doi = {https://doi.org/10.1016/j.mechrescom.2019.103469},
url = {https://www.sciencedirect.com/science/article/pii/S0093641319305051},
author = {W. Li and P. Suryanarayana and G. H. Paulino}
}

@misc{CEA_CAST3M_2023,
	author       = {{CEA}},
	title        = {{CAST3M}},
	year         = {2023},
	howpublished = {\url{https://www.cast3m.cea.fr/}},
}

@article{LRL_CMAME_2017,
author = {H. Liu and  I. Rami\`ere and F. Lebon},
title = {On the coupling of local multilevel mesh refinement and {ZZ} methods for unilateral frictional contact problems in elastostatics},
journal = {Computer Methods in Applied Mechanics and Engineering},
volume = {323},
pages = {1--26},
year = {2017},
issn = {0045-7825},
doi = {https://doi.org/10.1016/j.cma.2017.04.011},
url = {http://www.sciencedirect.com/science/article/pii/S0045782517301287},
}

@article{KRL_CMAME_2022,
title = {Fully automatic multigrid adaptive mesh refinement strategy with controlled accuracy for nonlinear quasi-static problems},
journal = {Computer Methods in Applied Mechanics and Engineering},
volume = {400},
pages = {115505},
year = {2022},
issn = {0045-7825},
doi = {10.1016/j.cma.2022.115505},
url = {https://www.sciencedirect.com/science/article/pii/S0045782522005151},
author = {Daria Koliesnikova and Isabelle Ramière and Frédéric Lebon},
}


\end{document}